\documentclass[notitlepage]{article}
\usepackage{authblk}
\usepackage{graphicx,url,color}
\usepackage{epstopdf}
\usepackage{amsfonts,bm}
\usepackage{graphicx}
\usepackage{amsmath}
\usepackage{setspace}
\usepackage{float}
\usepackage{subfigure,color,epsfig,multirow,bbm,graphicx,graphics, sidecap}
\usepackage{tabularx,ulem}
\usepackage[sort&compress,numbers]{natbib}
\usepackage{algorithm}
\usepackage{algpseudocode,ulem}
\usepackage[margin=1.1in]{geometry}
\providecommand{\keywords}[1]{\textbf{\textit{Keywords---}} #1}

\begin{document}

\thispagestyle{empty}
\setcounter{page}{0}

\begin{Huge}
\begin{center}
Computer Science Technical Report CSTR-{\tt TR 2/2014} \\
\today
\end{center}
\end{Huge}
\vfil
\begin{huge}
\begin{center}
{\tt R\u azvan \c Stef\u anescu, Adrian Sandu, and \\Ionel M. Navon}
\end{center}
\end{huge}

\vfil
\begin{huge}
\begin{it}
\begin{center}
``{\tt Comparison of POD reduced order strategies for the nonlinear 2D Shallow Water Equations}''
\end{center}
\end{it}
\end{huge}
\vfil

\begin{large}
\begin{center}
Computational Science Laboratory \\
Computer Science Department \\
Virginia Polytechnic Institute and State University \\
Blacksburg, VA 24060 \\
Phone: (540)-231-2193 \\
Fax: (540)-231-6075 \\
Email: \url{sandu@cs.vt.edu} \\
Web: \url{http://csl.cs.vt.edu}
\end{center}
\end{large}

\vspace*{1cm}

\begin{tabular}{ccc}
\includegraphics[width=2.5in]{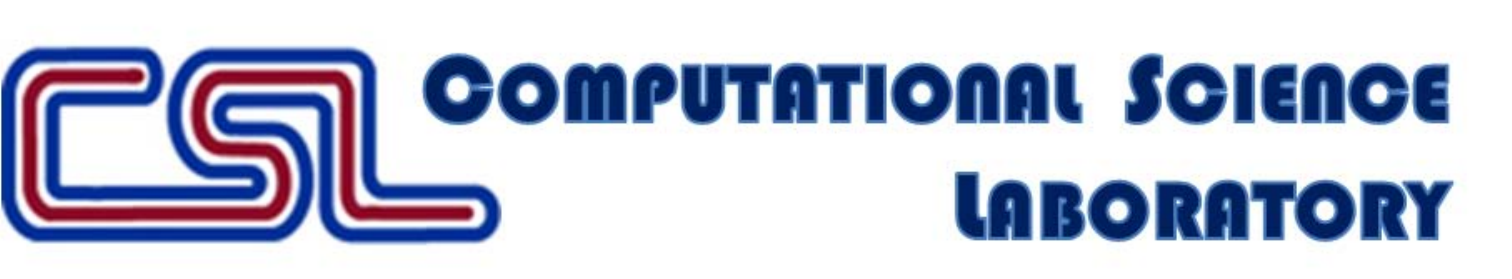}
&\hspace{2.5in}&
\includegraphics[width=2.5in]{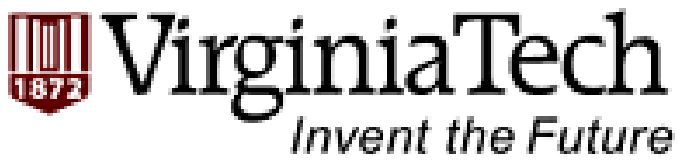} \\
{\bf\em Innovative Computational Solutions} &&\\
\end{tabular}

\newpage

\title{Comparison of POD reduced order strategies for the nonlinear 2D Shallow Water Equations}
\author[1]{R\u azvan \c Stef\u anescu \thanks{rstefane@vt.edu}}
\author[1]{Adrian Sandu \thanks{sandu@cs.vt.edu}}
\author[2]{Ionel M. Navon \thanks{inavon@fsu.edu}}
\affil[1]{Computational Science Laboratory, Department of
Computer Science, Virginia Polytechnic Institute and State University, Blacksburg, Virginia, USA, 24060}
\affil[2]{Department of Scientific Computing, The Florida State University, Tallahassee, Florida, USA, 32306}
\date{}
\maketitle
\begin{abstract}

This paper introduces tensorial calculus techniques in the framework of Proper Orthogonal Decomposition (POD) to reduce the computational complexity of the reduced nonlinear terms. The resulting method, named tensorial POD, can be applied to polynomial nonlinearities of any degree $p$. Such nonlinear terms have an on-line complexity of $\mathcal{O}(k^{p+1})$, where $k$ is the dimension of POD basis, and therefore is independent of full space dimension. However it is efficient only for quadratic nonlinear terms since for higher nonlinearities standard POD proves to be less time consuming once the POD basis dimension $k$ is increased. Numerical experiments are carried out with a two dimensional shallow water equation (SWE) test problem to compare the performance of tensorial POD, standard POD, and POD/Discrete Empirical Interpolation Method (DEIM). Numerical results show that tensorial POD decreases by $76\times$ times the computational cost of the on-line stage of standard POD for configurations using more than $300,000$ model variables. The tensorial POD SWE model was only $2-8\times$ slower than the POD/DEIM SWE model but the implementation effort is considerably increased. Tensorial calculus was again employed to construct a new algorithm allowing POD/DEIM shallow water equation model to compute its off-line stage faster than the standard and tensorial POD approaches.
\end{abstract}

\section{Introduction}

Modeling and simulation of multi-scale complex physical phenomena leads to large-scale systems of
coupled partial differential equations, ordinary differential equations,  and differential algebraic equations.  The high dimensionality of these models poses important mathematical and computational challenges. A computationally feasible approach to simulate, control, and optimize such systems is to simplify the models
by retaining only those state variables that are consistent with a particular phenomena of interest.

Reduced order modeling refers to the development of low-dimensional models that
represent important characteristics of a high-dimensional or infinite dimensional
dynamical system. The reduced order methods can be cast into three broad categories: Singular Values Decomposition (SVD) based methods, Krylov based methods and iterative methods combining aspects of both the SVD and Krylov methods (see e.g.  \citet{antoulas2005als}).

For linear models, methods like balanced truncation (\citet{Moore1981,Antoulas2009,Sorensen_Antoulas2002, Mullis_Roberts1976}) and moment matching (\citet{Freund2003,Feldmann_Freund1995,Grimme1997})  have been proving successful in developing reduced order models. However, balanced truncation doesn't extend easily for high-order systems, and several
grammians approximations were proposed leading to methods such as approximate subspace iteration (\citet{Bakeretal1996}),
least squares approximation (\citet{Hodel1991}), Krylov subspace methods (\citet{Jaimoukha_Kasenally1994} and
\citet{Gudmundsson_Laub1994} ) and balanced Proper Orthogonal Decomposition (\citet{Willcox02balancedmodel}).
Among moment matching  methods we mention partial realization (\citet{Gragg_Lindquist1983,Benner_Sokolov2006}),
Pad\'{e} approximation (\citet{Gragg1972,Gallivan94,Gutknecht1994,VanDooren_1995})
and rational approximation (\citet{Bultheel_Moor2000}).

While for linear models we are able to produce input-independent highly accurate reduced models,
in the case of general nonlinear systems, the transfer function approach is not yet applicable and input-specified semi-empirical methods are usually employed. Recently some encouraging research results using  generalized transfer functions and generalized moment matching have been obtained by \citet{MPIMD12-12} for nonlinear model order reduction but future investigations are required.

Proper Orthogonal Decomposition and its variants are also known as Karhunen-Lo\`{e}ve expansions
\cite{karhunen1946zss,loeve1955pt}, principal component analysis \citet{hotelling1939acs}, and empirical orthogonal functions \citet{lorenz1956eof} among others. It is the most prevalent basis selection method for nonlinear problems
and, among other requirements, relies on the fact that the desired simulation is well simulated in the input collection. Data analysis using POD is conducted to
extract basis functions, from experimental data or detailed simulations of high-dimensional systems (method of snapshots introduced by \citet{Sir87a, Sir87b, Sir87c}), for subsequent use in Galerkin projections that yield low dimensional dynamical models. Unfortunately the standard POD approach displays a major disadvantage since its nonlinear reduced terms still have to be evaluated on the original state space making the simulation of the reduced-order system too expensive. There exist several ways to avoid this problem such as the empirical interpolation method (EIM) \citet{BMN2004} and its discrete variant DEIM \citet{Cha2008,ChaSor2012,ChaSor2010}, best points interpolation method \citet{NPP2008}. Missing point estimation \citet{Astrid_2008} and Gauss-Newton with approximated tensors
\citet{Amsallem_2011,Carlberg_2011,Carlberg2_2011,Carlberg_2012, Carlberg_2012} methods are relying upon the gappy POD technique \citet{Everson_1995} and were developed for
the same reason. Reduced basis methods have been recently developed and utilize on greedy algorithms to efficiently compute numerical solutions for parametrized applications \citet{BMN2004,grepl2005posteriori,patera2007reduced,rozza2008reduced,Dihlmann_2013}.

Dynamic mode decomposition is a relatively recent development in the field of modal decomposition (\citet{Rowley2009,Schmid2010,Tissot2013}) and in comparison with POD approximates the temporal dynamics by a high-degree polynomial.
Trajectory piecewise linear method proposed by \citet{Rewienski2001} follows a different strategy where the nonlinear system
is represented by a piecewise-linear system which then can be efficiently approached by the standard linear reduction method.
Parameter model reduction has emerged recently as an important research direction and \citet{MPIMD13-14} highlights the major contribution in the field.

This paper combines  standard POD and tensor calculus techniques to reduce the on-line computational complexity of the reduced nonlinear terms for a shallow water equations model. Tensor based calculus was already applied
by \citet{Kunisch_Volkwein_1999,Kunisch_Volkwein_2004} to represent quadratic nonlinearities of reduced order POD models.
We show that the tensorial POD (TPOD) approach can be applied to polynomial nonlinearities of any degree $p$,
and the its representation has a complexity of $\mathcal{O}(k^{p+1})$, where $k$ is the dimension of POD subspace. This complexity is independent of the full space dimension.
For $k$ between $10$ and $50$ and $p=2$ the number of floating-point operations required to calculate the tensorial POD quadratic terms is $10$--$40\times$ lower than in the case of standard POD, and $10$--$20\times$ higher than for the POD/DEIM. However, CPU time for solving the TPOD SWE model (on-line stage) is only $2$--$8\times$ times slower than POD/DEIM SWE model for $10^3$--$10^5$ grid points, $k\leq 50$, and number of DEIM interpolation points $m\leq 180$. For example, for an integration interval of $3$h, $10^5$ mesh points, $k=50$, and $m=70$, tensorial POD and POD/DEIM are $76\times$ and $450\times$ faster than standard POD, but the implementation effort of POD/DEIM is considerably increased. Many useful models are characterized by quadratic nonlinearities in both fluid dynamics and geophysical fluid flows including SWE model. In the case of cubic or higher polynomial nonlinearities the advantage of tensorial POD is lost and its nonlinear computational complexity is similar or larger than the computational complexity of the standard POD approach. This proves that for models depending only on quadratic nonlinearities, the tensorial POD represents a solid alternative to POD/DEIM where the implementation effort is considerably larger. We also propose a fast algorithm to pre-compute the reduced order coefficients for polynomial nonlinearities of order $p$ which allows the POD/DEIM SWE model to to compute its off-line stage faster than the standard and tensorial
POD approaches despite additional SVD calculations and reduced coefficients computations.

The paper is organized as follows. Section \ref{sec:ROM} reviews the reduced order modeling methodologies used in this work:  standard, tensorial, and DEIM POD. Section \ref{sec:complexity} analyses the computational complexity of the reduced order polynomial nonlinearities for all three methods, and introduces a new DEIM based algorithm to efficiently compute the coefficients needed for reduced Jacobians. Section \ref{sec:SWE} discusses the shallow water equations model and its full implementation, and Section \ref{sec:ROMS_SWE} describes the construction of reduced models.
Results of extensive numerical experiments are discussed in Section \ref{sec:Numerical_resuls} while conclusions are drawn in Section \ref{sec:Conclusions}.

\section{Reduced Order Modeling}\label{sec:ROM}

For highly efficient flows simulations, reduced order modeling is a powerful tool for representing the dynamics of large-scale dynamical systems using only a smaller number of variables and reduced order basis functions.
Three approaches will be considered in this study: standard Proper Orthogonal Decomposition (POD), tensorial POD (TPOD), and POD/Discrete Empirical Interpolation Method (POD/DEIM). They are discussed below.
The tensorial POD approach proposed herein is different than the method of \citet{Belzen2012} which makes use of tensor decompositions for generating  POD bases.
\subsection{Standard Proper Orthogonal Decomposition}
Proper Orthogonal Decompositions has been used successfully in numerous applications such as compressible flow \citet{Rowley2004}, computational
fluid dynamics \citet{Kunisch_Volkwein_POD2002,Rowley2005,Willcox02balancedmodel}, aerodynamics  \cite{Bui-thanh04aerodynamicdata}. It can be thought of as a Galerkin approximation in
the spatial variable built from functions corresponding to the solution of the physical system at specified time instances. \citet{Noack2010}
proposed a system reduction strategy for Galerkin models of fluid flows leading to dynamic models of lower order based on a partition in slow, dominant and fast modes. \citet{San_Iliescu2013} investigate several closure models for POD reduced order modeling of fluids flows and benchmarked against the fine resolution numerical simulation.

In what follows, we will only work with discrete inner products (Euclidian dot product) though continuous products may be employed too.  Generally, an atmospheric or oceanic model is
usually governed by the following semi--discrete dynamical system
\begin{equation}\label{eqn::-4}
 \frac{d{\bf x}(t)}{dt} = {\bf F}({\bf x},t),~~~~{\bf x}(0) = {\bf x}_0 \in \mathbb{R}^n.
\end{equation}
From the temporal-spatial flow ${\bf x}(t) \in \mathbb{R}^n $, we select an ensemble of $ N_t $ time instances ${\bf x}_1,...,{\bf x}_{N_t}  \in \mathbb{R}^n,~n$ being the total number of discrete model variables per time step  and $N_t \in \mathbb{N},~N_t>0$. Let us define
the centering trajectory, shift mode, or mean field correction (\citet{NAMTT03}) ${\bf \bar x} =\frac{1} N_t \sum_{i=1}^{N_t} {\bf x}_i$. The method of POD consists in choosing a
complete orthonormal basis $U=\{{\bf u}_{i}\},~i=1,..,k;~k>0;~u_i \in \mathbb{R}^n;~U \in \mathbb{R}^{n\times k}$ such that the mean square error between ${\bf x}(t)$ and
POD expansion ${\bf x}^{POD}(t) = {\bf \bar x} + U{\bf \tilde x}(t),~{\bf \tilde x}(t) \in \mathbb{R}^k$ is minimized on average. The POD dimension $k \ll n$ is appropriately chosen to capture the dynamics of the flow as follows:

\begin{algorithm}
 \begin{algorithmic}[1]
 \State Calculate the mean ${\bf \bar x} =\frac{1}{N_t} \sum_{i=1}^{N_t} {\bf x}_i$.
 \State Set up the correlation matrix $K=[k_{ij}]_{i,j=1,..,n}$ where $k_{ij} = \langle{\bf x}_i - {\bf \bar x},{\bf x}_j - {\bf \bar x}\rangle$, and $\langle\cdot, \cdot\rangle$ being the Euclidian dot product.
 \State Compute the eigenvalues $\lambda_1\geq \lambda_2\geq ...\lambda_n\geq 0$ and the corresponding orthogonal eigenvectors ${\bf v}^1,{\bf v}^2,..,{\bf v}^n \in \mathbb{R}^n$ of $K$.
 \State Set ${\bf u}_{i} = \langle{\bf v}^i,{\bf x}_i - {\bf \bar x}\rangle$, $i=1,..,n.$ Then, ${\bf u}_{i},~i=1,..,n$ are normalized to obtain an orthonormal basis.
\State Define $I(m)={\frac {\sum_{i=1}^m \lambda_i}{\sum_{i=1}^n \lambda_i}}$ and choose $k$ such that $ k=\min \{I(m):I(m)\geq \gamma\}$ where $0 \leq \gamma \leq 1$ is the percentage of
total informations captured by the reduced space $\textrm{span}\{{\bf u}_1,{\bf u}_2,...,{\bf u}_k\}.$ Usually $\gamma$ is taken $0.99$.

 \end{algorithmic}
 \caption{POD basis construction}
 \label{euclid}
\end{algorithm}

To obtain the reduced model of (\ref{eqn::-4}), we first employ a numerical scheme to solve the full model for a set of snapshots and follow the above procedure, then use a Petrov--Galerkin (PG) projection
of the full model equations onto the space $\mathcal{X}^k$ spanned by the POD basis elements
\begin{equation}\label{eqn::-3}
 \frac{d{\bf \tilde x}(t)}{dt} = W^T{\bf F}\bigg({\bf \bar x}+U{\bf \tilde x}(t),t\bigg),~~~{\bf \tilde x}(0) = W^T\bigg({\bf x}(0)-{\bf \bar x}\bigg),
\end{equation}
where $W \in \mathbb{R}^{n \times k}$ contains the discrete test functions from the PG projection, i.e. $W^TU = I \in \mathbb{R}^k$. The Galerkin projection may be also a choice being just a particular case of PG ($W = U$).

The efficiency of the POD-Galerkin techniques is limited to linear or bilinear terms, since the projected
nonlinear term at every discrete time step still depends on the number of variables of the full model:
$$\tilde N ({\bf \tilde x}) =\underbrace{W^T}_{k \times n} \underbrace{{\bf F}({\bar {\bf x}}+U{\bf \tilde x}(t))}_{n \times 1}. $$

To be precise, consider a steady polynomial nonlinearity $x^p$. A POD expansion involving mean ${\bf \bar x}$ will unnecessarily complicate the description of tensorial POD representation of a $p^{\rm th}$ order polynomial nonlinearity. Moreover the terms depending on ${\bf \bar x}$ are just a particular case of the term depending only on $U{\bf \tilde x}$ since vector componentwise multiplication is distributive over vector addition. Consequently the expansion ${\bf \bar x}\approx U{\bf \tilde x}$ will not decrease the generality of the reduced nonlinear term.
In the finite difference case, the standard POD projection is
described as follows
\begin{equation}\label{POD_nonlinearity}
\tilde N({\bf \tilde x}) = \underbrace{W^T}_{k \times n} \,  \underbrace{\left(U{\bf \tilde x} \right)^p}_{n \times 1}
\end{equation}
where vector powers are taken component-wise.

To mitigate this inefficiency we propose two approaches: {\bf (1)} Tensorial POD and {\bf (2)} Discrete Empirical Interpolation Method. The former approach is able to calculate the reduced polynomial nonlinearities independent of $n$, while the latter method can handle efficiently all type of nonlinearities.
\subsection{Tensorial POD}
Tensorial POD technique employs the simple structure of the polynomial nonlinearities to remove the dependence
on the dimension of the original discretized system by manipulating the order of computing. It can be successfully used in a POD framework for finite
difference (FD), finite element (FE) and finite volume (FV) discretization methods and all other type of discretization methods that engage in spectral expansions. Tensorial POD separates the full spatial variables from reduced variables allowing fast nonlinear terms computations in the on--line stage. For time dependent nonlinearities this implies separation of spatial variables from reduced time variables. Thus, the reduced nonlinear term evaluation requires a tensorial Frobenius dot-product computation between rank $p$ tensors, where $p$ is the order of polynomial nonlinearity. The projected spatial variables are stored into tensors and calculated off--line. These are also used for reduced Jacobian computation in the on-line stage.

The tensorial POD representation of (\ref{POD_nonlinearity}) is given by vector
\begin{equation}
 \label{tensor_POD_nonlinearity}
 {\bf \mathcal{M}} = \big[{\bf \mathcal{M}}^i\big]_{i=1,2,..,k};
 \quad {\bf \mathcal{M}}^i = \left\langle {\bf M}^i,{\bf \tilde X}\right\rangle_{\textrm{Frobenius}}\in \mathbb{R},~i=1,2,..,k;
 \quad {\bf \mathcal{M}}\in \mathbb{R}^k,
\end{equation}
where $p-$order tensors ${\bf \tilde X}$ and ${\bf M}^i,~i=1,2,..,k$  are defined as
\begin{equation}
\begin{split}
& {\bf \tilde X} = \big[{\bf \tilde X}_{i_1i_2..i_p}\big]_{i_1,i_2,..,i_p=1,..,k}\in \mathbb{R}^{\underbrace{k\times  ... \times k}_{\textrm{p times}}};
\quad {\bf \tilde X}_{i_1i_2..i_p}={\bf \tilde x}_{i_1}{\bf \tilde x}_{i_2}...{\bf \tilde x}_{i_p} \in \mathbb{R}; \\
& {\bf M}^i = \big[{{\bf M}^i}_{i_1i_2..i_p}\big]_{i_1,i_2,..,i_p=1,2,..k},\quad i=1,2,..,k;~{\bf M}^i\in \mathbb{R}^{\underbrace{k\times  ... \times k}_{\textrm{p times}}};\\
 & {{\bf M}^i}_{i_1i_2..i_p} = \sum_{l=1}^nW_{li}U_{li_1}U_{li_2}...U_{li_p} \in \mathbb{R},
\end{split}
\label{tensor_POD}
\end{equation}
and ${\bf \tilde x}_{i_j},~W_{li},~U_{li_j}$ are just entries of POD reduced order solution ${\bf \tilde x}$, POD test functions basis $W$
and POD trial functions basis $U$. The tensorial Frobenius dot product is defined as
\begin{eqnarray*}
&& \langle\cdot,\cdot\rangle_{\textrm{Frobenius}} :
\mathbb{R}^{\underbrace{k\times ... \times k}_{\textrm{p times}}} \times \mathbb{R}^{\underbrace{k\times  ... \times k}_{\textrm{p times}}} \rightarrow \mathbb{R},\\
&& \langle{\bf A},{\bf B}\rangle_{\textrm{Frobenius}} = {\bf A:B}
= \sum_{i_1,i_2,..,i_p = 1}^k
A_{i_1i_2..i_p}B_{i_1i_2,..i_p}\in \mathbb{R}.
\end{eqnarray*}

We note that ${\bf M}^i,~i=1,2,..,k$ are $p^{th}$ order tensors computed in the off-line stage and their dimensions do not depend on the full space dimension. For finite element and finite volume the tensorial  POD representations are the same except for the type of products used in computation of ${\bf M}^i_{i_1i_2..i_p}$ in (\ref{tensor_POD}) which now are continuous and replace the
sum of products used in the finite difference case.

Reduced nonlinearities depending on space derivatives are treated similarly as in equation (\ref{POD_nonlinearity} - \ref{tensor_POD})
since POD expansion of $x_x$ (space derivative of $x$) is $U_x\tilde x$, where $U_x \in \mathbb{R}^{n \times k}$ contains the space derivatives of POD basis
functions ${\bf u}_i,~i=1,2,..,k$, ${\bf u}_i\in\mathbb{R}^n$.
\subsection{Standard POD and Discrete Empirical Interpolation Method}
The empirical interpolation method (EIM) and its discrete version (DEIM), were developed to approximate the nonlinear term allowing an effectively affine offline–-online computational decomposition. Both interpolation methods provide an efficient way to approximate nonlinear functions.
They were successfully used in a standard POD framework for finite difference (FD),  finite element (FE) and  finite volume (FV) discretization methods. A description of EIM in connection with the reduced basis framework and a posteriori error bounds can be found in \citet{Mad_Nguy_Pat_Pau,Grep_Mad_Nguy_Pat}.

The DEIM implementation is based on a POD approach combined with a greedy algorithm while the EIM implementation relies on a greedy algorithm  \citet{LW2012}.

For $m\ll n$ the finite difference POD/DEIM nonlinear term approximation is
$$\tilde N({\bf \tilde x}) \approx \underbrace{W^TV(P^TV)^{-1}}_{{\rm precomputed}~k\times m} \underbrace{{\bf F}(P^T({\bar {\bf x}}+U{\bf\tilde x}))}_{m \times 1},$$
where $V\in \mathbb{R}^{n \times m}$ gathers the first $m$ POD basis modes of nonlinear function ${\bf F}$ while $P \in \mathbb{R}^{n \times m}$
is the DEIM interpolation selection matrix.

The POD/DEIM approximation of ($\ref{POD_nonlinearity}$) is
\begin{equation}
\tilde N({\bf \tilde x}) \approx \underbrace{W^TV(P^TV)^{-1}}_{{\rm precomputed}~k\times m} \underbrace{\left(P^TU{\bf \tilde x}\right)^p}_{m \times 1}
\label{POD_DEIM_nonlinearity}
\end{equation}
where vector powers are taken component-wise and $P^TV$ $\in \mathbb{R}^{m \times m}$, $P^TU$ $\in \mathbb{R}^{m\times k}$. $P^TU$ is also recommended for pre-computation in the off-line stage.

DEIM has developed in several research directions, e.g., rigorous state space error bounds \citet{ChaSor2012}, a posteriori error estimation \citet{Wirtz2012}, 1D FitzHugh-Nagumo model \citet{ChaSor2010}, 1D simulating neurons model \citet{ChaSor2_2010}, 1D nonlinear thermal model \citet{HBW2011}, 1D Burgers equation \citet{Aanonsen2009,Cha2008}, 2D nonlinear miscible viscous fingering in porous medium \cite{ChaSor2011}, oil reservoirs models \citet{EkaSuwartadi2012}, and 2D SWE model \citet{Stefanescu2013}.
We emphasize that only few POD/DEIM studies with FE and FV methods were performed, e.g., for electrical networks \citet{HK2012} and for a
2D ignition and detonation problem \citet{Nguyen_Van_Bo}. Flow simulations past a cylinder using a hybrid reduced approach combining the quadratic
expansion method and DEIM are available in \citet{Xiao2014}.

\section{Computational Complexity of the reduced $p^{th}$ order nonlinear representations. ROMs off-line stage discussion.}\label{sec:complexity}

We will focus on finite difference reduced order $p^{th}$ order polynomial nonlinearities (\ref{POD_nonlinearity},\ref{tensor_POD_nonlinearity},\ref{POD_DEIM_nonlinearity}). We begin with an important observation. For POD ROMs construction the usually approach is to store each of the state variables separately and to project every model equations to a different POD basis corresponding to a state variable whose time derivative is present. This is also the procedure we employed for this study. In this context $n$ doesn't denote the total number of state variables but only the number of state variables of the same type which most of the time is equal with the number of mesh points. Consequently from now on we refer to $n$ as the number of spatial points.

Standard POD representation is computed with a complexity of
$\mathcal{O}\big(p\times k \times n + (p-1) \times n + k \times n\big)$ and the POD/DEIM term requires $\mathcal{O}\big(p \times k \times m + (p-1)\times m + k \times m\big)$ basic operations
in the on-line stage. Tensorial POD nonlinear term has a complexity of $\mathcal{O}(k^{p+1})$. While standard POD computational complexity still depends on the full space
dimension the other twos tensorial POD and POD/DEIM don't. Table \ref{table1} describes the number of operations required to compute the projected $p^{th}$ order polynomial
nonlinearity for each of the three ROMs approaches and various values of $n,k,m,p$.

\begingroup
\begin{table}[h]
\begin{center}
\begin{tabular}[h]{|l|l|l|l|l|l|l|}
\hline
$n$&$k$&$m$&$p$&\textrm{POD} & \textrm{POD/DEIM} & Tensorial POD  \\ \hline
$10^3$& $10$& $10$& $2$ & 31,000 & 310 & 2,990 \\ \hline
$10^3$& $10$& $10$& $3$ & 42,000 & 420 & 29,990 \\ \hline
$10^3$& $10$& $10$& $4$ & 53,000 & 530 & 299,990 \\ \hline
$10^4$& $30$& $50$& $2$ & 910,000 & 4550 & 80,970 \\ \hline
$10^4$& $30$& $50$& $3$ & 1,220,000 & 6100 & 2,429,970 \\ \hline
$10^5$& $50$& $100$& $2$ & 15,100,000 & 15,100 & 374,950 \\ \hline
$10^5$& $50$& $100$& $3$ & 20,200,000 & 20,200 & 18,749,950 \\ \hline
$10^5$& $50$& $100$& $4$ & 25,300,000 & 25,300 & 937,499,950\\ \hline
\hline
\end{tabular}
\end{center}
\caption{Number of floating-point operations in
the on--line stage for different numbers of spatial points $n$, POD modes $k$, DEIM points $m$, and polynomial orders $p$. }\label{table1}
\end{table} %
\endgroup%

Clearly POD/DEIM provides the fastest nonlinear terms computations in the on-line stage. For quadratic nonlinearities,
i.e. $p=2$, and $n = 10^5$, POD/DEIM outperforms POD and POD tensorial by $10^3\times$ and $25\times$ times. But these performances are not necessarily translated into the same CPU time rates for solving
the ROMs solutions since other more time consuming calculations may be needed (reduced Jacobians computations and their LU decompositions). It was already proven in \citet{Stefanescu2013} that for a SWE model DEIM decreases the computational complexity of the standard POD by $60 \times$ for full space dimensions $n \geq 60,000$,
and leads to a CPU time reduction proportional to $n$.
CPU times and error magnitudes comparisons will be discussed in Numerical Results Section \ref{sec:Numerical_resuls}.

For cubic nonlinearities ($p=3$), the computational complexities are almost similar for both tensorial and standard POD while for higher nonlinearities (e.g. $p=4$) tensorial POD cost becomes prohibitive.

In the context of reduced optimization, the off--line stage computational complexity weights heavily in the final CPU time costs since several POD bases updates and DEIM
interpolation points recalculations are needed during the minimization process. Since the proposed schemes are implicit in time we need to compute the reduced
Jacobians as a part of a Newton type solver. For the current study we choose to calculate derivatives exactly for all three ROMs. Consequently, some reduced coefficients
such as tensors ${\bf M}^i$ defined in (\ref{tensor_POD}) must be calculated for all three reduced approaches including POD/DEIM in the off-line stage.

A simple evaluation suggests that POD/DEIM off-line stage will be slower than the corresponding tensorial POD and POD stages since more SVD computations are required in addition to particular POD/DEIM coefficients and DEIM index points calculations (see Table \ref{table2}). At a more careful examination we noticed that we can exploit the structure of POD/DEIM nonlinear term (\ref{POD_DEIM_nonlinearity}) like in the tensorial POD approach (\ref{tensor_POD_nonlinearity},\ref{tensor_POD}) and provide a fast calculation for ${\bf M}^i$.

Thus, let us denote the precomputed term and $P^TU$ in (\ref{POD_DEIM_nonlinearity}) by $E=W^TV(P^TV)^{-1}\in\mathbb{R}^{k \times m}$ and
$U^m=P^TU \in \mathbb{R}^{m \times k}$, where $m$ is the numeber of DEIM points.\newline
The $p$-tensor ${\bf M}^i$ can be computed as follows:
\begin{equation}
\begin{split}
& {\bf M}^i_{i_1i_2..i_p} = \sum_{l=1}^mE_{il}U^m_{li_1}U^m_{li_2}...U^m_{li_p} \in \mathbb{R},\quad, i,i_1,i_2,..,i_p=1,2,..,k.
\end{split}
\label{DEIM_M_i}
\end{equation}
Clearly, this estimation is less computationally expensive then (\ref{tensor_POD})  since the summation stops at $m \ll n$.

During the numerical experiments we observed that tensors ${\bf M}^i,~i=1,2,..,k$ calculated in $POD/DEIM$ off-line stage (\ref{DEIM_M_i}) are different than ${\bf M}^i,~i=1,2,..,k$ obtained in the tensorial POD case
(\ref{tensor_POD}), but $\bf \mathcal{M}$ the reduced nonlinear term estimations of $\tilde N({\bf \tilde x})$ are accurate for both methods. We also mention for both standard POD and POD/DEIM approaches
terms as ${\bf M}^i,~i=1,2,..,k$ are used only for reduced Jacobian computations. In the case of POD/DEIM method, this leads to different derivatives values than in the case of tensorial POD or standard POD but the output solution error results are accurate as will see in Section $6$.

\section{The Shallow Water Equations}\label{sec:SWE}

In meteorological and oceanographic problems, one is often not interested in small time steps because the discretization error in time is small compared to the discretization error in space. SWE can be used to model Rossby and Kelvin waves in the atmosphere, rivers, lakes and oceans as well as gravity waves in a smaller domain. The alternating direction fully implicit (ADI) scheme \citet{Gus1971} considered in this paper is first order in both time and space and it is stable for large CFL condition numbers (we tested the stability of the scheme for a CFL condition number equal up to $8.9301$). It was also proved that the method is unconditionally stable for the linearized version of the SWE model. Other research work on this topic include efforts of \citet{FN1980,NVG1986}).

We are solving the SWE model using the $\beta$-plane approximation on a rectangular domain \citet{Gus1971}
\begin{equation}\label{eqn:swe-pde}
\frac{\partial w}{\partial t}=A(w)\frac{\partial w}{\partial x}+B(w)\frac{\partial w}{\partial y}+C(y)w,
\quad (x,y) \in [0,L] \times [0,D], \quad t\in(0,t_{\rm f}],
\end{equation}
where $w=(u,v,\phi)^T$ is a vector function, $u,v$ are the velocity components in the $x$ and $y$ directions, respectively, $h$ is the depth of the fluid, $g$ is the acceleration due to gravity, and $\phi = 2\sqrt{gh}$.

The matrices $A$, $B$ and $C$ are
\[
A=-\left(\begin{array}{ccc}
           u&0&\phi/2\\
           0&u&0\\
           \phi/2&0&u \end{array}\right), \quad
B=-\left(\begin{array}{ccc}
           v&0&0\\
           0&v&\phi/2\\
           0&\phi/2&v \end{array}\right), \quad
C=\left(\begin{array}{rrr}
           0&f&0\\
           -f&0&0\\
           0&0&0 \end{array}\right),
\]
where $f$ is the Coriolis term
\[
f=\hat f + \beta(y-D/2),~\beta=\frac{\partial f}{\partial y}, \quad \forall\, y,
\]
with $\hat f$ and $\beta$ constants.

We assume periodic solutions in the $x$ direction for all three state variables
while in the $y$ direction
$$v(x,0,t)=v(x,D,t)=0,~x\in[0,L],~t\in(0,t_{\rm f}]$$
and Neumann boundary condition are considered for $u$ and $\phi$.

Initially $w(x,y,0)=\psi(x,y),~\psi:\mathbb{R}\times\mathbb{R}\rightarrow \mathbb{R},~(x,y)\in[0,L]\times[0,D]$.
Now we introduce a mesh of $n = N_x\cdot N_y$ equidistant points on $[0,L]\times[0,D]$, with $\Delta x=L/(N_x-1),~\Delta y=D/(N_y-1)$. We also discretize the time interval $[0,t_{\rm f}]$ using $N_t$ equally distributed points and $\Delta t=t_{\rm f}/(N_t-1)$. Next we define vectors of unknown variables of dimension $n$ containing approximate solutions such as

$${\boldsymbol w}(t_N)\approx [w(x_i,y_j,t_N)]_{i=1,2,..,N_x,~j=1,2,..,N_y} \in \mathbb{R}^{n},~N=1,2,..N_t. $$

The semi-discrete equations of SWE \eqref{eqn:swe-pde} are:
\begin{eqnarray*}\label{eqn:swe-sd}
 {\bf u}' & = & -F_{11}({\bf u},\bm{{\phi}})-F_{12}({\bf u},{\bf v}) + {\bf F}\odot {\bf v}, \\
  {\bf v}' & = & -F_{21}({\bf u},{\bf v})-F_{22}({\bf v},\bm{{\phi}}) - {\bf F}\odot {\bf u}, \\
  {\bm{{\phi}}}' & = & -F_{31}({\bf u},\bm{{\phi}})-F_{32}({\bf v},\bm{{\phi}}),
\end{eqnarray*}
where $\odot$ is the Matlab componentwise multiplication operator, ${\bf u}'$, ${\bf v}'$, ${\bm{{\phi}}}'$ denote semi-discrete time derivatives, ${\bf F} =[\underbrace{{\bf f},{\bf f},..,{\bf f}}_{N_x}]$ stores Coriolis components ${\bf f} = [f(y_j)]_{j=1,2,..,N_y}$ while the nonlinear terms $F_{i1}$ and $F_{i2}$,~$i=1,2,3$, involving derivatives in $x$ and $y$ directions, respectively, are defined as follows:

\begin{eqnarray*}\label{eqn:swe-sd}
  && F_{i1},F_{i2}: \mathbb{R}^{n} \times \mathbb{R}^{n} \rightarrow \mathbb{R}^{n},~i=1,2,3,~F_{11}({\bf u},{\boldsymbol \phi})={\boldsymbol u}\odot A_x{\boldsymbol u}+\frac{1}{2}{\boldsymbol \phi}\odot A_x{\boldsymbol\phi},\\
  && F_{12}({\boldsymbol u},{\boldsymbol v})={\boldsymbol v}\odot A_y{\boldsymbol u}, F_{21}({\boldsymbol u},{ \boldsymbol v})={\boldsymbol u}\odot A_x{\boldsymbol v}, F_{22}({\boldsymbol v},{\boldsymbol\phi})={\boldsymbol v}\odot A_y{\boldsymbol v}+\frac{1}{2}{\boldsymbol \phi}\odot A_y{\boldsymbol\phi}, \\
  &&F_{31}({\boldsymbol u},{\boldsymbol \phi})=\frac{1}{2}{\boldsymbol\phi} \odot A_x {\boldsymbol u}+{\boldsymbol u} \odot {A_x\boldsymbol \phi},~F_{32}({\boldsymbol v},{\boldsymbol \phi})=\frac{1}{2}{\boldsymbol\phi}\odot A_y{\boldsymbol v}+{\boldsymbol v} \odot A_y{\boldsymbol\phi}.\\
\end{eqnarray*}

Here $A_x,A_y\in \mathbb{R}^{n\times n}$ are constant coefficient matrices for discrete first-order and second-order differential operators which take into account the boundary conditions.

The numerical scheme was implemented in Fortran and uses a sparse matrix environment. For operations with sparse matrices we employed SPARSEKIT library \citet{Saad1994}
and the sparse linear systems obtained during the quasi-Newton iterations were solved using MGMRES library \citet{Barrett94,Kelley95,Saad2003}. Here we didn't decouple the model equations like in Stefanescu and Navon \cite{Stefanescu2013} where the Jacobian is either block cyclic tridiagonal or
block tridiagonal. We followed this approach since we plan to implement a 4D-Var data assimilation system based on ADI SWE and the adjoints of the
decoupled systems can't be solved with the same implicit scheme applied for solving the forward model.

\section{Reduced Order Shallow Water Equation Models}\label{sec:ROMS_SWE}

Here we will not describe the entire standard POD SWE, tensorial POD SWE and POD/DEIM SWE discrete models but only we introduce the projected nonlinear term $F_{11}$
for all three ROMs. ADI discrete equations were projected onto reduced POD subspaces and a detailed description of the reduced equations for standard POD and POD/DEIM
is available in \citet{Stefanescu2013}.

Depending on the type of reduced approaches, the Petrov-Galerkin projected nonlinear term $\tilde F_{11}$ has the following form

\paragraph{Standard POD.}
\begin{equation}
 \tilde F_{11} = W^TF_{11}=\underbrace{W^T}_{k \times n} \biggl(\underbrace{(U\tilde {\bf u})\odot(U_x\tilde {\bf u})}_{n \times 1}\biggr)+\frac{1}{2}\underbrace{W^T}_{k \times n}
 \biggl(\underbrace{(\Phi\tilde {\boldsymbol \phi})\odot(\Phi_x\tilde {\boldsymbol \phi})}_{n \times 1}\biggr),
\end{equation}
where $U$ and $\Phi$ contains the POD bases corresponding to state variables $u$ and $\phi$ while the POD basis derivatives are included in $U_x=A_xU \in \mathbb{R}^{n \times k}$ and
$\Phi_x=A_x\Phi\in \mathbb{R}^{n \times k}$.

\paragraph{Tensorial POD.}
\begin{equation}
 \tilde F_{11} = W^TF_{11} \in \mathbb{R}^k;~\big[\tilde F_{11}\big]_i= \langle{\bf M}^i_1,{\bf \tilde U}\rangle_{Frobenius} + \langle{\bf M}^i_2,\tilde {\boldsymbol\Phi}\rangle_{Frobenius},~i=1,2,..,k;
\end{equation}
\begin{equation*}
\begin{split}
& \tilde {\bf U} = \big[{\tilde {\bf U}}_{i,j}\big]_{i,j=1,..,k}\in \mathbb{R}^{k\times k};~\tilde {\bf U}_{i,j}=\tilde {\bf u}_{i}\tilde{\bf u}_{j} \in \mathbb{R},~\tilde {\boldsymbol \Phi} = \big[{\tilde {\boldsymbol \Phi}}_{i,j}\big]_{i,j=1,..,k}\in \mathbb{R}^{k\times k};~\tilde {\boldsymbol \Phi}_{i,j}=\tilde {\boldsymbol \phi}_{i}\tilde{\boldsymbol \phi}_{j} \in \mathbb{R},
\end{split}
\end{equation*}
where $\tilde{\bf u }\in \mathbb{R}^k$ and $\tilde{\boldsymbol \phi}\in \mathbb{R}^k$ are reduced state variables.
\begin{equation}
\begin{split}
& {\bf M}^i_1 = \big[{\bf M}^i_{1i_1i_2}\big]_{i_1,i_2=1,..,k}\in \mathbb{R}^{k\times k};~{\bf M}^i_{1i_1i_2} = \sum_{l=1}^nW_{li}U_{li_1}{U_x}_{li_2} \in \mathbb{R}\\
& {\bf M}^i_2 = \big[{\bf M}^i_{2i_1i_2}\big]_{i_1,i_2=1,..,k}\in \mathbb{R}^{k\times k};~{\bf M}^i_{2i_1i_2} = \sum_{l=1}^nW_{li}\Phi_{li_1}{\Phi_x}_{li_2} \in \mathbb{R},
\end{split}\label{tensor_pod_coeff}
\end{equation}
and $U_x$ and $\Phi_x$ were defined above.

\paragraph{POD/DEIM.}
\begin{equation}
\tilde F_{11} \approx \underbrace{W^TV_{F_{11}}(P_{F_{11}}^TV_{F_{11}})^{-1}}_{{\rm precomputed}~k\times m} \bigg(\underbrace{(P_{F_{11}}^TU\tilde{\bf u})\odot(P_{F_{11}}^TU_x)\tilde{\bf u}}_{m \times 1} + \underbrace{(P_{F_{11}}^T\Phi\tilde {\boldsymbol \phi})\odot(P_{F_{11}}^T\Phi_x\tilde {\boldsymbol \phi})}_{m \times 1}\bigg),
\label{POD_DEIM_nonlinearity1}
\end{equation}
where $V_{F_{11}}\in \mathbb{R}^{n \times m}$ collects the first $m$ POD basis modes of nonlinear function ${ F_{11}}$ while $P_{F_{11}} \in \mathbb{R}^{n \times m}$ is the DEIM interpolation selection matrix. Let us denote the precomputed term by $E_{11}= W^TV_{F_{11}}(P_{F_{11}}^TV_{F_{11}})^{-1}$.

Tensors like ${\bf M}^i_1$ and ${\bf M}^i_2$ \eqref{tensor_pod_coeff} must also be computed in the case of standard POD and POD/DEIM since the analytic form of reduce Jacobian was employed. This approach reduces the CPU time of standard POD since usually the reduced Jacobians are obtained by projecting the full Jacobian at every time step. A generalization of DEIM to approximate operators is not been yet developed but has the ability to decrease more the computational complexity of POD/DEIM approach. Some related work includes \citet{Tonn_2011} who developed Multi-Component EIM for deriving affine approximations for continuous vector valued functions. \citet{Wirtz2012} introduced the matrix-DEIM  approach to approximate the Jacobian of a nonlinear function.
\citet{Cha2008} proposed a sampling strategy centered on the trajectory of the nonlinear functions in order to approximate the reduced Jacobian. An extension for nonlinear problems that do not have component-wise dependence on the state has been introduced in \citet{Zhou_2012}.

Table \ref{table2} contains the procedure list required by all three algorithms in the off-line stage. ${\bf M}_1^i$, ${\bf M}_2^i$ and $E_{11}$ are POD and POD/DEIM coefficients related to nonlinear term $F_{11}$, similar coefficients being required for computation of other reduced nonlinear terms.
\begingroup
\begin{table}[h]
\begin{tabular}{|p{0.3\textwidth}|p{0.3\textwidth}|p{0.3\textwidth}|}\hline
 Standard POD & Tensorial POD & POD/DEIM\\ \hline
 Generate snapshots & Generate snapshots & Generate snapshots\\ \hline
 SVD for $u,v,~\phi$ & SVD for $u,v,~\phi$ & SVD for $u,v,~\phi$\\ \hline
 -- & -- & SVD for all nonlinear terms\\ \hline
 -- & -- & DEIM index points for all nonlinear terms\\ \hline
 Calc. POD coefficients ${\bf M}^i_j$ (\ref{tensor_pod_coeff}) & Calc. POD coefficients ${\bf M}^i_j$ (\ref{tensor_pod_coeff}) & Calc. POD coefficients ${\bf M}^i_j$ (\ref{DEIM_M_i})\\
 (reduced Jac. calc.) & (reduced Jac. and right-hand side terms calc.) & (reduced Jac. calc.) \\  \hline
 -- & -- & Calc. all POD/DEIM coef. such as $E_{11}$\\ \hline
\end{tabular}
\caption{\label{table2} ROMs off-line stage procedures - POD coefficients ${\bf M}^i_1,~{\bf M}^i_2$ are required for reduced Jacobian calculation. Only tensorial POD uses them also for right-hand side terms computations during the quasi-Newton iterations required by Gustafsson's nonlinear ADI finite difference scheme.}
\end{table} %
\endgroup%


\section{Numerical Results }\label{sec:Numerical_resuls}

For all tests we derived the initial conditions from the initial height condition No. 1 of Grammeltvedt 1969 \cite{Gram1969} i.e.
\begin{equation*}
\hspace{-10mm}h(x,y,0)=H_0+H_1+\tanh\biggl(9\frac{D/2-y}{2D}\biggr)+H_2\textrm{sech}^2\biggl(9\frac{D/2-y}{2D}\biggr)\sin\biggl(\frac{2\pi x }{L}\biggr),
\end{equation*}
The initial velocity fields are derived from the initial height field using the geostrophic relationship

\[
u = \biggl(\frac{-g}{f}\biggr)\frac{\partial h}{\partial y}, \quad v = \biggl(\frac{g}{f}\biggr)\frac{\partial h}{\partial x}.
\]

\begin{figure}[h]
  \centering
  \subfigure[Geopotential height field ] {\includegraphics[scale=0.35]{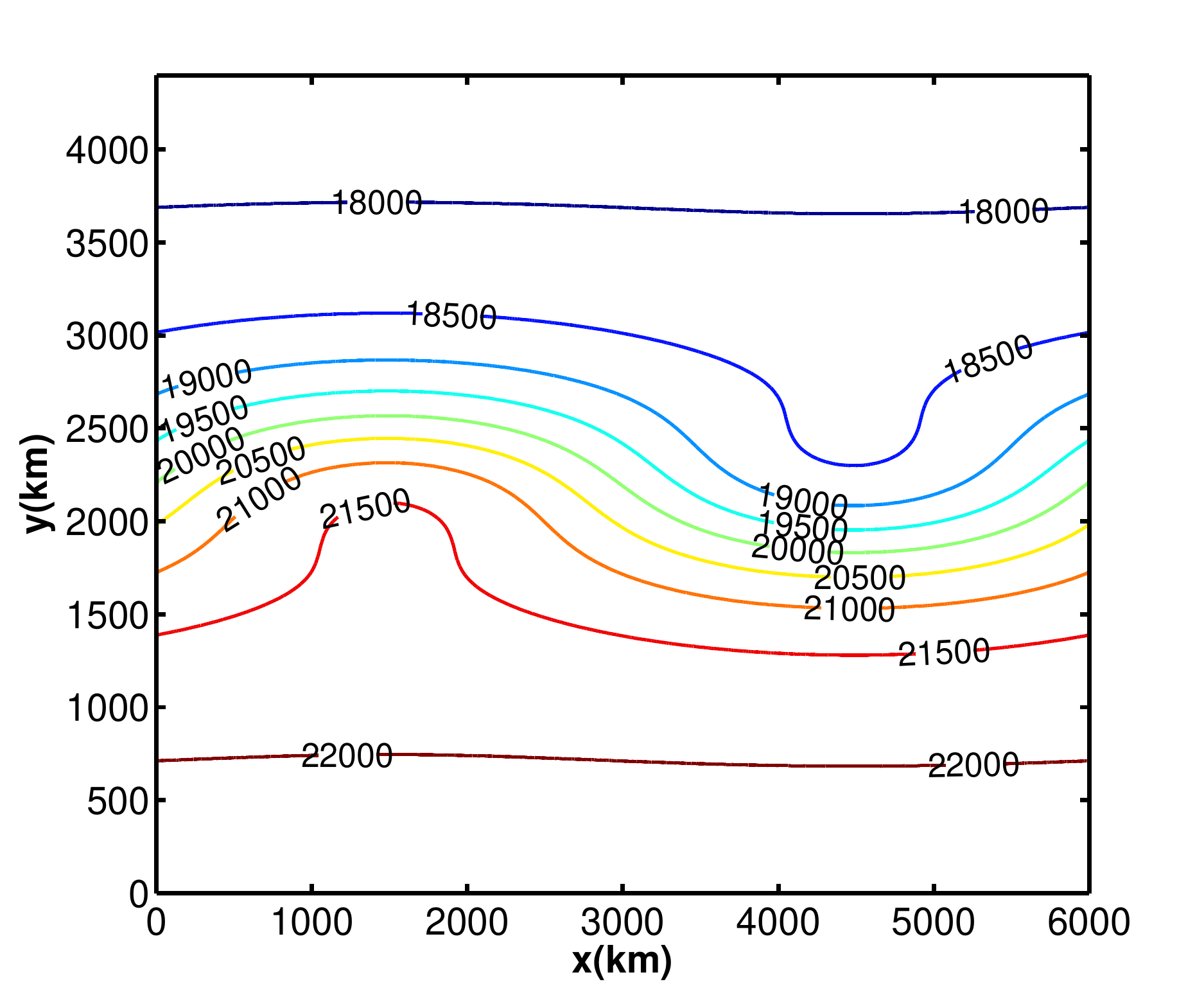}}
  \subfigure[Windfield]{\includegraphics[scale=0.35]{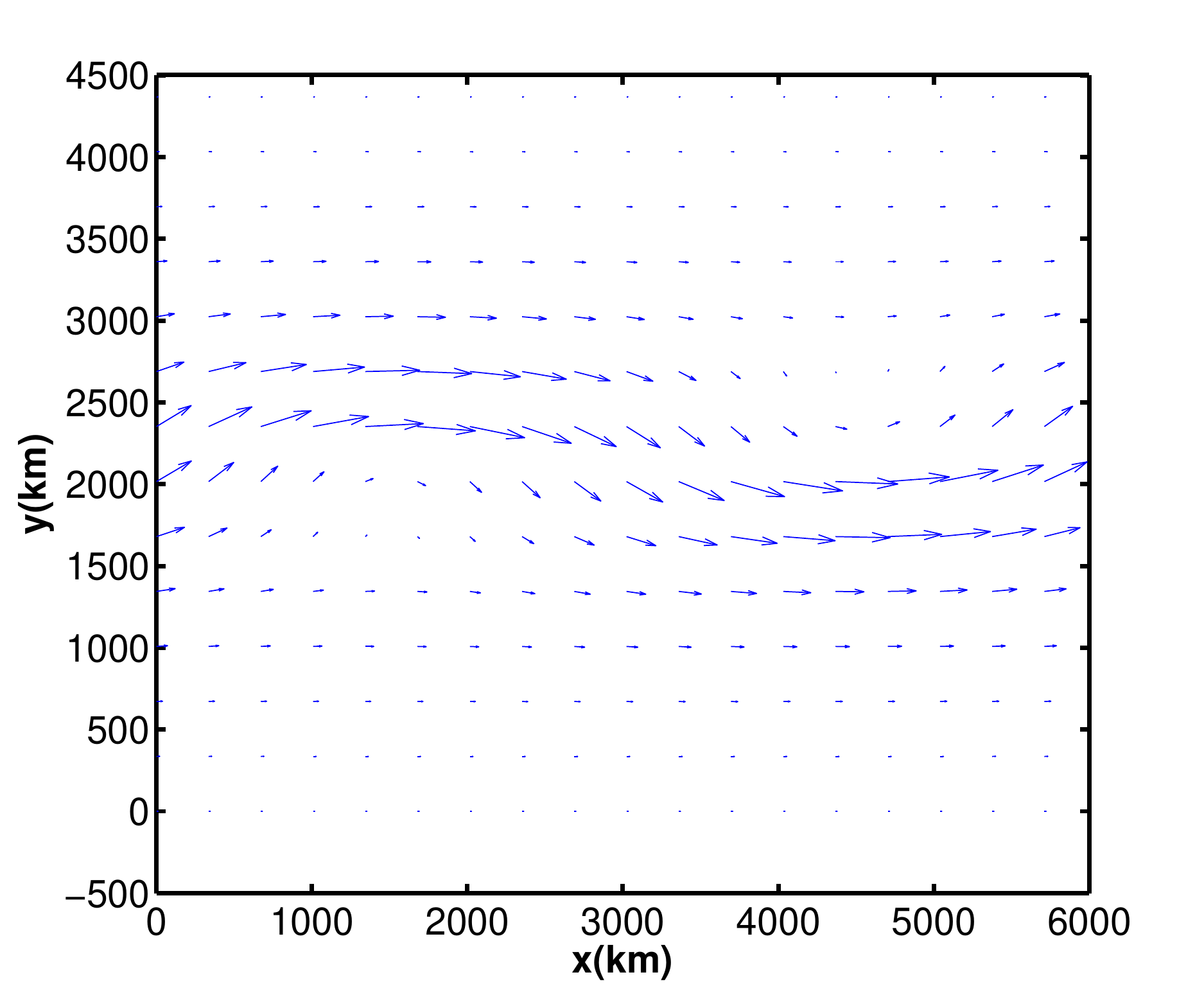}}
\caption{Initial condition: Geopotential height field for the Grammeltvedt initial condition and windfield (the velocity unit is 1km/s)
calculated from the geopotential field using the geostrophic approximation. \label{Fig::1}}
\end{figure}

We use the following constants $ L=6000km,~D=4400km,~\hat f=10^{-4}s^{-1}~,\beta=1.5\cdot10^{-11}s^{-1}m^{-1},~g=10 m s^{-2},~H_0=2000m,~H_1=220m,~H_2=133m.$
Figure \ref{Fig::1} depicts the initial geopotential isolines and the geostrophic wind field.


Most of the depicted results are obtained in the case when the domain is discretized using a mesh of  $376\times276 = 103,776$ points, with $\Delta x = \Delta y = 16$km. We select two integration time windows of $24$h and $3$h and we use $91$ time steps ($NT=91$) with $\Delta t = 960$s and $\Delta t = 120$s.

ADI FD SWE scheme proposed by Gustafsson 1971 in \cite{Gus1971} is first employed in order to obtain the numerical solution of the SWE model. The implicit scheme allows us to integrate in time at a Courant-Friedrichs-Levy (CFL) condition of
$\sqrt{gh}(\Delta t/\Delta x) < 8.9301.$

The nonlinear algebraic systems of ADI FD SWE scheme is solved using quasi - Newton method, and the LU decomposition is performed every $6$ time steps.

\begin{figure}[h]
  \centering
  \subfigure[State variables $u,~v,~\phi$] {\includegraphics[scale=0.35]{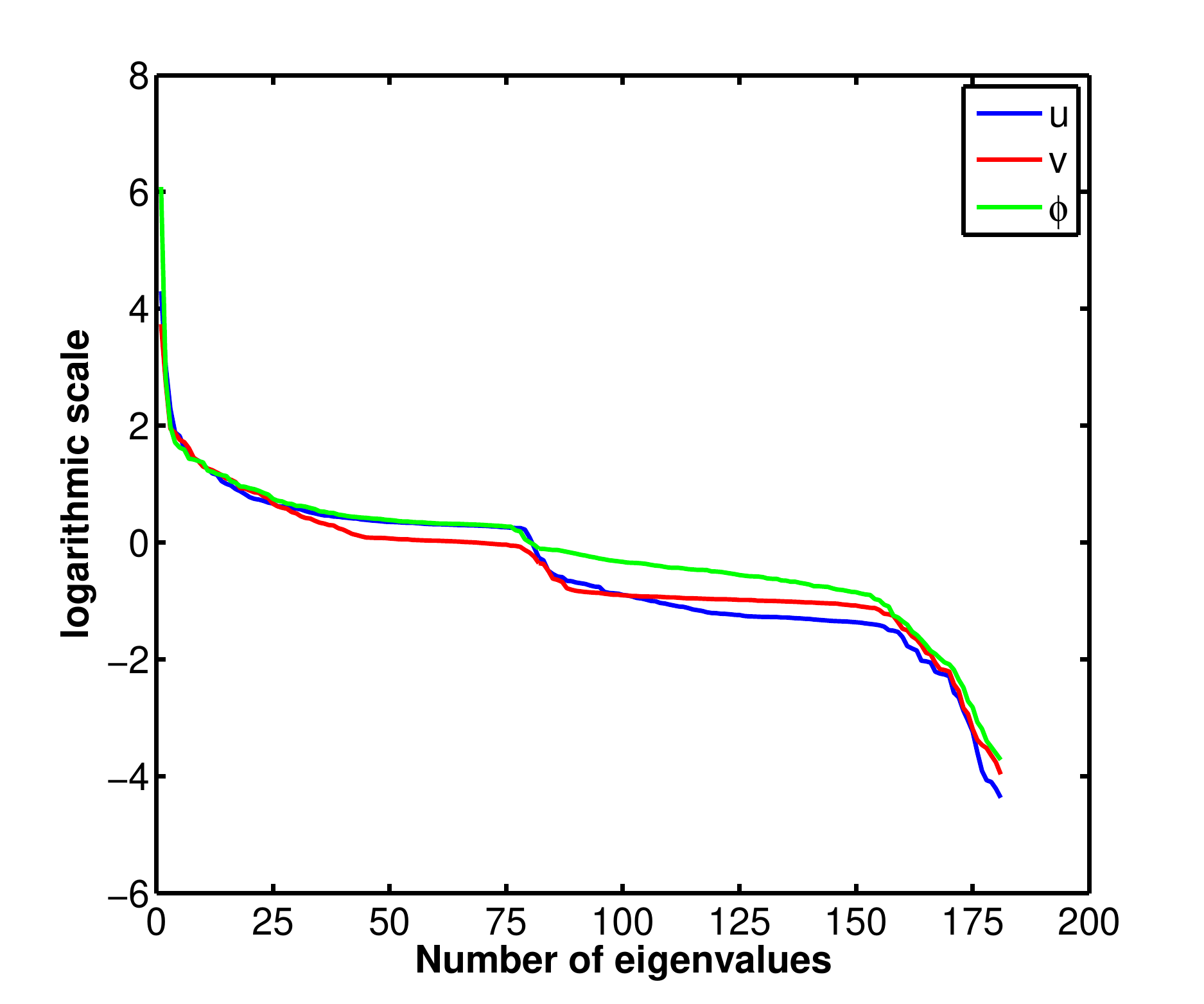}}
  \subfigure[Nonlinear terms]{\includegraphics[scale=0.35]{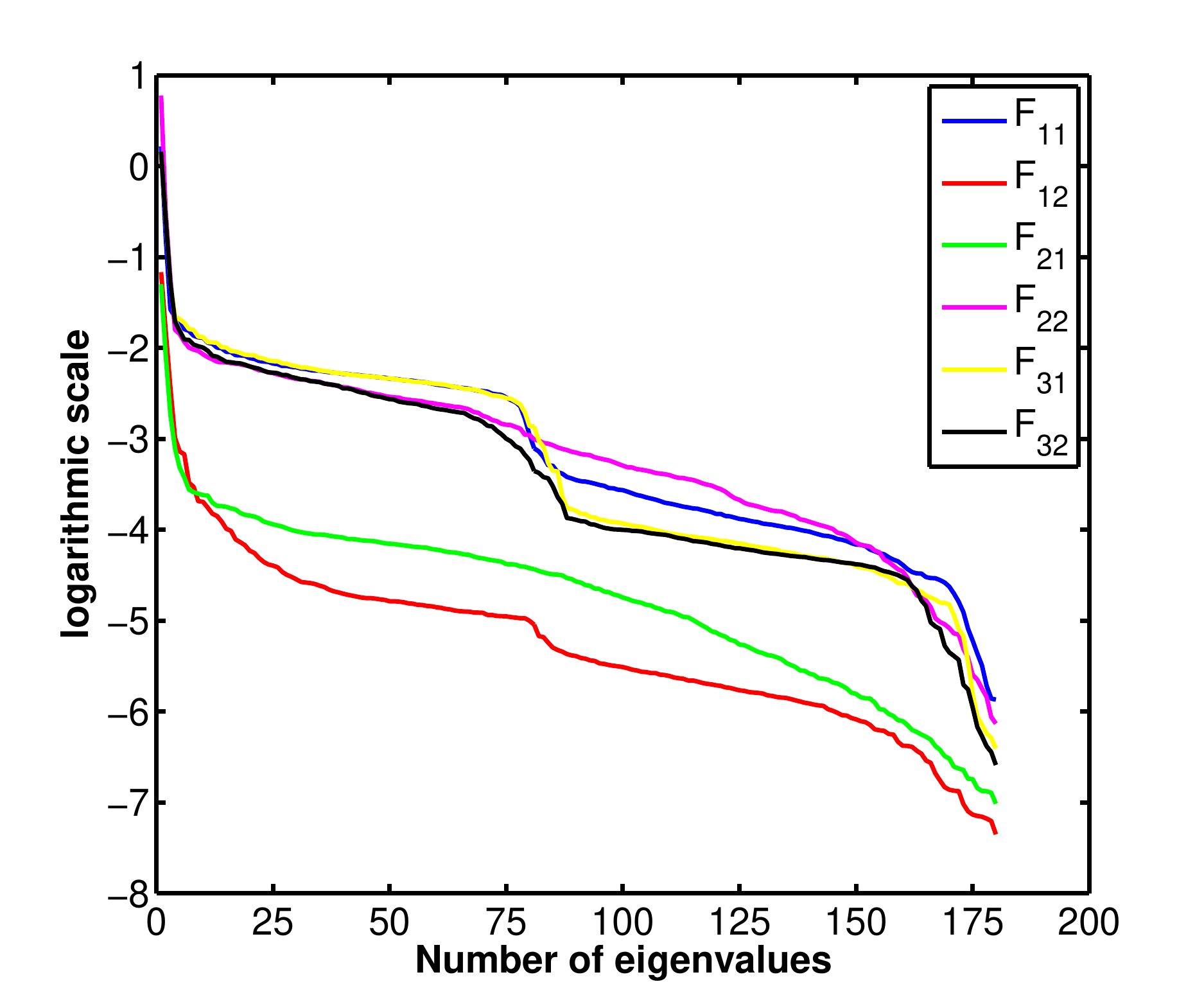}}
\caption{The decay around the singular values of the snapshots solutions for $u,~v,~\phi$ and nonlinear terms for $\Delta t = 960s$ and integration time window of $24h$ .\label{Fig::3}}
\end{figure}

\begin{figure}[h]
  \centering
  \subfigure[State variables $u,~v,~\phi$] {\includegraphics[scale=0.35]{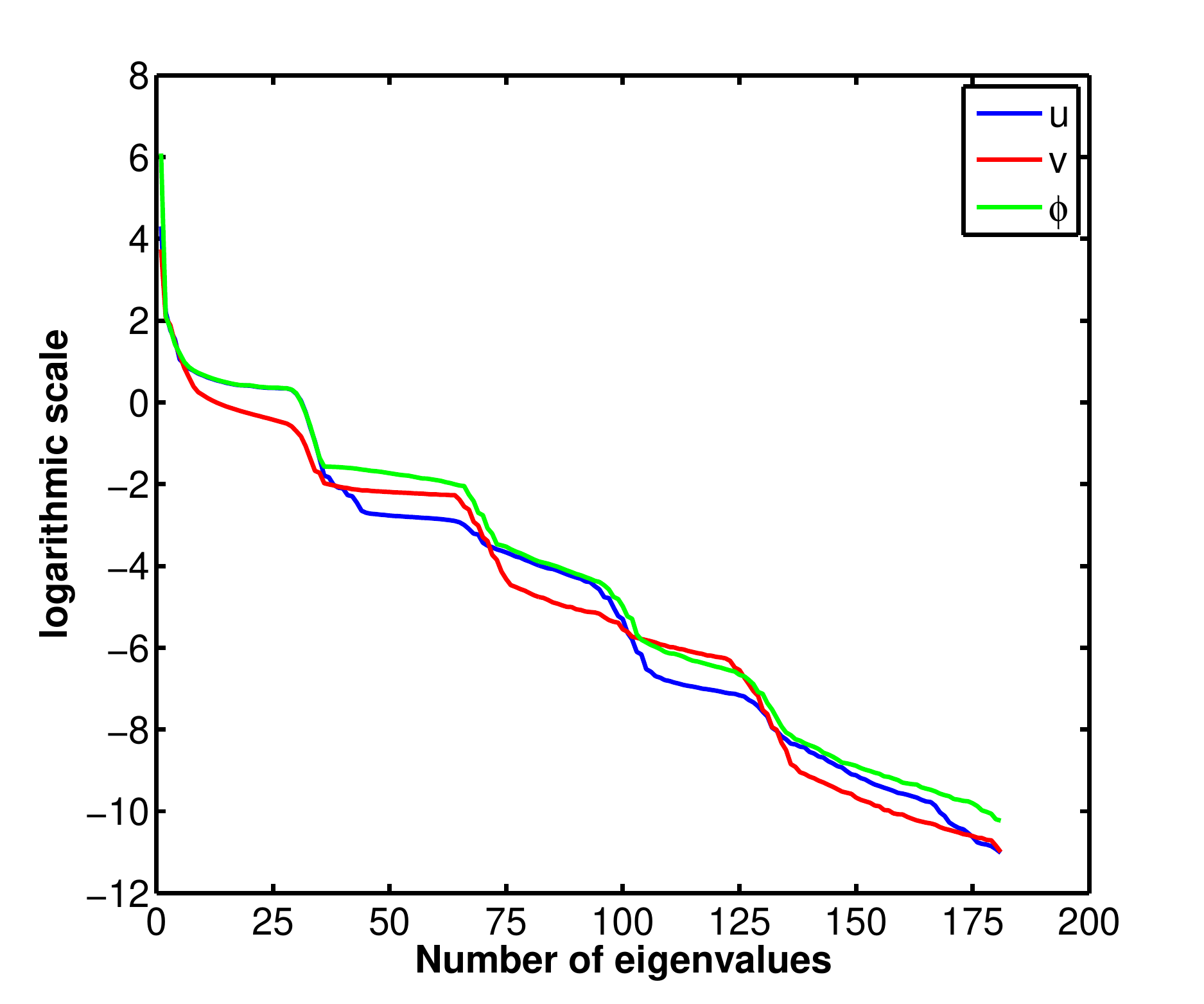}}
  \subfigure[Nonlinear terms]{\includegraphics[scale=0.36]{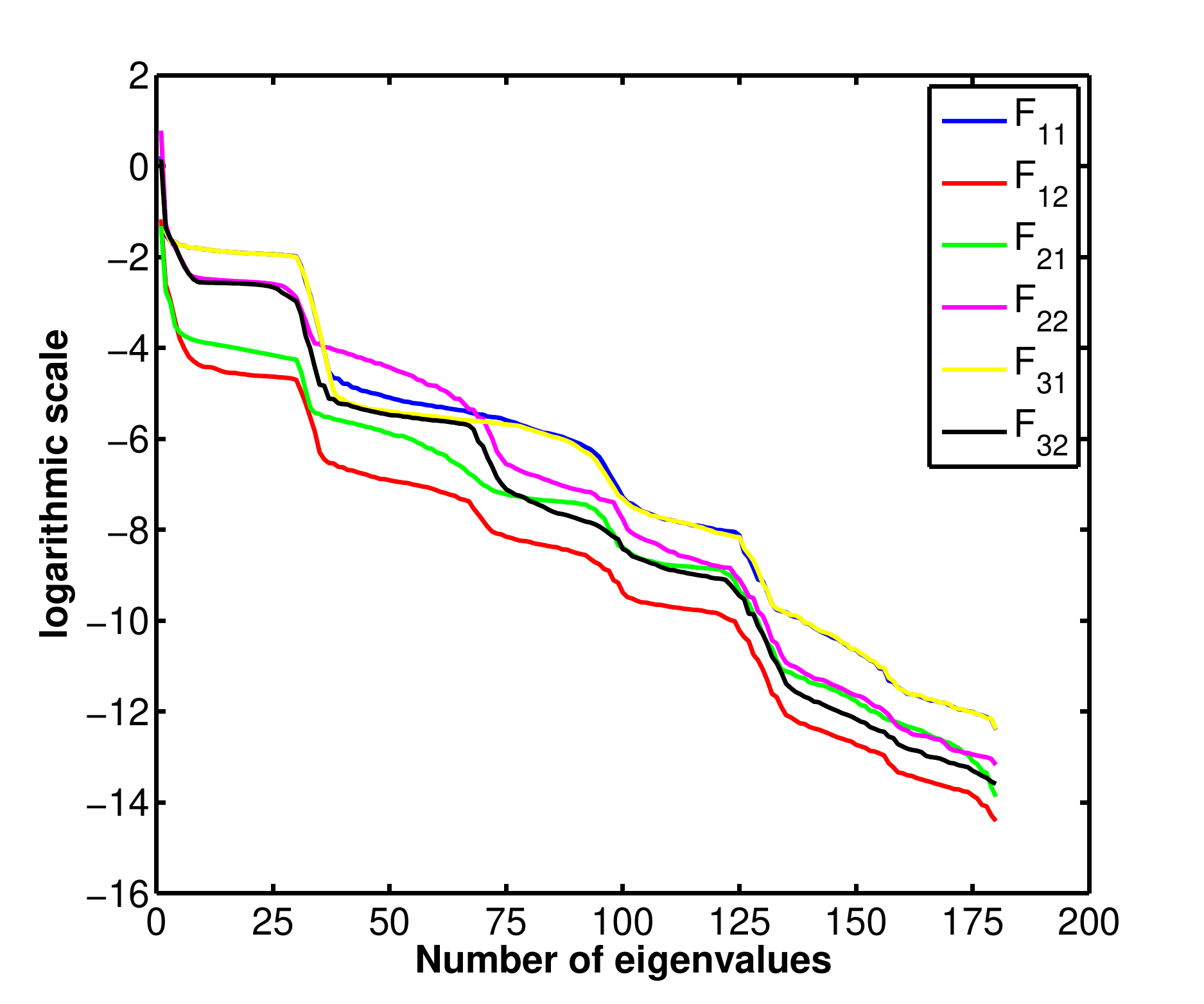}}
\caption{The decay around the singular values of the snapshots solutions for $u,~v,~\phi$ and nonlinear terms for $\Delta t = 120s$ and a time integration window of $3h$ .\label{Fig::4}}
\end{figure}

We derive the reduced order models by employing a Galerkin projection.The POD basis functions are constructed using $91$ snapshots (number of snapshots equal with the number of time steps $N_t$) obtained from the numerical solution of the full - order ADI FD SWE model at equally spaced time steps for each time interval $[0,24h]$ and $[0,3h]$. Figures \ref{Fig::3},{\ref{Fig::4}} show the decay around the eigenvalues of the snapshot solutions for $u,~v,~\phi$ and the nonlinear snapshots $F_{11},~F_{12},$ $F_{21},~F_{22}$, $F_{31},~F_{32}$. We notice that the singular values decay much faster when the model is integrated for $3$h. Consequently this translates in a more accurate solution representations for all three ROM methods using the same number of POD modes. For both time configurations and all tests in this study, the dimensions of the POD bases for each variable is taken to be $50$, capturing more than $99\%$ of the system energy. The largest neglected eigenvalues corresponding to state variables $u,~v,~\phi$ are $2.23,~1.16$ and $2.39$ for $t_{\rm f}=24$h and $0.0016,~0.0063$ and $0.0178$ for $t_{\rm f}=3$h, respectively.


%

\begin{figure}[h]
  \centering
  \subfigure[Nonlinear term $F_{11}$ ]{\includegraphics[scale=0.222]{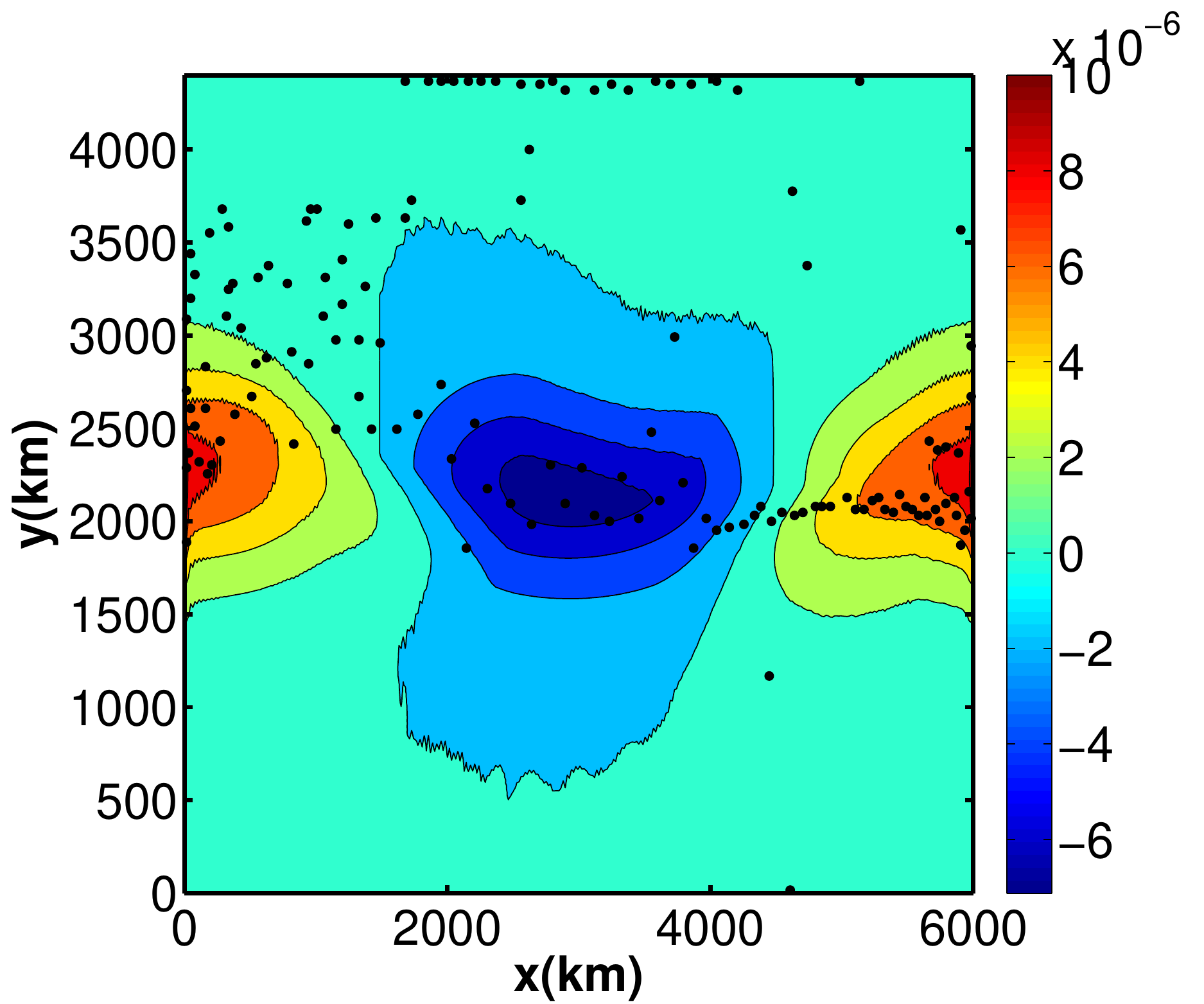}}\label{Fig::5a}
  \subfigure[Nonlinear term $F_{12}$]{\includegraphics[scale=0.222]{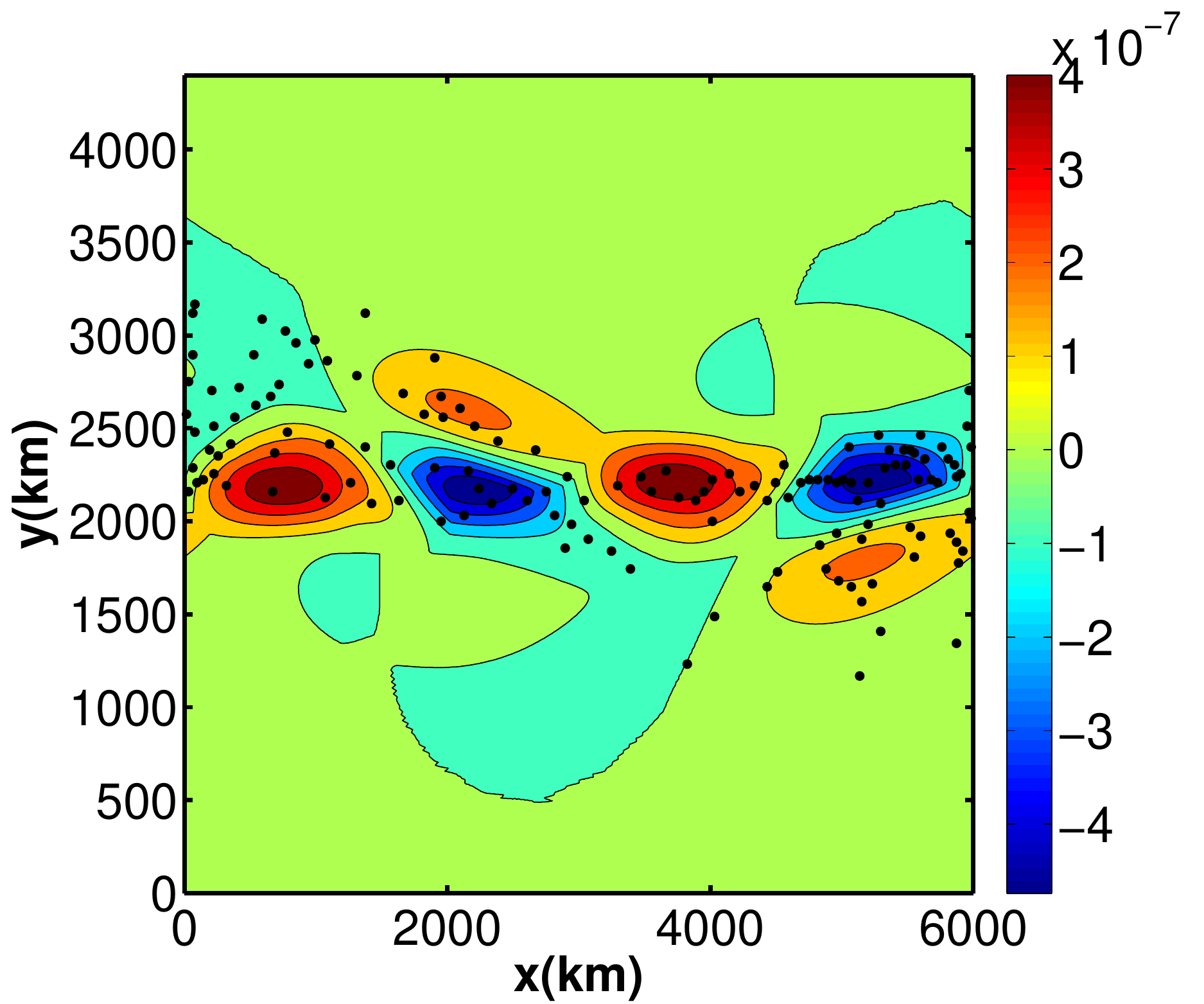}}\label{Fig::5b}
  \subfigure[Nonlinear term $F_{21}$]{\includegraphics[scale=0.222]{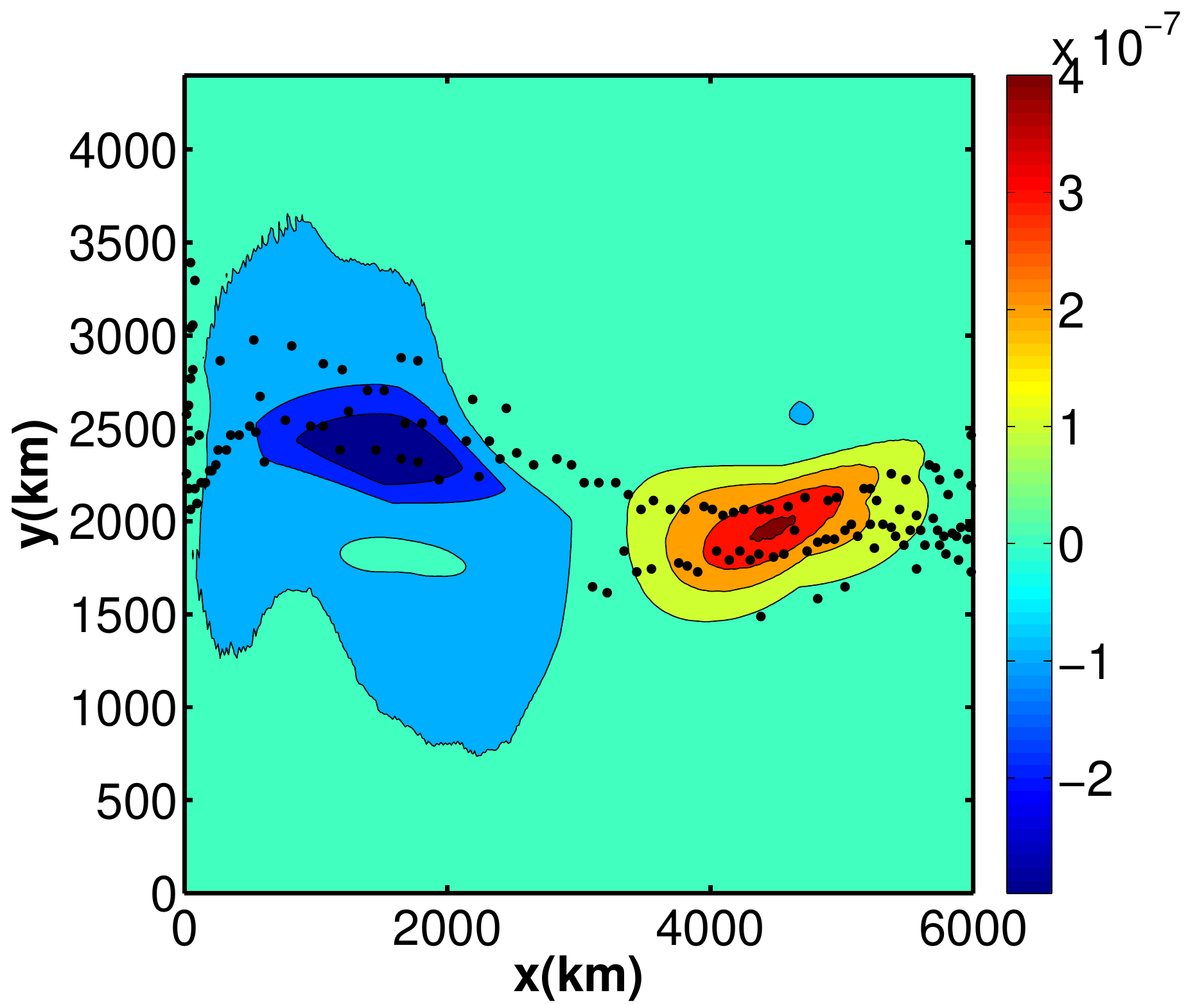}}\label{Fig::5c}
  \subfigure[Nonlinear term $F_{22}$]{\includegraphics[scale=0.222]{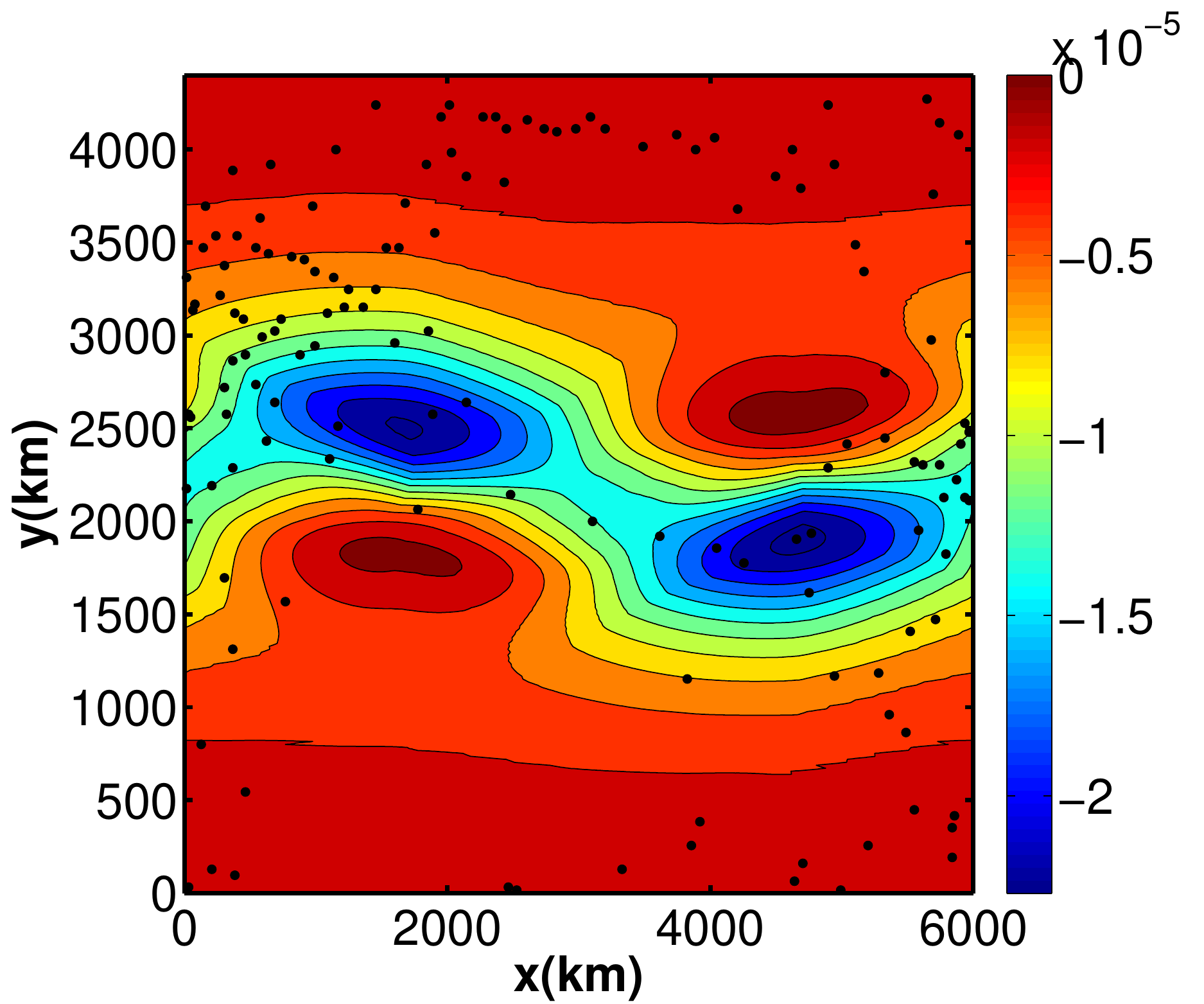}}\label{Fig::5d}
  \subfigure[Nonlinear term $F_{31}$]{\includegraphics[scale=0.222]{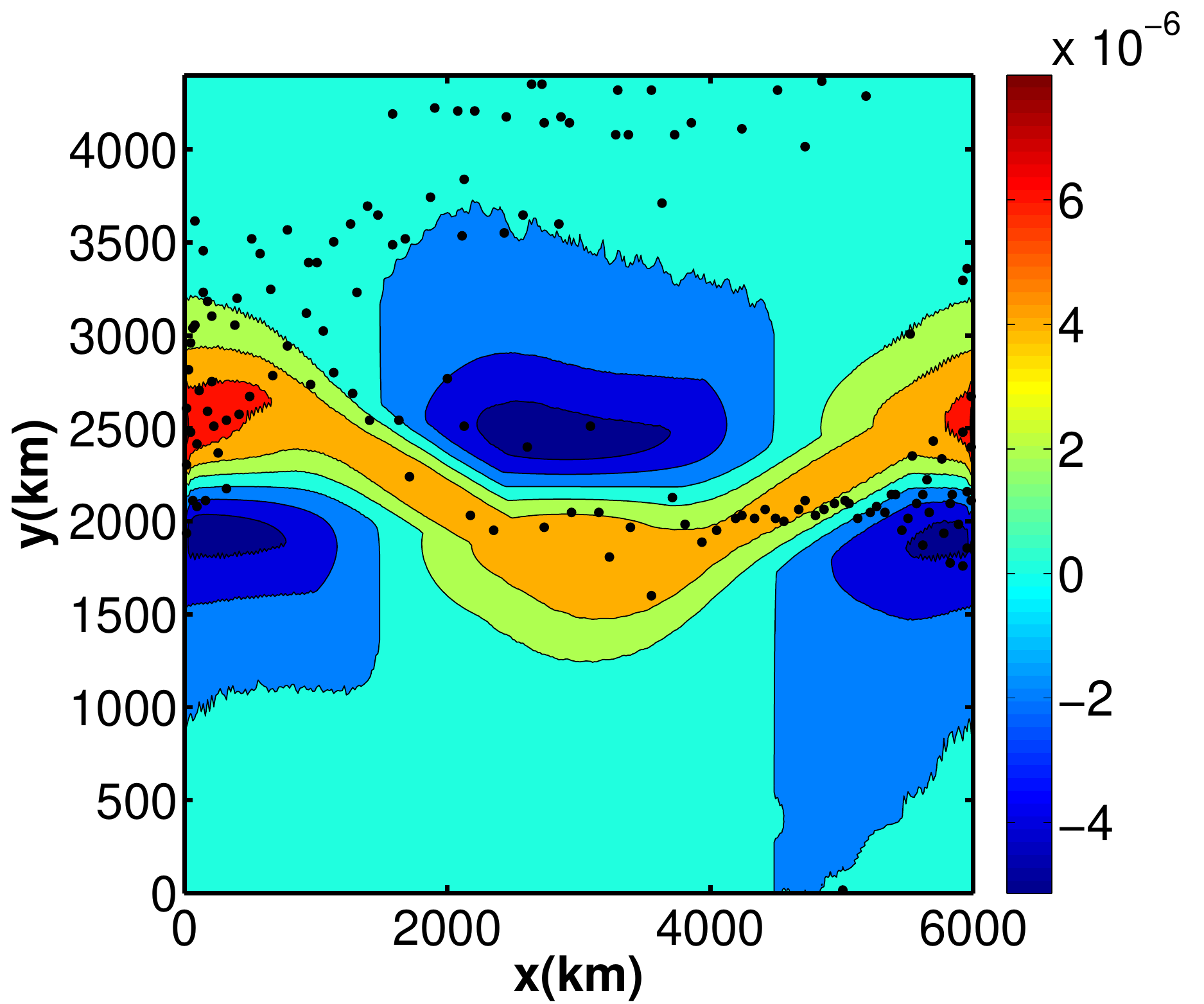}}\label{Fig::5e}
  \subfigure[Nonlinear term $F_{32}$]{\includegraphics[scale=0.222]{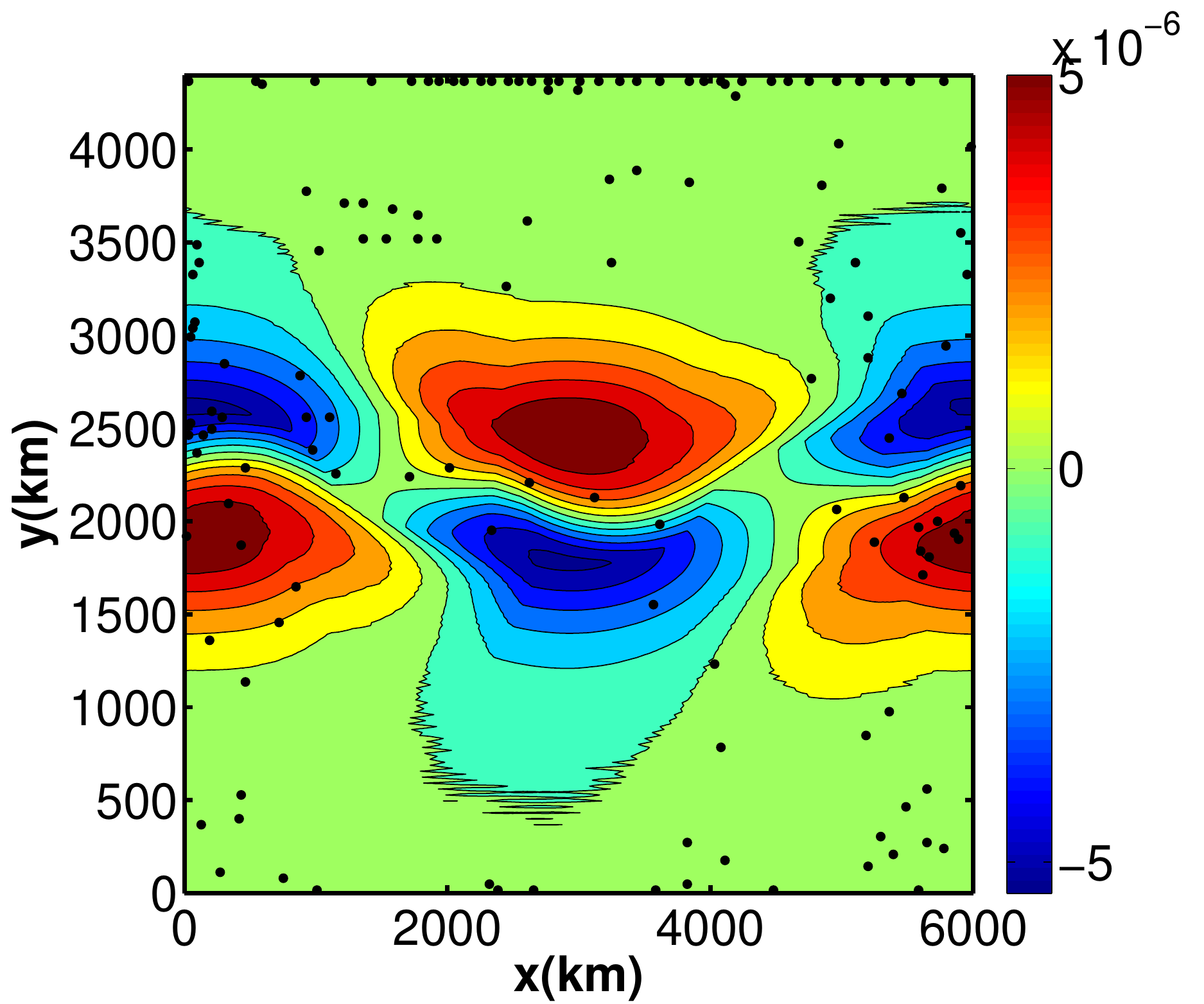}}\label{Fig::5f}
  \caption{The first 100 DEIM interpolation points corresponding to all nonlinear terms  in the SWE model for time integration window of $24h$. The background consists in isolines of the maximum values of the nonlinear terms over time.  \label{Fig::5}}
\end{figure}

\begin{figure}[h]
  \centering
  \subfigure[Nonlinear term $F_{11}$ ]{\includegraphics[scale=0.222]{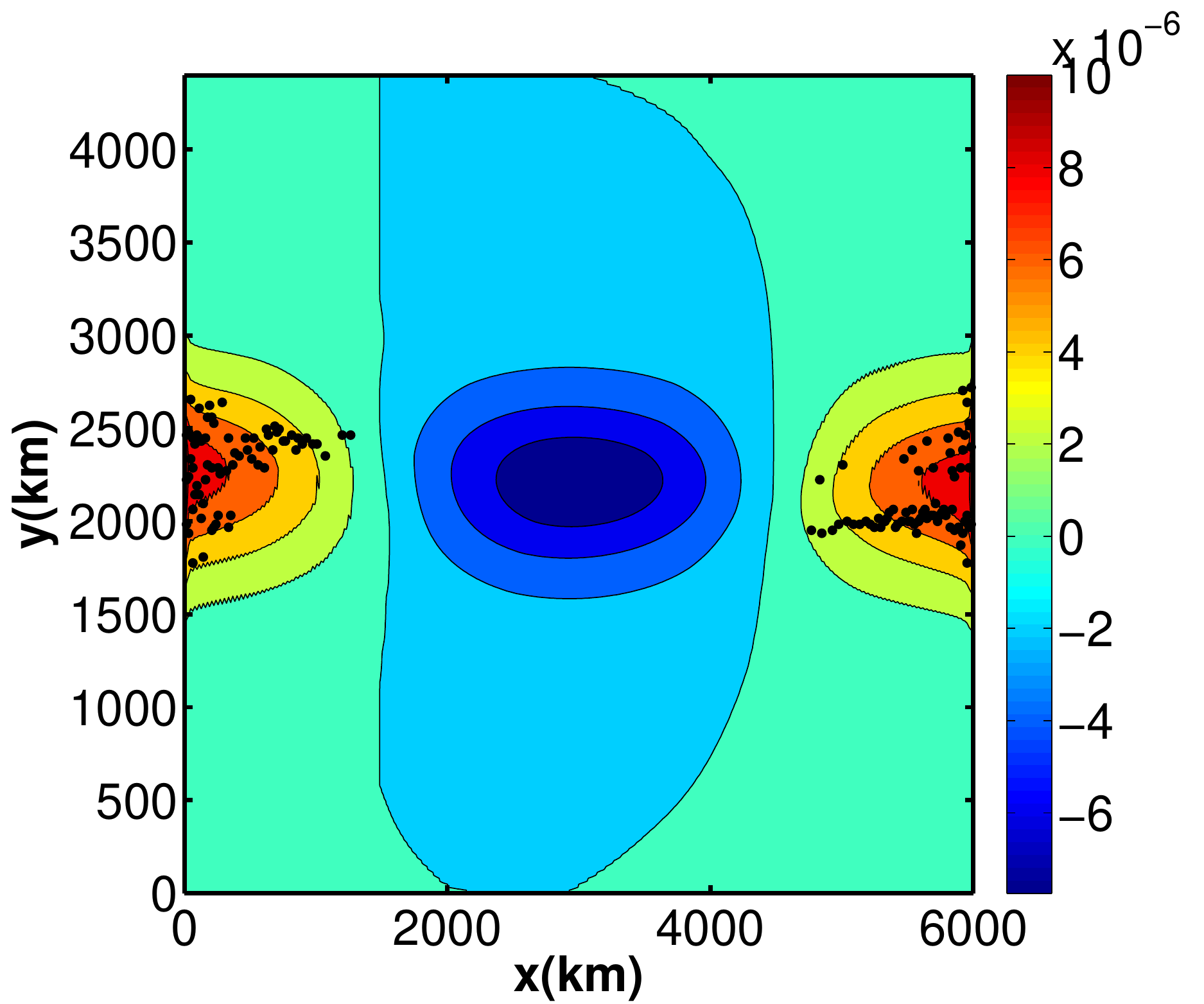}}\label{Fig::6a}
  \subfigure[Nonlinear term $F_{12}$]{\includegraphics[scale=0.222]{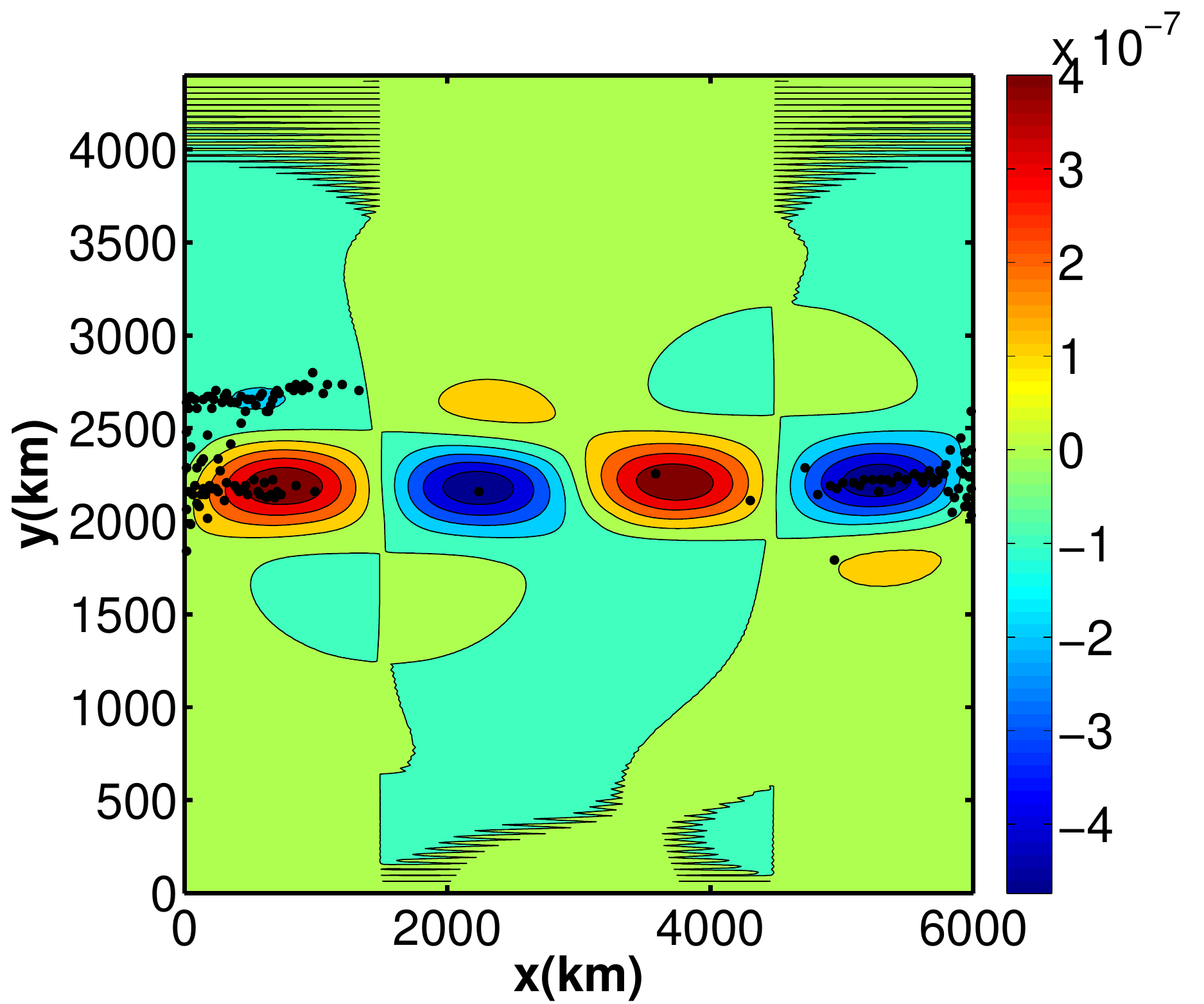}}\label{Fig::6b}
  \subfigure[Nonlinear term $F_{21}$]{\includegraphics[scale=0.222]{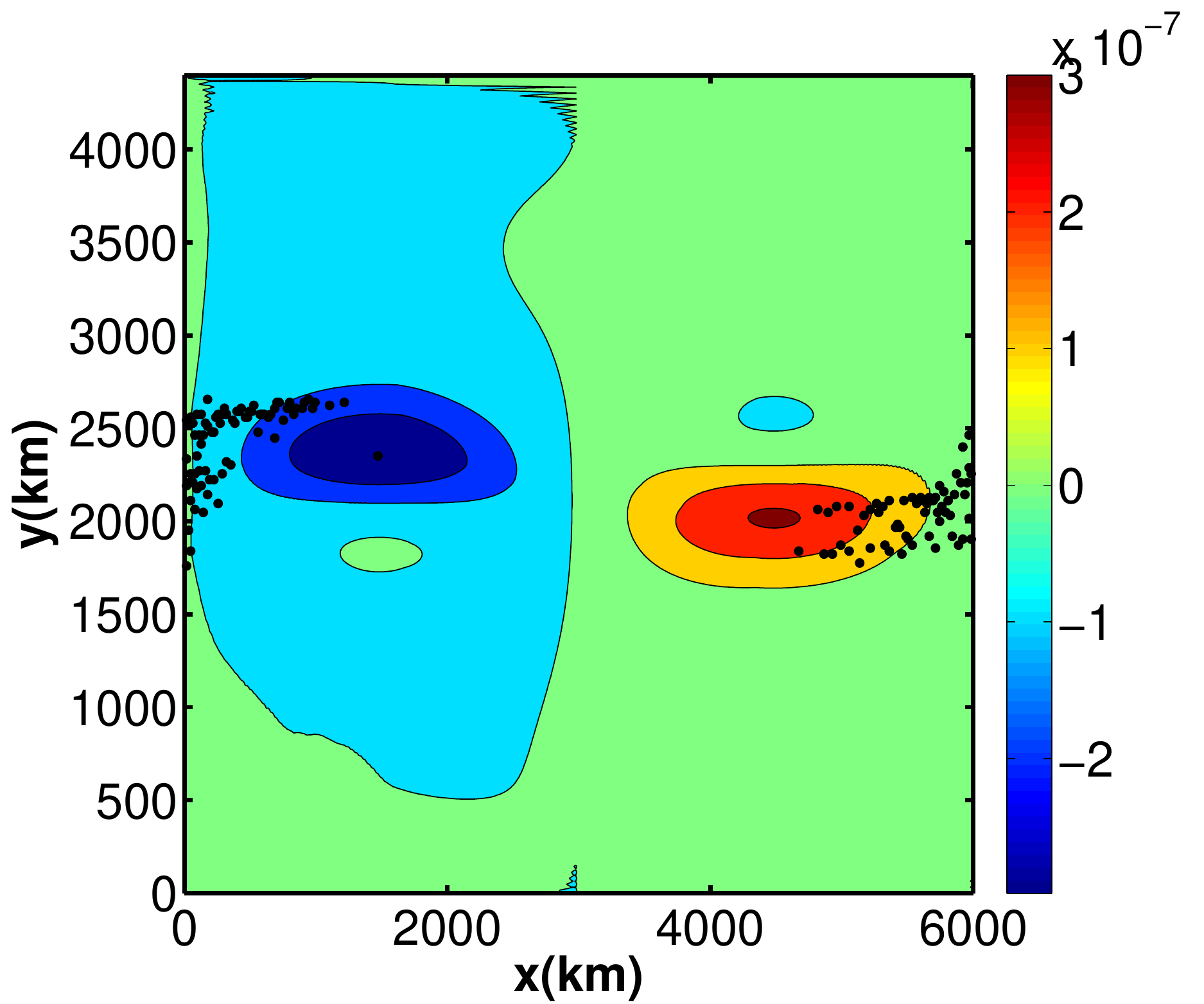}}\label{Fig::6c}
  \subfigure[Nonlinear term $F_{22}$]{\includegraphics[scale=0.222]{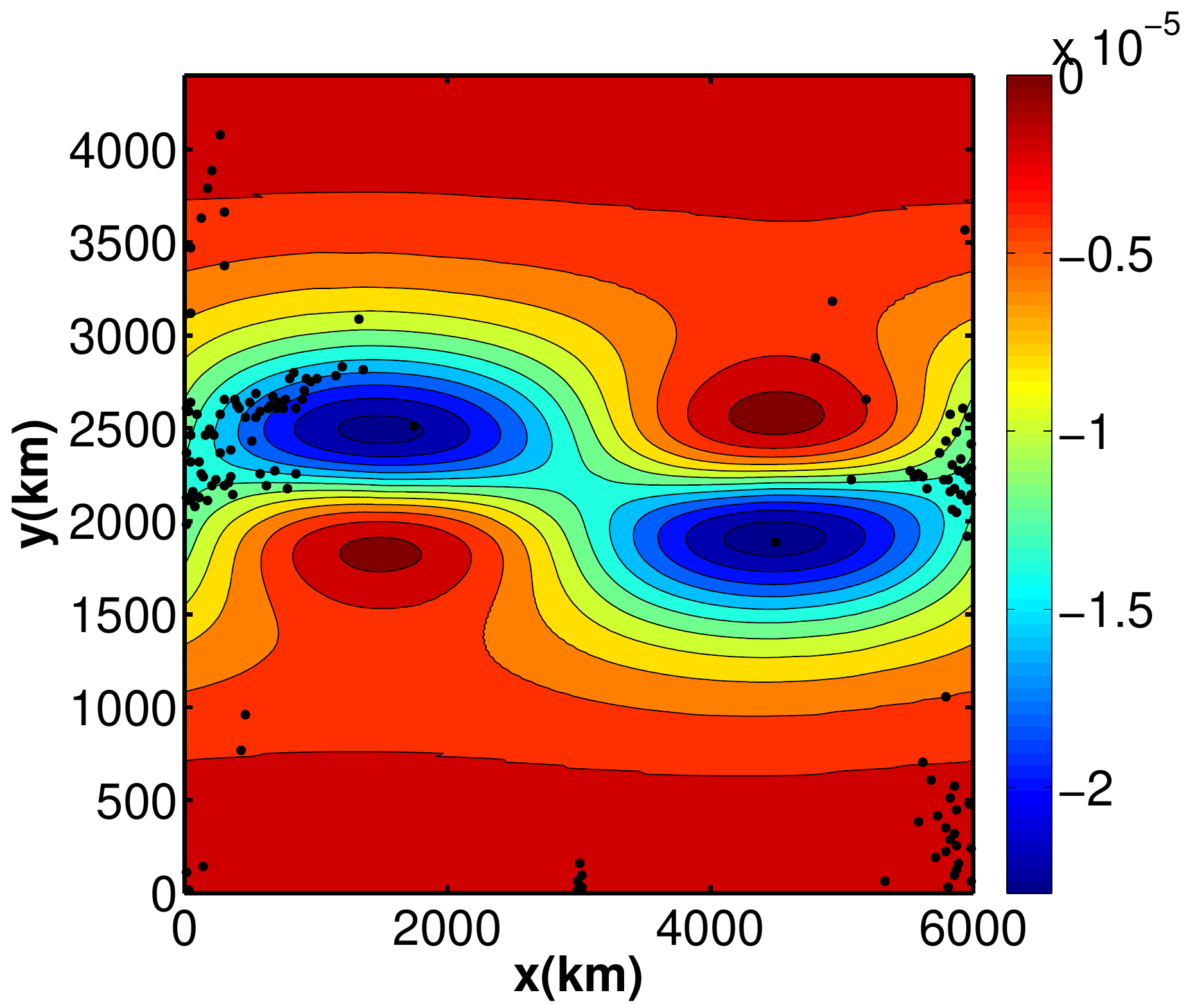}}\label{Fig::6d}
  \subfigure[Nonlinear term $F_{31}$]{\includegraphics[scale=0.222]{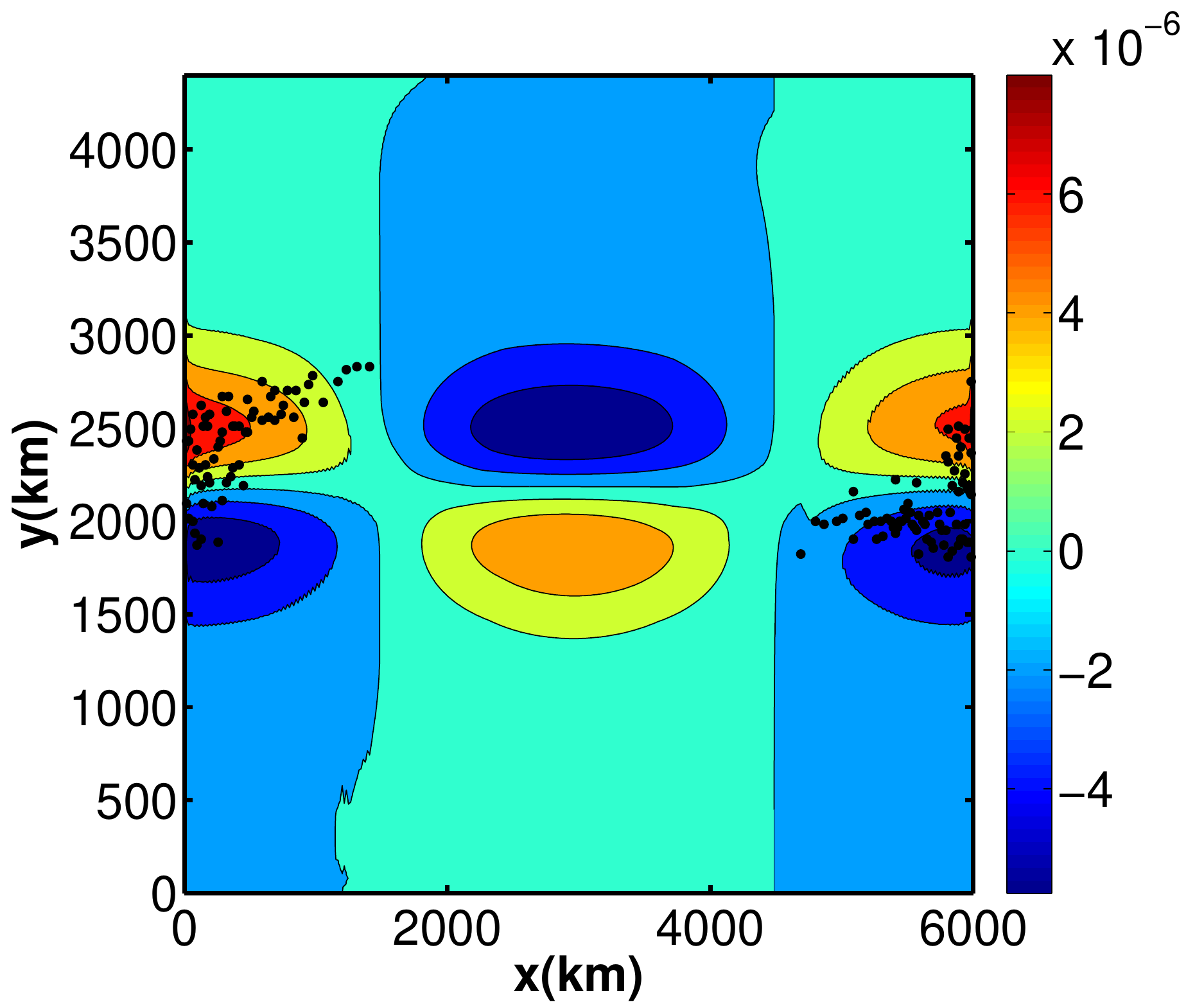}}\label{Fig::6e}
  \subfigure[Nonlinear term $F_{32}$]{\includegraphics[scale=0.222]{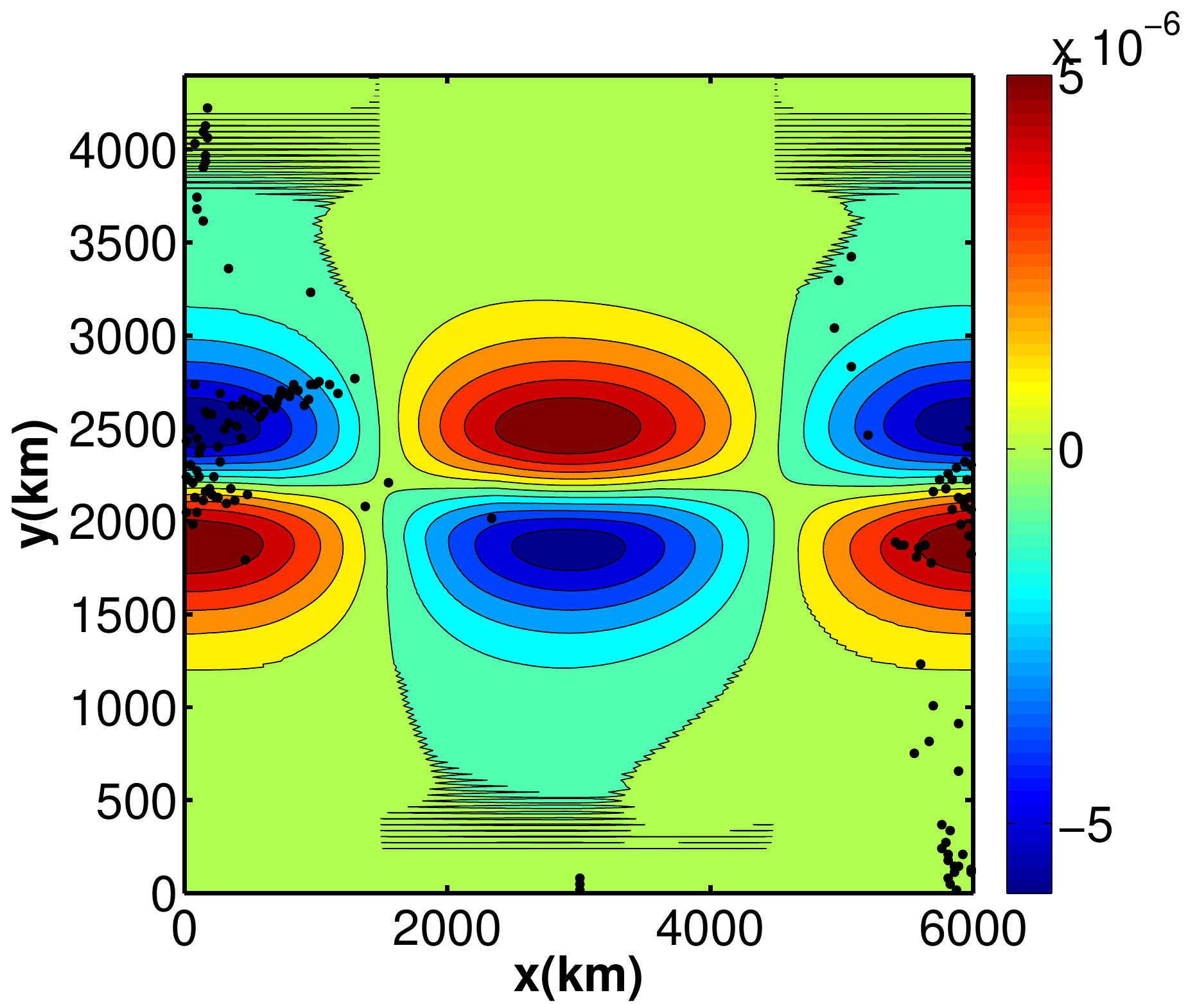}}\label{Fig::6f}
  \caption{100 DEIM interpolation points corresponding to all nonlinear terms  in the SWE model for time integration window of $3h$. The background consists in isolines of the maximum values of the nonlinear terms over time. Most of the points are concentrated in the region with larger errors depicted in Figure \ref{Fig::8}. \label{Fig::6}}
\end{figure}

Next we apply DEIM algorithm and calculate the interpolation points to improve the efficiency of the
standard POD approximation and to achieve a complexity reduction of the nonlinear terms with a complexity
proportional to the number of reduced variables, as in the case of tensorial POD. Figures \ref{Fig::5},\ref{Fig::6} illustrate the distribution of the first $100$ spatial points selected by the DEIM algorithm together with the isolines of the nonlinear terms statistics. Each of these statistics contain in every space location the maximum values of the corresponding nonlinear term over time. Maximum is preferred instead of time averaging since a better correlation between location of DEIM points and physical structures was observed in the former case.

%

However, in most cases, the spatial positions of the interpolation points don't follow the nonlinear statistics structures. This is more visible in Figure \ref{Fig::6} for $t_{\rm f}=3$h. The exceptions are $F_{12}$ and $F_{21}$ (Figures 5b,c), nonlinear terms depending only on velocity components, where DEIM interpolation points target better the underlying physical structures. This proves that DEIM algorithm doesn't particularly take into account the physical structures of the nonlinear terms but search (in a greedy manner) to minimize the error (residual) between each column of the input basis (POD basis of the nonlinear term snapshots) and its proposed low-rank approximations \citet[p.16]{Stefanescu2013}.

\begin{figure}[ht]
  \centering
  \subfigure[ $u_\textnormal{\sc pod}-u_\textnormal{\sc full}$ ]{\includegraphics[scale=0.222]{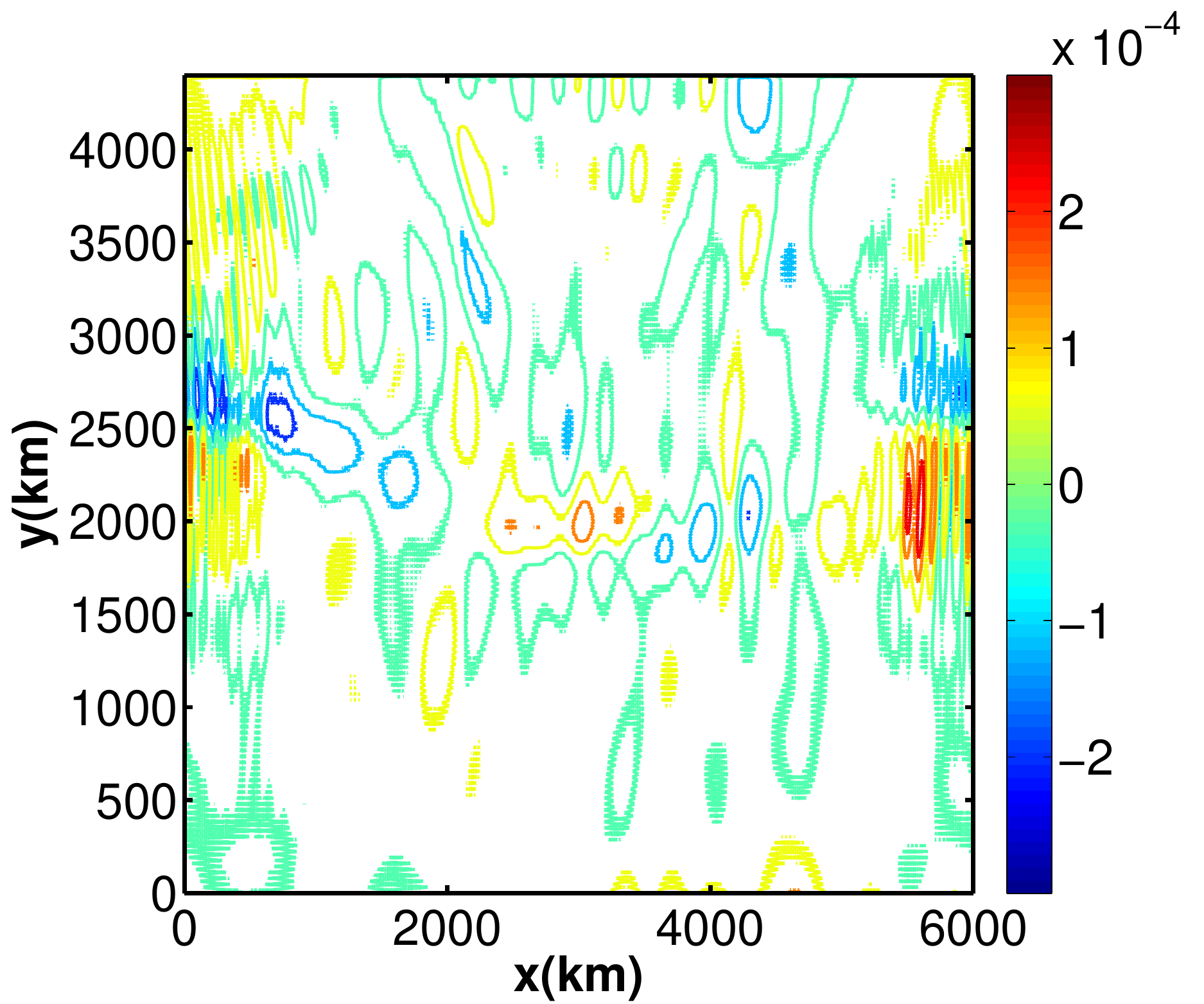}}\label{Fig::7a}
  \subfigure[ $v_\textnormal{\sc pod}-v_\textnormal{\sc full}$]{\includegraphics[scale=0.222]{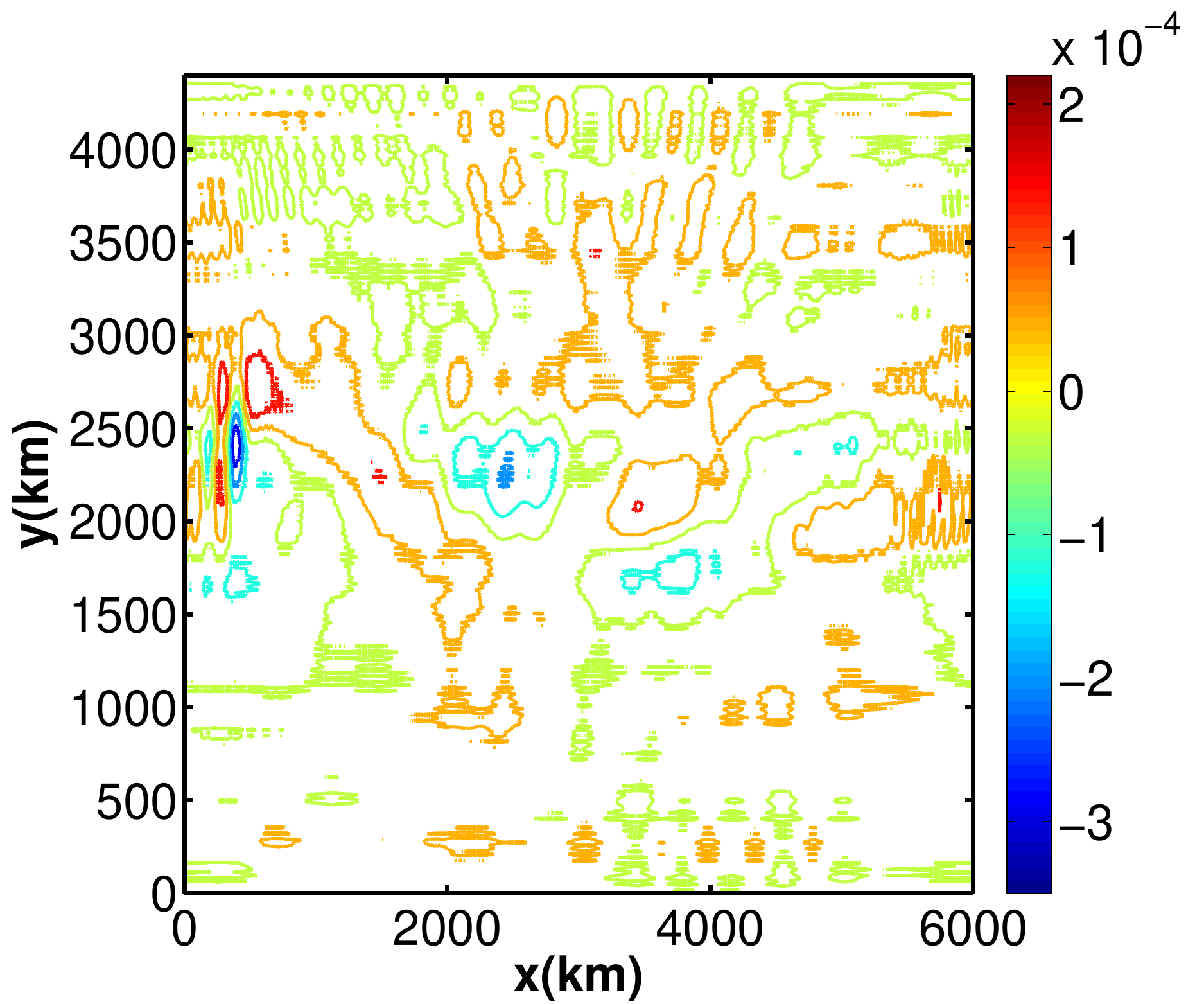}}\label{Fig::7b}
  \subfigure[ $\phi_\textnormal{\sc pod}-\phi_\textnormal{\sc full}$]{\includegraphics[scale=0.222]{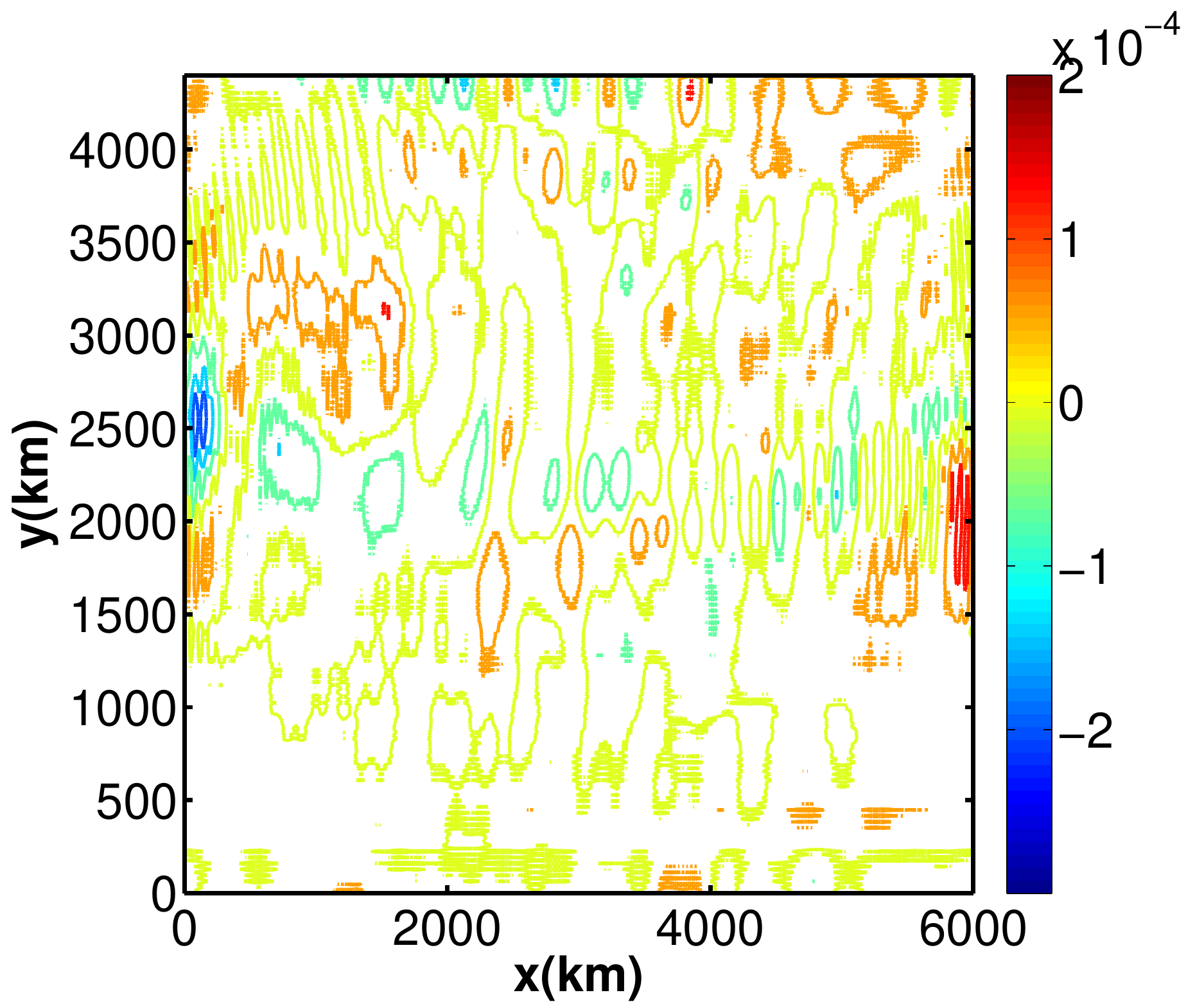}}\label{Fig::7c}
  \subfigure[ $u_\textnormal{\sc tpod}-u_\textnormal{\sc full}$]{\includegraphics[scale=0.222]{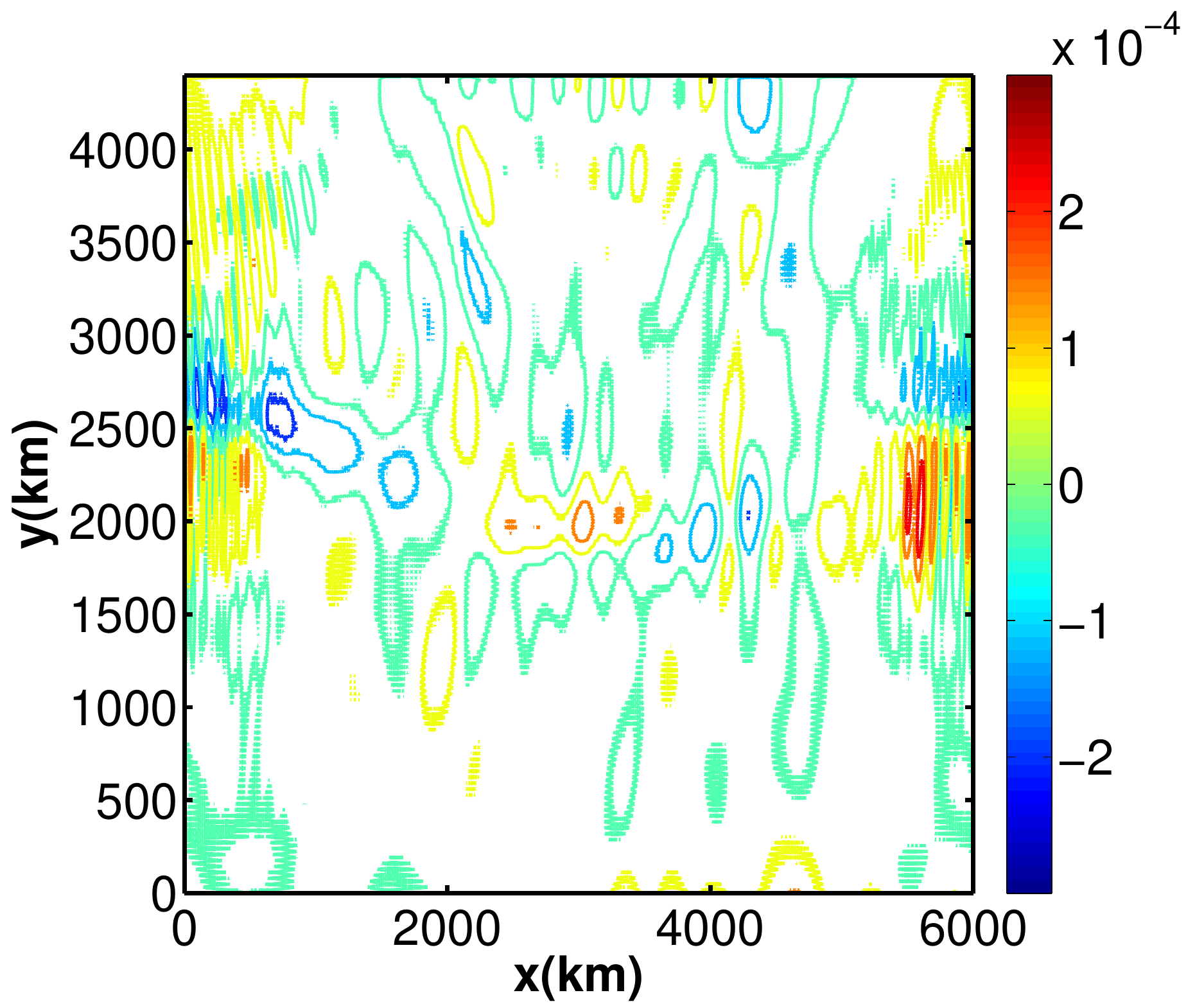}}\label{Fig::7d}
  \subfigure[ $v_\textnormal{\sc tpod}-v_\textnormal{\sc full}$]{\includegraphics[scale=0.222]{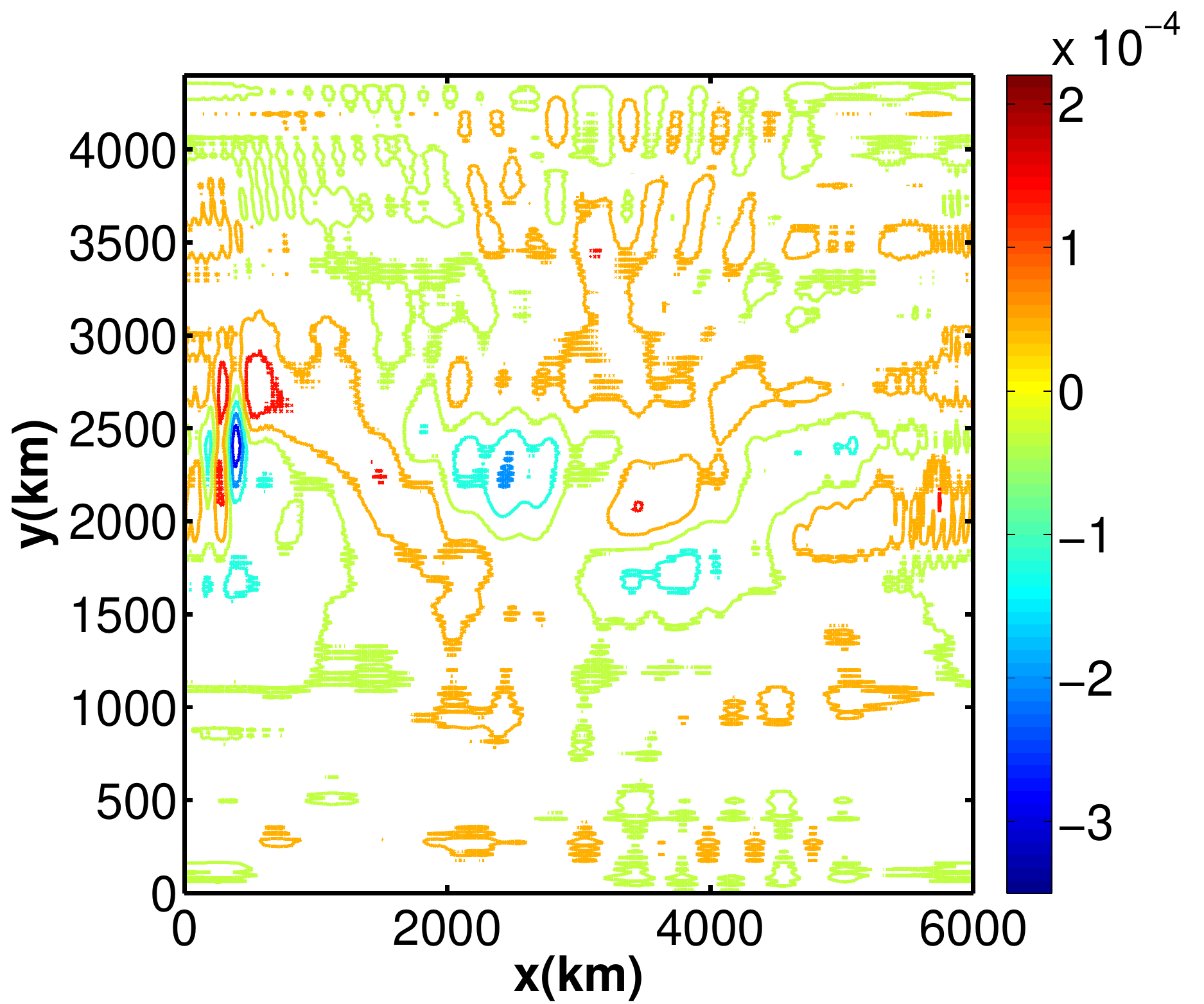}}\label{Fig::7e}
  \subfigure[ $\phi_\textnormal{\sc tpod}-\phi_\textnormal{\sc full}$]{\includegraphics[scale=0.222]{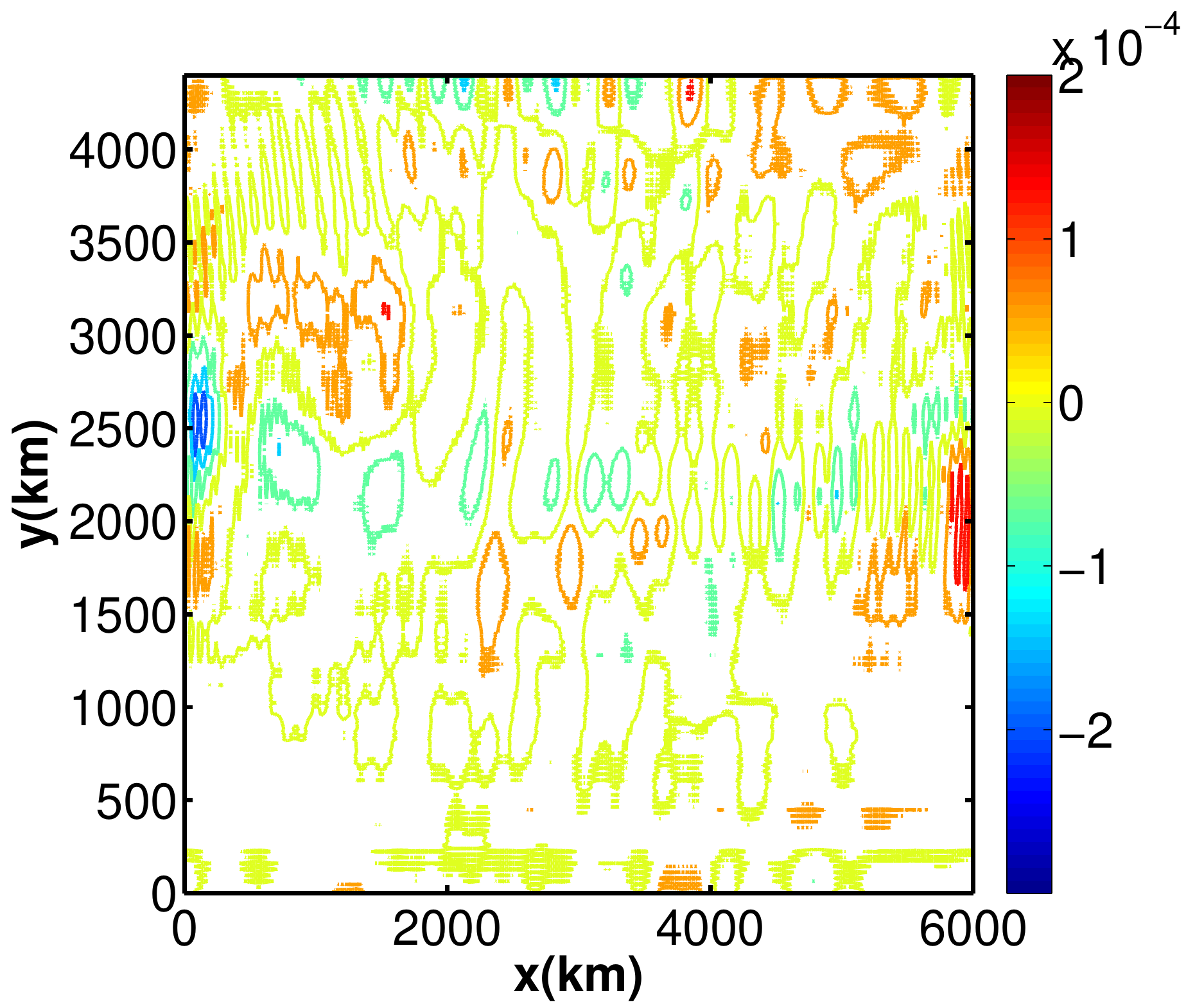}}\label{Fig::7f}
  \subfigure[   $u_\textnormal{\sc pod/deim}-u_\textnormal{\sc full}$]{\includegraphics[scale=0.222]{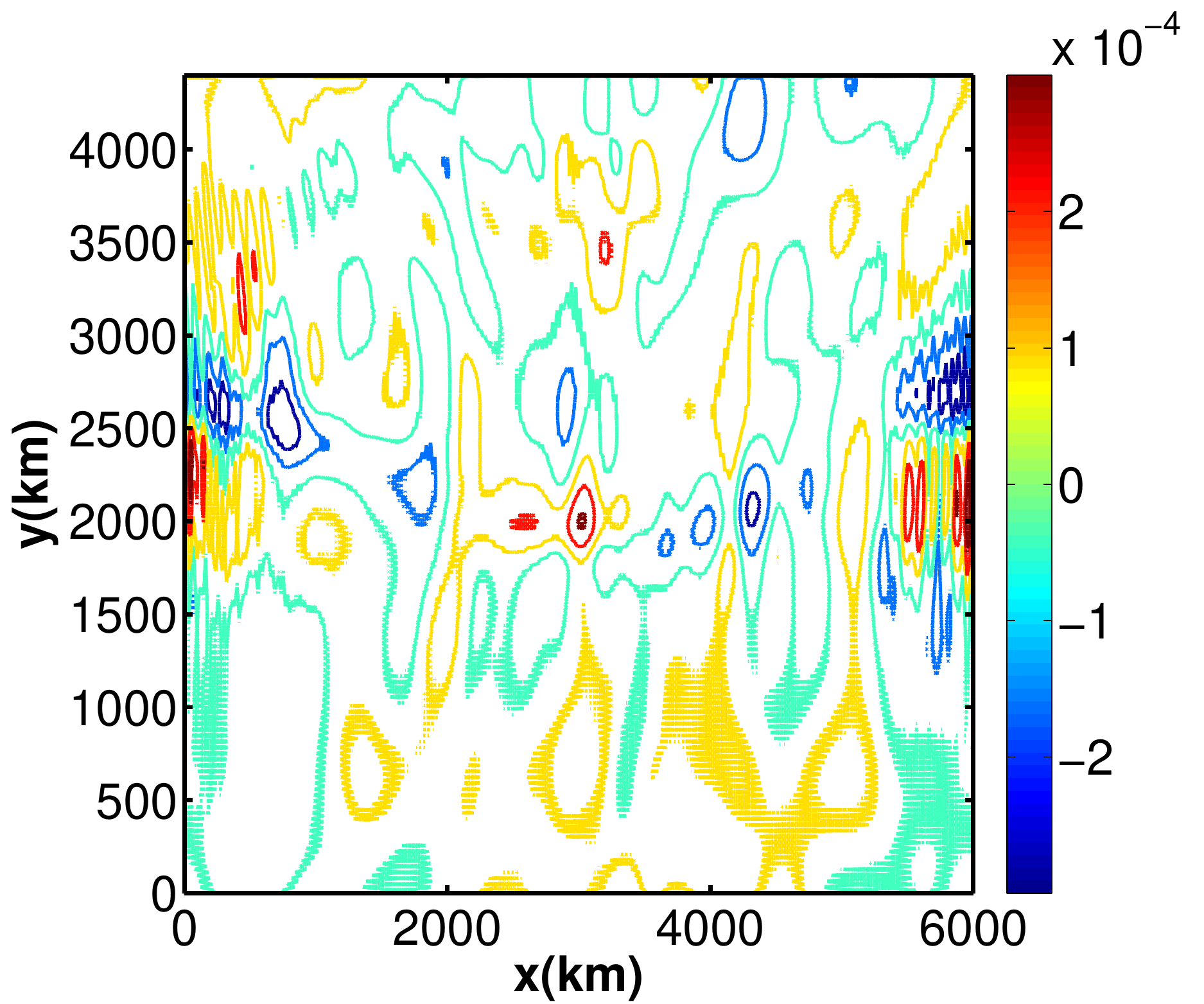}}\label{Fig::7g}
  \subfigure[   $v_\textnormal{\sc pod/deim}-v_\textnormal{\sc full}$]{\includegraphics[scale=0.222]{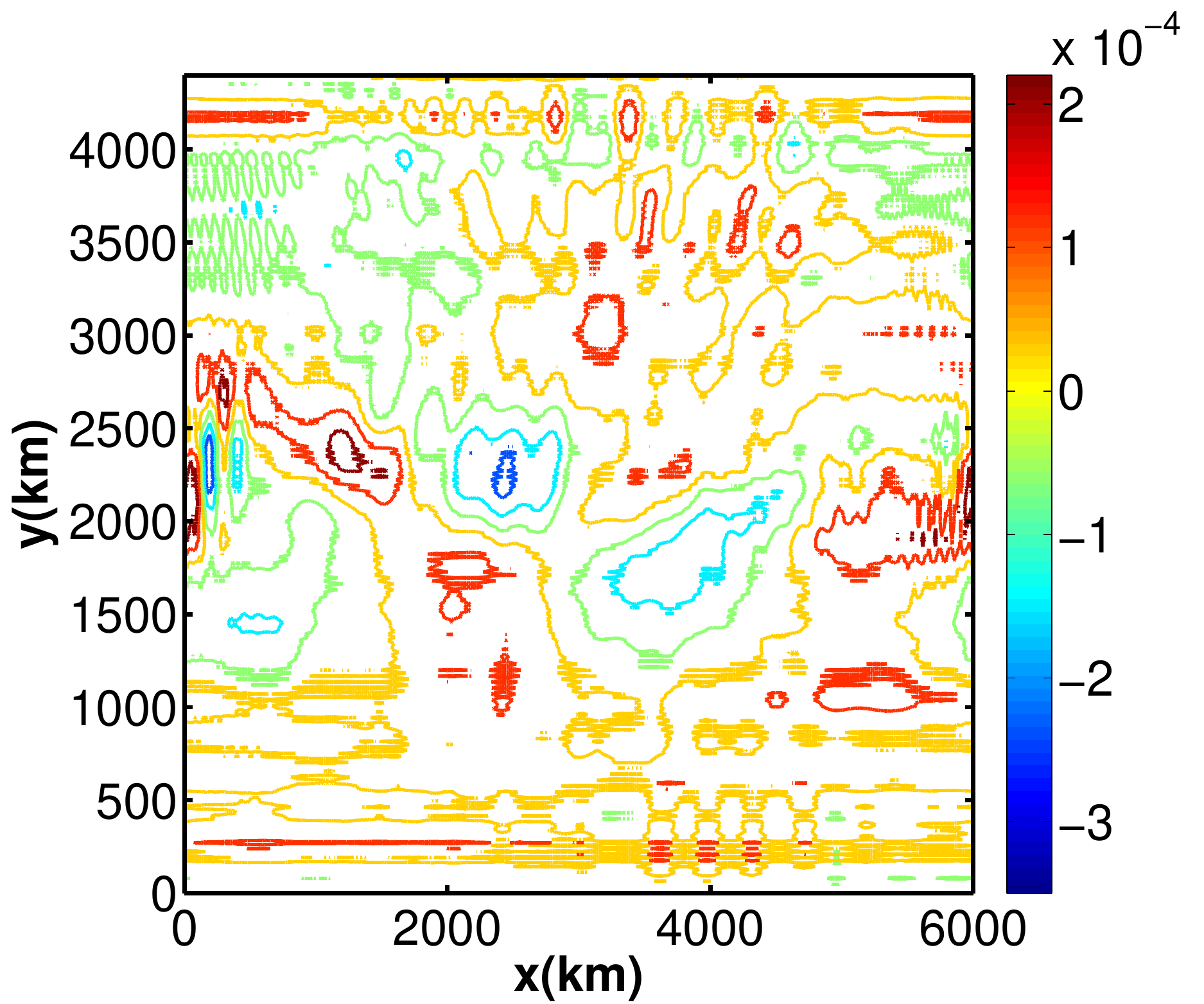}}\label{Fig::7h}
  \subfigure[   $\phi_\textnormal{\sc pod/deim}-\phi_\textnormal{\sc full}$]{\includegraphics[scale=0.222]{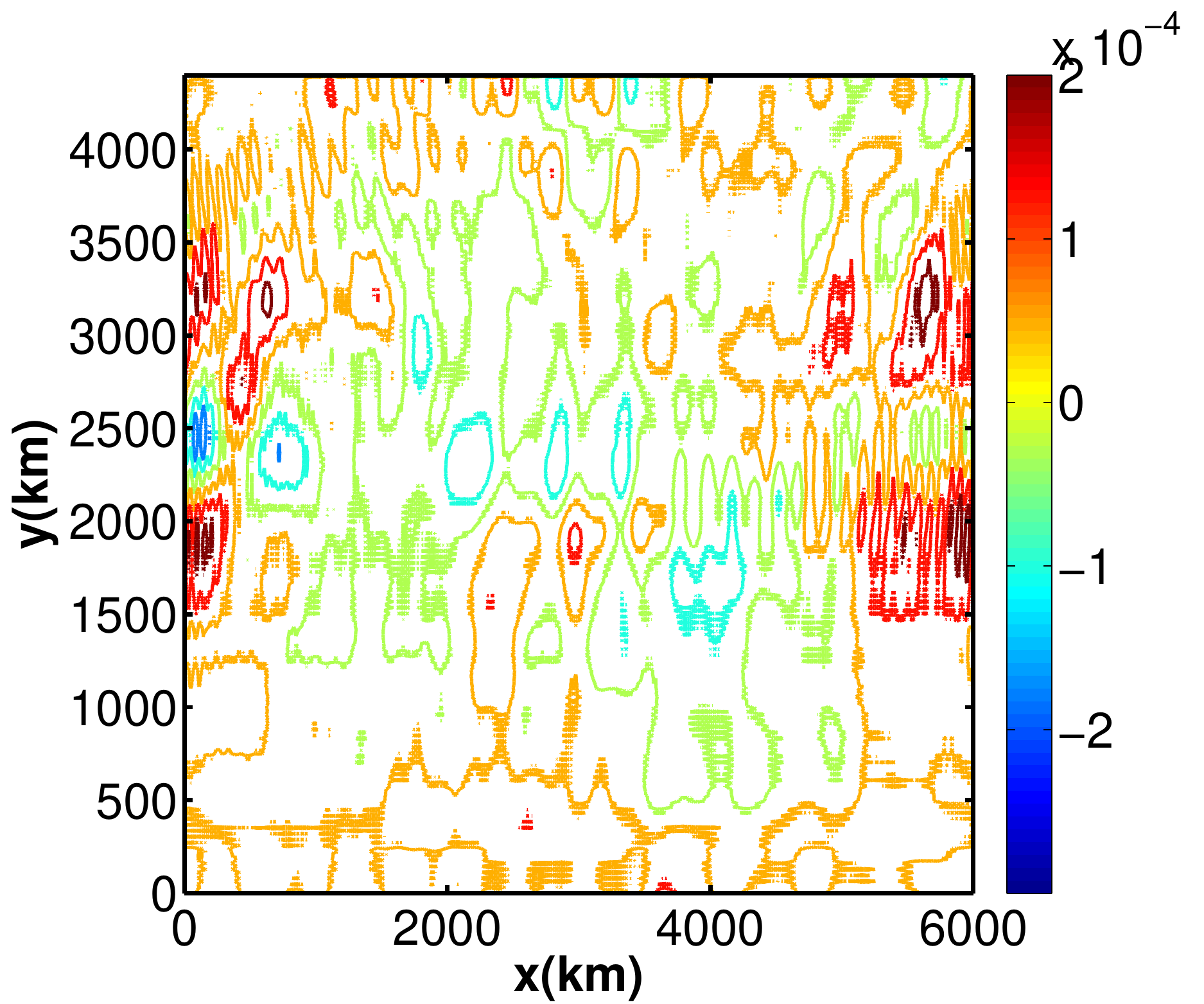}}\label{Fig::7i}
  \caption{Absolute errors between standard POD, tensorial POD and POD/DEIM solutions and the full trajectories
at t = 24h ($\Delta t$ = 960s). The number of DEIM points was taken 180\label{Fig::7}}
\end{figure}

\begin{figure}[ht]
  \centering
  \subfigure[ $u_\textnormal{\sc pod}-u_\textnormal{\sc full}$ ]{\includegraphics[scale=0.222]{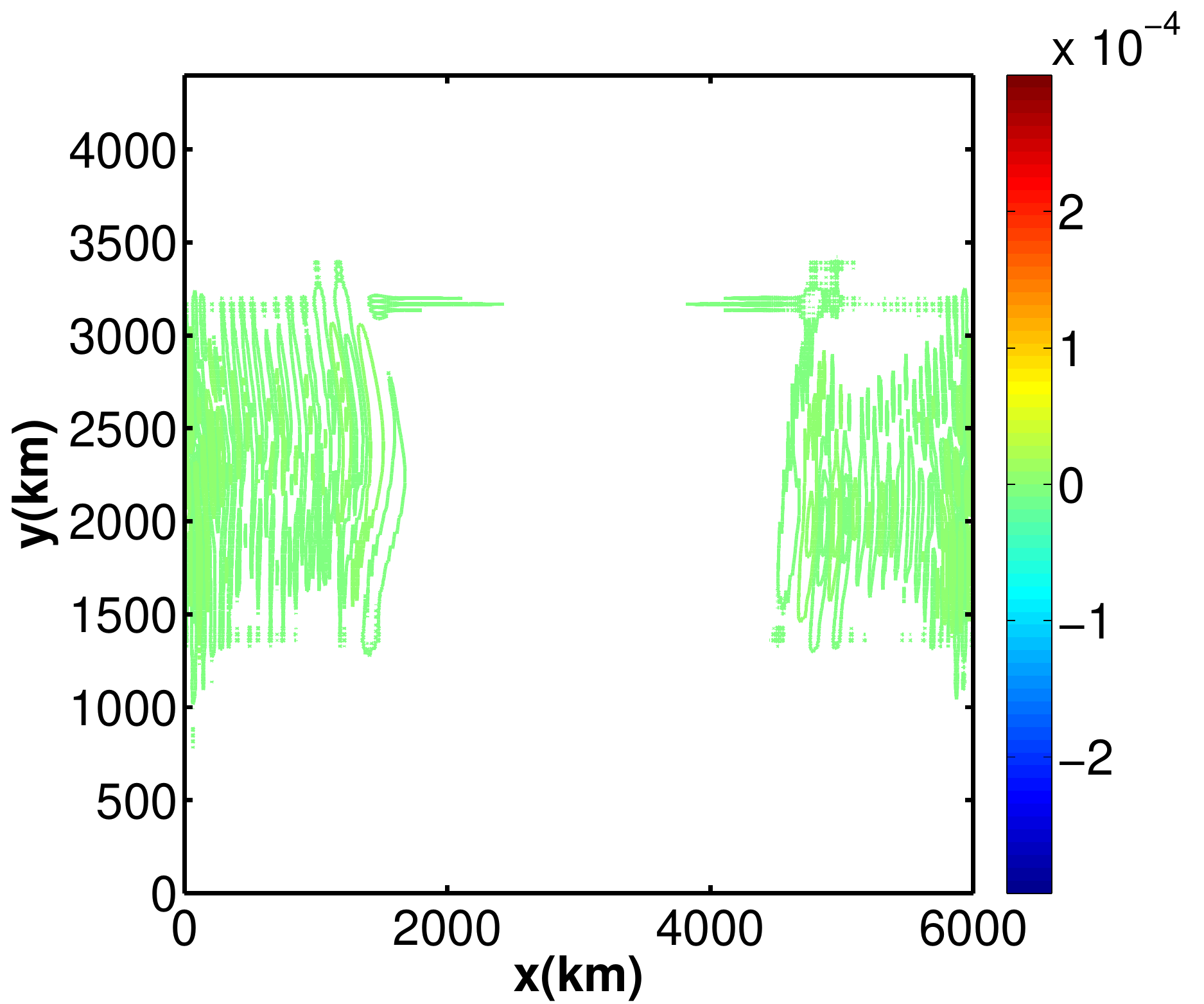}}\label{Fig::8a}
  \subfigure[ $v_\textnormal{\sc pod}-v_\textnormal{\sc full}$]{\includegraphics[scale=0.222]{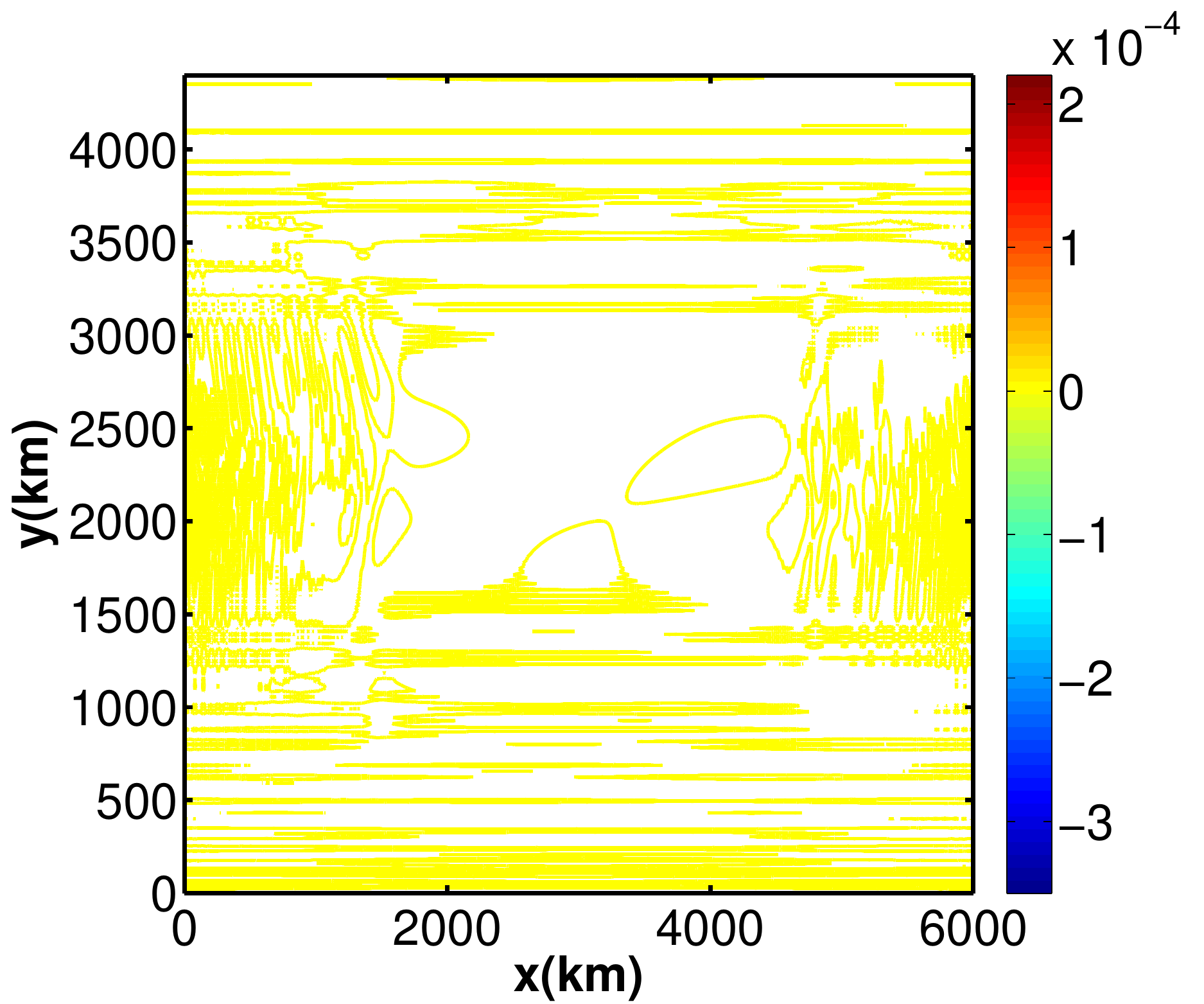}}\label{Fig::8b}
  \subfigure[ $\phi_\textnormal{\sc pod}-\phi_\textnormal{\sc full}$]{\includegraphics[scale=0.222]{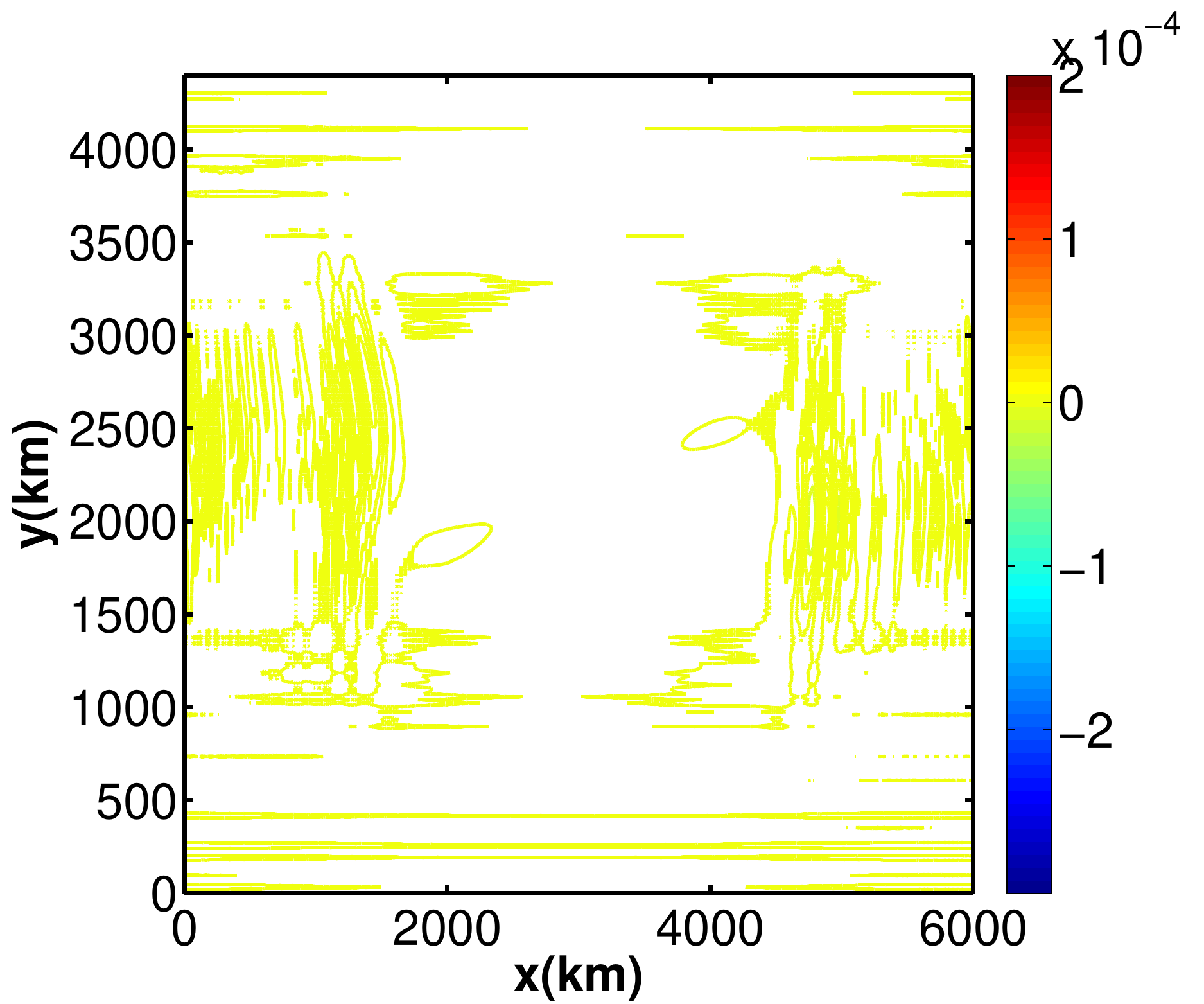}}\label{Fig::8c}
  \subfigure[ $u_\textnormal{\sc tpod}-u_\textnormal{\sc full}$]{\includegraphics[scale=0.222]{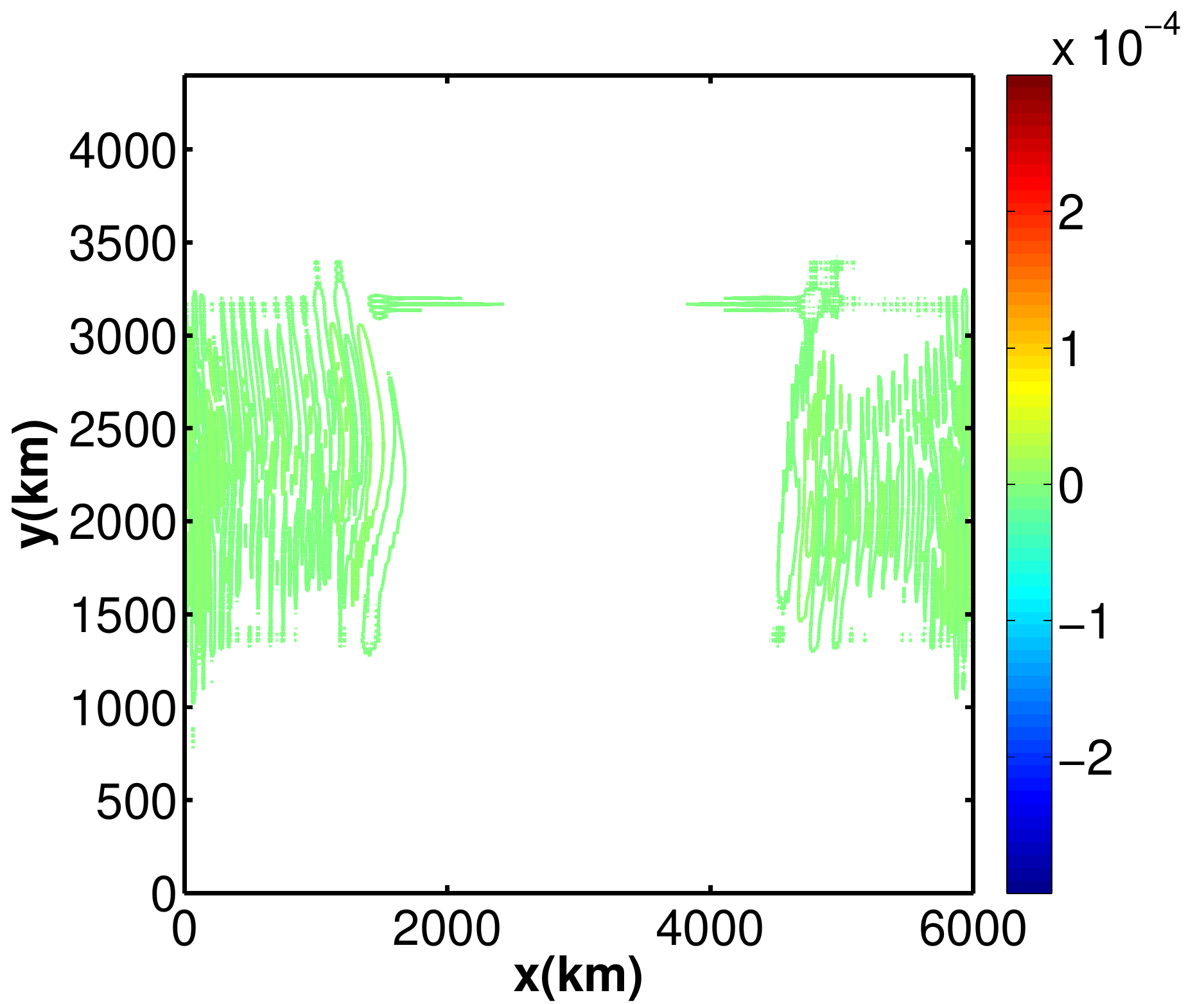}}\label{Fig::8d}
  \subfigure[ $v_\textnormal{\sc tpod}-v_\textnormal{\sc full}$]{\includegraphics[scale=0.222]{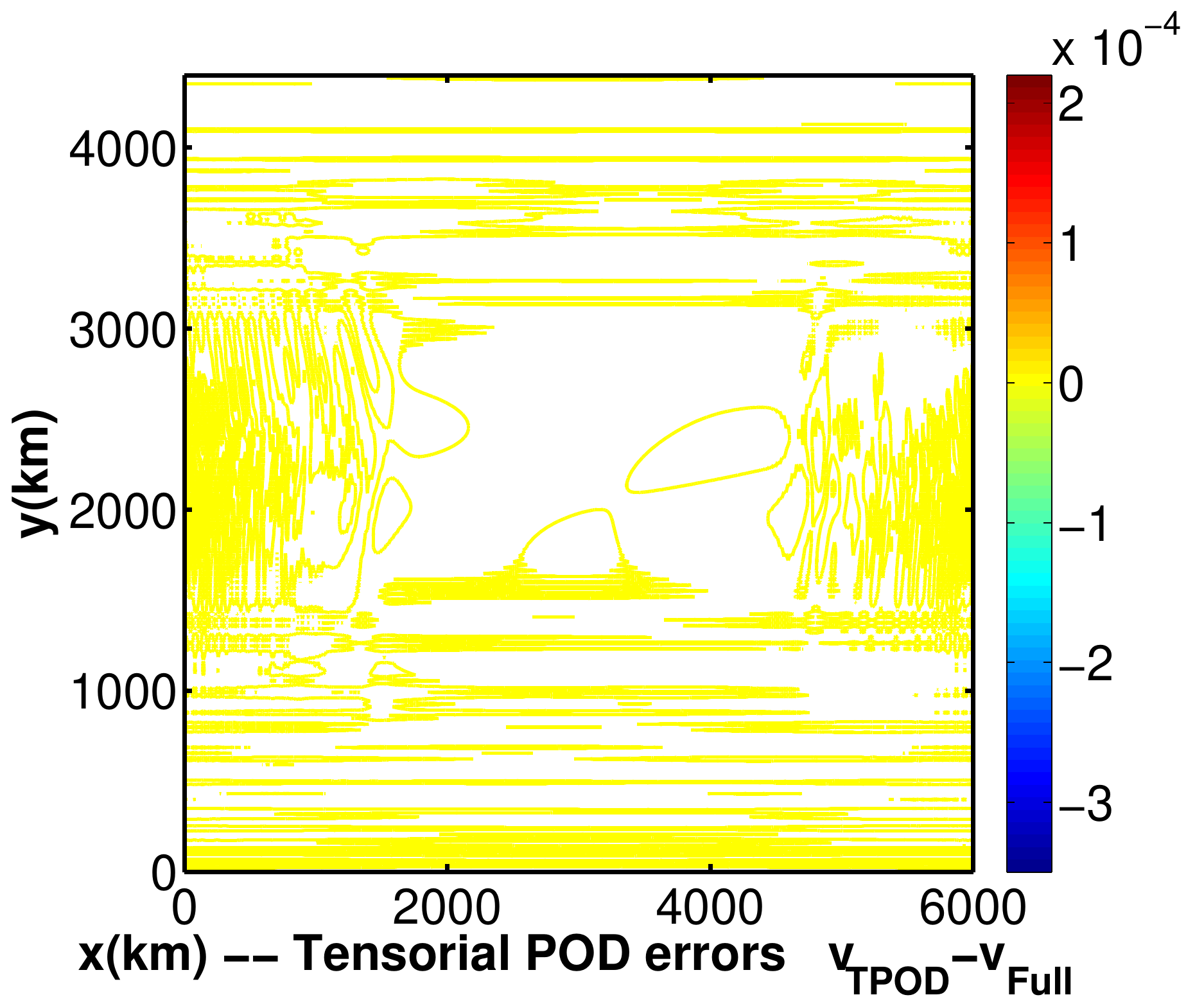}}\label{Fig::8e}
  \subfigure[ $\phi_\textnormal{\sc tpod}-\phi_\textnormal{\sc full}$]{\includegraphics[scale=0.222]{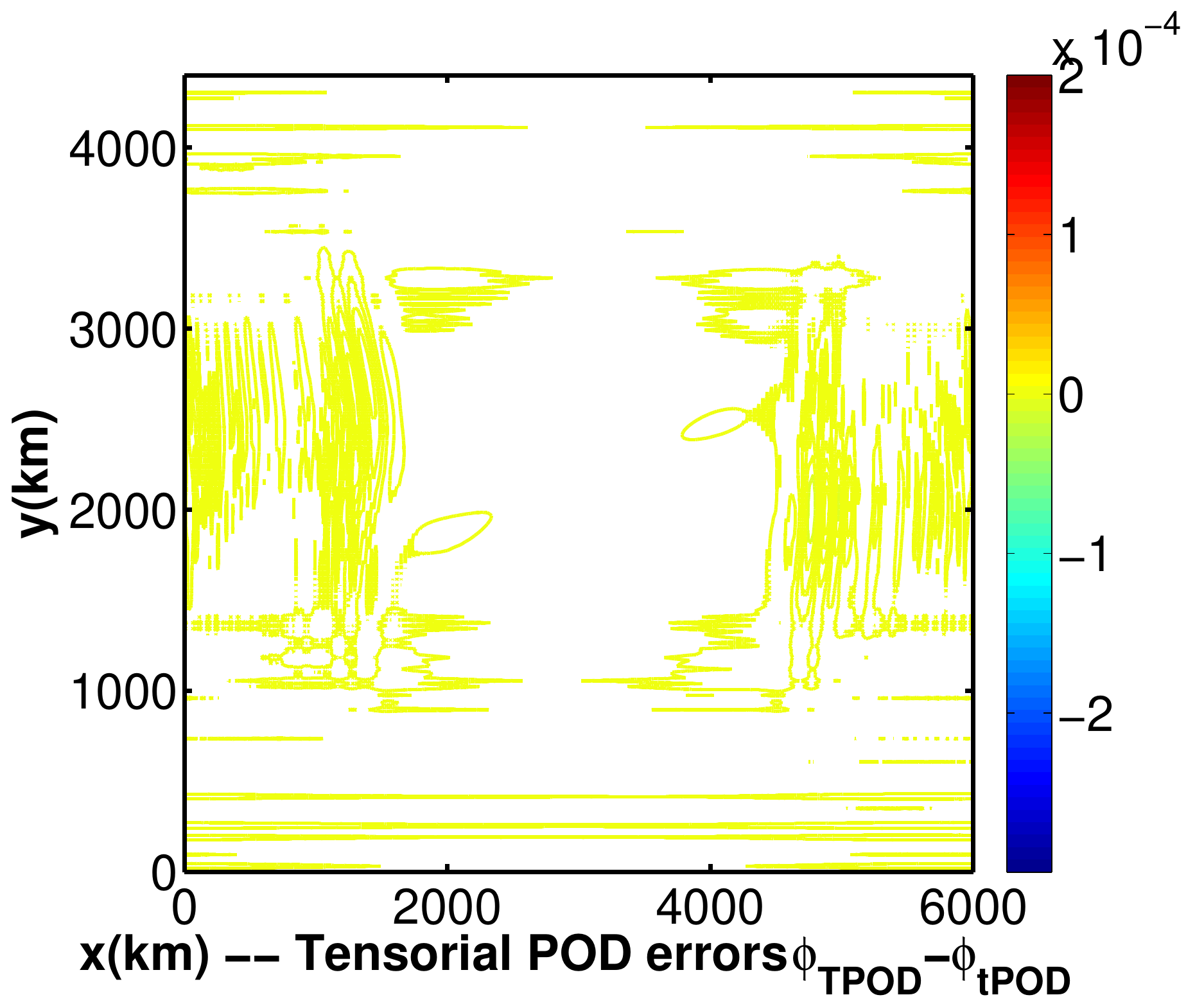}}\label{Fig::8f}
  \subfigure[   $u_\textnormal{\sc pod/deim}-u_\textnormal{\sc full}$]{\includegraphics[scale=0.222]{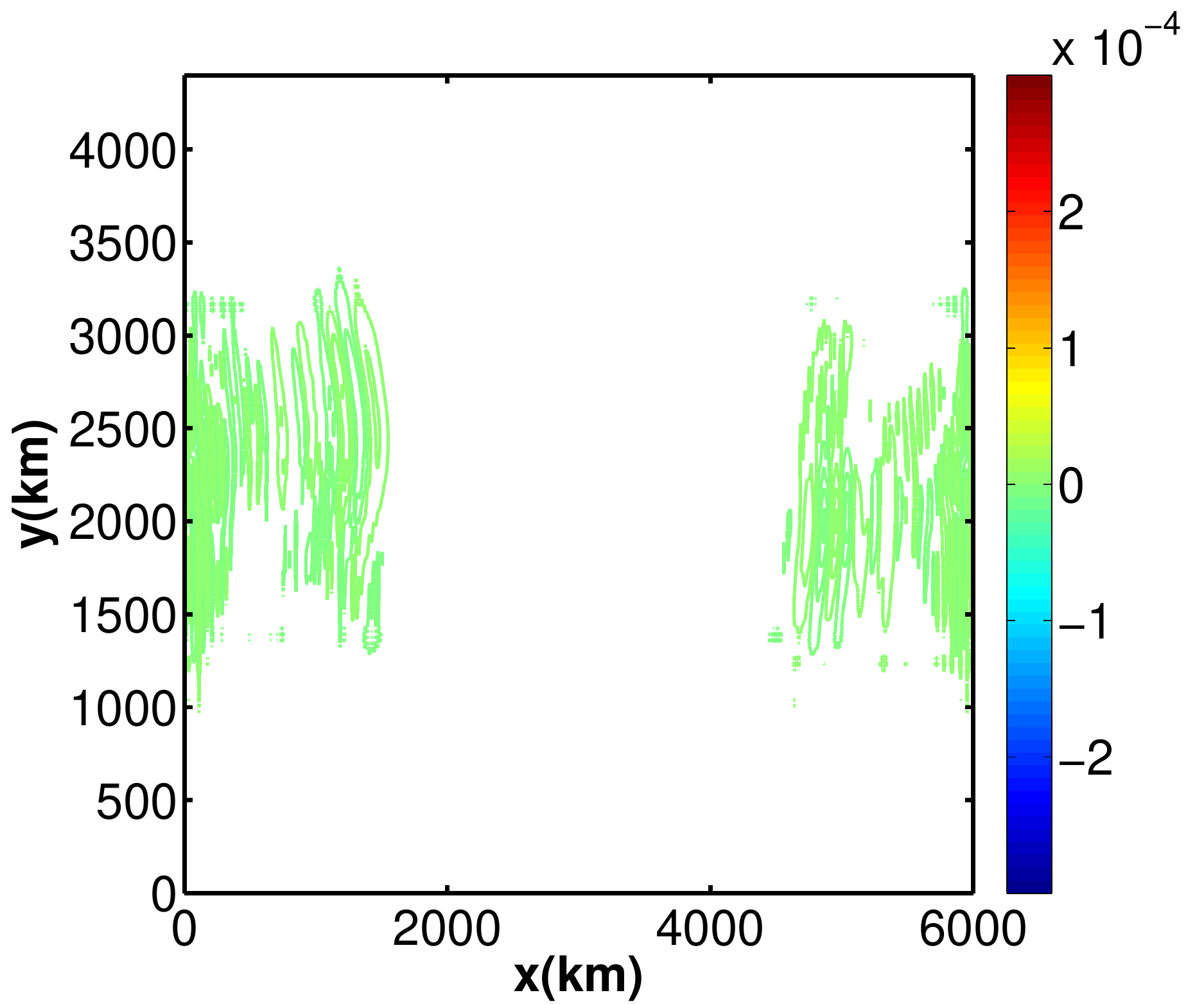}}\label{Fig::8g}
  \subfigure[   $v_\textnormal{\sc pod/deim}-v_\textnormal{\sc full}$]{\includegraphics[scale=0.222]{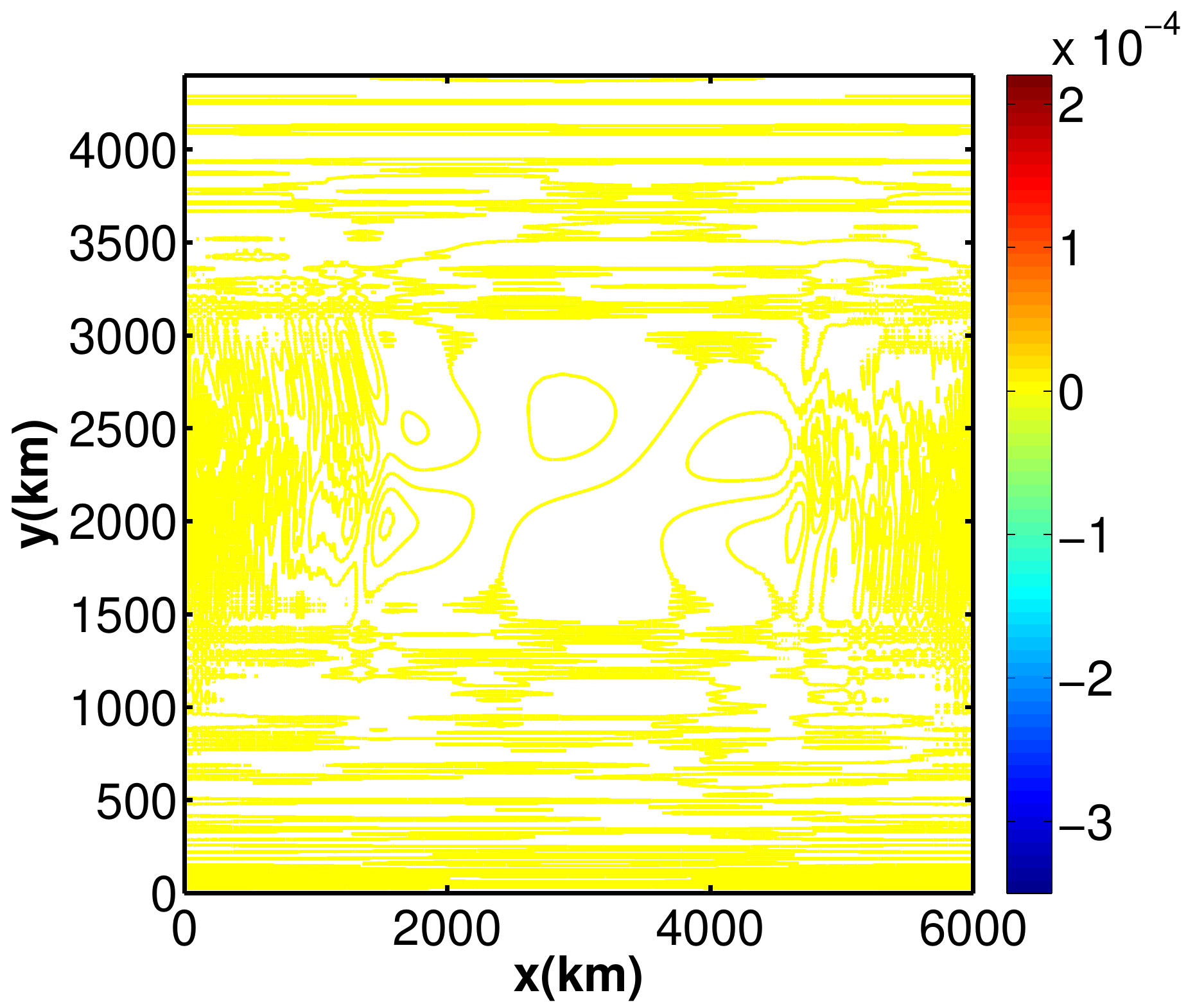}}\label{Fig::8h}
  \subfigure[   $\phi_\textnormal{\sc pod/deim}-\phi_\textnormal{\sc full}$]{\includegraphics[scale=0.222]{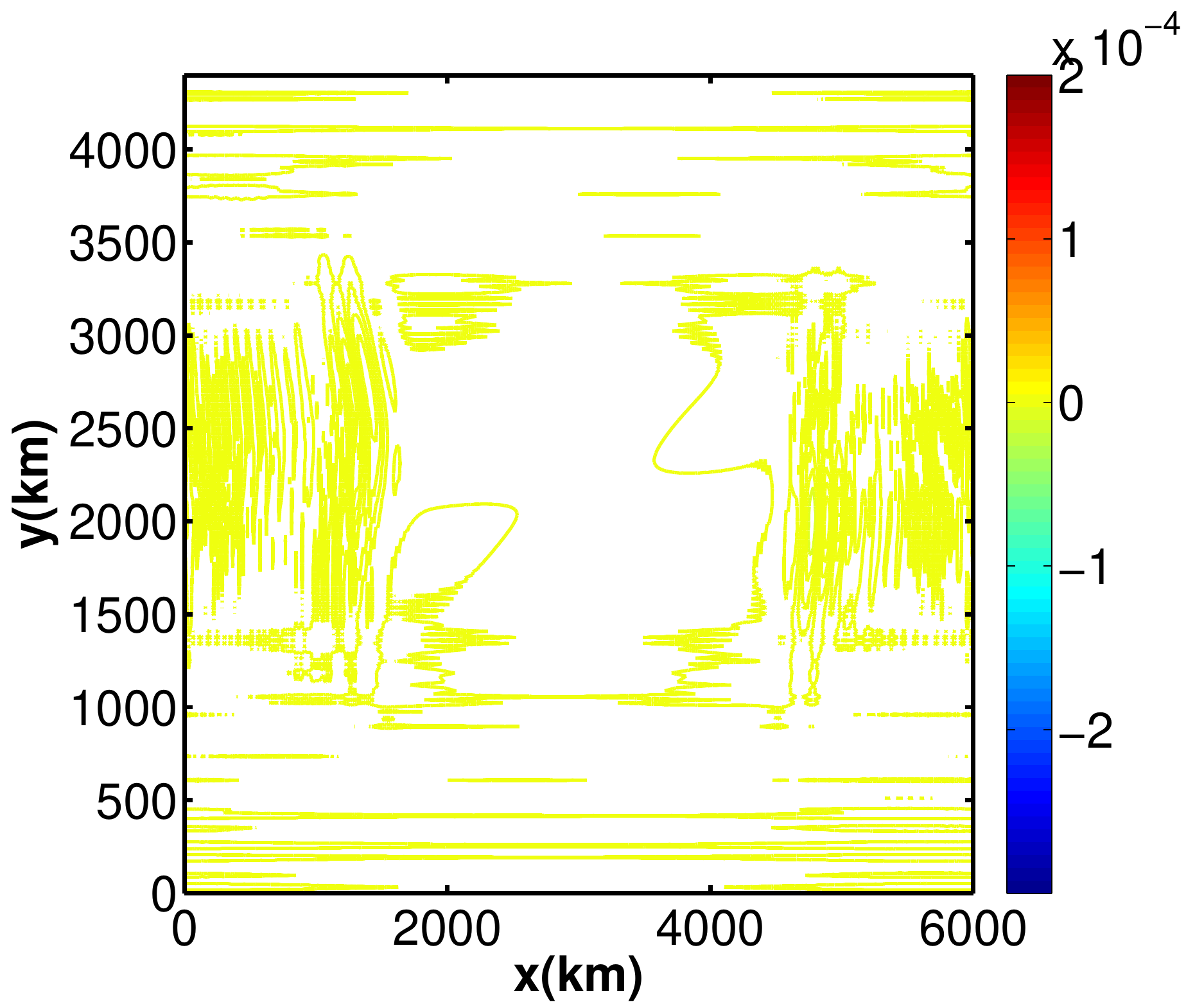}}\label{Fig::8i}
\caption{Absolute errors between standard POD, tensorial POD and POD/DEIM solutions and the full trajectories
at t = 3h ($\Delta t$ = 120s). The number of DEIM points was taken 180\label{Fig::8}}
\end{figure}

Figures \ref{Fig::7},\ref{Fig::8} depict the grid point absolute error of the standard POD, tensorial POD and POD/DEIM solutions with respect to the full solutions. For POD/DEIM reduced order model we use $180$ interpolation points. The magnitude of the errors are similar for each of the method proposed in this study. Moreover, we observed that error isolines distribution in Figure \ref{Fig::8} is well correlated with the location of interpolation points illustrated in Figure \ref{Fig::6} underlying the empirical characteristics of DEIM.

%

%

In addition, we propose two metrics to quantify the accuracy level of standard POD, tensorial POD and standard POD/DEIM approaches. First, we use the following norm
\[
 \frac{1}{N_t} \sum_{i=1}^{t_{\rm f}}\frac{||w^\textnormal{\sc full}(:,t_i)-w^\textnormal{\sc rom}(:,t_i)||_2}{||w^\textnormal{\sc full}(:,t_i)||_2}
 \]
  \[
 \frac{1}{N_t} \sum_{i=1}^{N_t}\frac{\left\Vert w^\textnormal{\sc full}(:,t_i)-w^\textnormal{\sc rom}(:,t_i) \right\Vert_2}{\left\Vert w^\textnormal{\sc full}(:,t_i)  \right\Vert_2}
 \]
$i=1,2,..,t_{\rm f}$ and calculate the relative errors for all three variables of SWE model $w=(u,v,\phi)$. The results are presented in Table \ref{table3}. We perform numerical experiments using two choices for number of DEIM points $70$ and $180$. For $24$h tests we notice that more than $70$ number of DEIM points are needed for convergence of quasi-Newton method for POD/DEIM SWE scheme explaining the absence of numerical results in Table \ref{table3} (left part).

\begingroup
\begin{table}[h]
\centerline{
\scalebox{0.9}{
\begin{tabular}{|c|c|c|c|}\hline
 &  Standard & Tensorial & POD/DEIM \\
 &  POD        & POD & m=180 \\ \hline
$u$ & 1.276e-3 & 1.276e-3 & 1.622e-3\\ \hline
$v$ & 3.426e-3 & 3.426e-3 & 4.639e-3\\ \hline
$\phi$ & 2.110e-5 & 2.110e-5 & 2.489e-5\\ \hline
\end{tabular}}
\scalebox{0.9}{
\begin{tabular}{|c|c|c|c|c|}\hline
 &  Standard & Tensorial & POD/DEIM & POD/DEIM \\
 & POD & POD & m = 180 & m = 70 \\ \hline
$u$ & 7.711e-6 & 7.711e-6 & 7.965e-6 & 9.301e-6\\ \hline
$v$ & 1.665e-5 & 1.666e-5 & 1.73e-5 & 1.975e-5\\ \hline
$\phi$ & 1.389e-7 & 1.389e-7 & 1.426e-7 & 1.483e-7 \\ \hline
\end{tabular}}}
\caption{\label{table3}Relative errors for each of the model variables for $t_f=24$h (left) and $t_f=3$h (right). The POD bases dimensions were taken $50$. For $24$h experiments we display only the results for $180$ number of DEIM points while in the case of $3$h time integration window tests with $180$ and $70$ numbers of DEIM points are shown.}
\end{table} %
\endgroup%

Root mean square error is also employed to compare the reduced order models. Table \ref{table4} and \ref{table5} show the RMSE for final times together with the CPU times of the on-line stage of ROMs.

\begingroup
\begin{table}[h]
\centerline{
\scalebox{0.9}{
\begin{tabular}{|c|c|c|c|c|}\hline 
 & Full ADI SWE & Standard POD & Tensorial POD & POD/DEIM m=180 \\ \hline
CPU time & 1813.992s & 191.785 & 2.491& 1.046\\ \hline
$u$ & - & 9.095e-3 & 9.095e-3& 1.555e-2\\ \hline
$v$ & - & 8.812e-3 & 8.812e-3& 1.348e-2\\ \hline
$\phi$ & - & 6.987e-3e & 6.987e-3& 1.13e-2\\ \hline
\end{tabular}}}%
\caption{\label{table4}CPU time gains and the root mean square errors for each of the model variables at $t_f=24$h. Number of POD modes was $k=50$ and we choose $180$ number of DEIM points.}
\end{table} %
\endgroup%
\begingroup
\begin{table}[h]
\centerline{
\scalebox{0.87}{
\begin{tabular}{|c|c|c|c|c|c|}\hline 
 & Full ADI SWE & Standard POD & Tensorial POD & POD/DEIM m=180 & POD/DEIM m=70 \\ \hline
CPU time  & 950.0314s & 161.907 & 2.125& 0.642 & 0.359\\ \hline
$u$ & - & 5.358e-5 & 5.358e-5& 5.646e-5 & 7.453e-5\\ \hline
$v$ & - & 2.728e-5 & 2.728e-5& 3.418e-5 & 4.233e-5\\ \hline
$\phi$ & - & 8.505e-5e & 8.505e-5& 8.762e-5 & 9.212e-5\\ \hline
\end{tabular}}}
\caption{\label{table5}CPU time gains and the root mean square errors for each of the model variables at $t_f=3$h for a $3$h time integration window. Number of POD modes was $k=50$ and two tests with different number of DEIM points $m = 180,~70$ were simulated.}
\end{table} %
\endgroup%

Thus, for $103,776$ spatial points, tensorial POD method reduces the computational complexity of the nonlinear terms in comparison with the POD ADI SWE model and overall decreases the computational time with a factor of $77\times$ for $24$h time integration and $76\times$ for a $3$h time window integration. POD/DEIM  outperforms standard POD being $450\times$ and $250 \times$ time faster for $70$ and $180$ DEIM interpolation points and a time integration window of $3h$. For $t_{\rm f}=24$h and $m=180$, POD/DEIM SWE model is $183\times$ time faster than standard POD SWE model. In terms of CPU time, tensor POD SWE model is only $2.38\times$ slower than POD/DEIM SWE model for $m=180$ and $t_{\rm f}=24$h while for $t_{\rm f}=3$h the new tensorial POD scheme is $5.9\times$ and $3.3 \times$ less efficient than POD/DEIM SWE model for $m=70$ and $m=180$. This suggests that operations like Jacobian computations and its LU decomposition required by both reduced order approaches weight more in the overall CPU time cost since the quadratic nonlinear complexity of the tensorial POD requires $374,950$ floating-point operations and POD/DEIM only $27,180$ ($m=180$, see Section $3$). Given that the implementation effort is much reduced, in the cases of models depending only on quadratic nonlinearities, the tensorial POD poses the appropriate characteristics of a reduced order method and represent a solid alternative to the POD/DEIM approach.

For cubical nonlinearities and larger, tensorial POD loses its ability to deliver fast calculations (see Table \ref{table1}), thus the POD/DEIM should be employed.

In our case the Jacobians are calculated analytically and its computations depend only on the reduced space dimension $k$. However, more gain can be obtain if DEIM would be applied to approximate the reduced Jacobians but this is subject of future research.

The computational savings and accuracy levels obtained by the ROMs studied in this paper depend on the number of POD modes and number of DEIM points. These numbers may be large in practice in order to capture well the full model dynamics. For exemple, in the case of a time window integration of $24$h, if someone would ask to increase the ROMs solutions accuracy with only one order of magnitude, the POD basis dimension must be at least larger than $100$ which will drastically compromise the time performances of ROMs methods. Elegant solutions to this problem were proposed by  \citet{Rapun_2010,Peherstorfer_2013} where local POD and local DEIM versions were proposed. Machine learning techniques such as $K$-means \citet{Lloyd_1957,MacQueen_1967,Steinhaus_1956}  can be used for both time and space partitioning. A recent study investigating cluster-based reduced order modeling was proposed by \citet{Kaiseretal2013}.

\begin{figure}[h]
  \centering
  \subfigure[On-line stage ] {\includegraphics[scale=0.35]{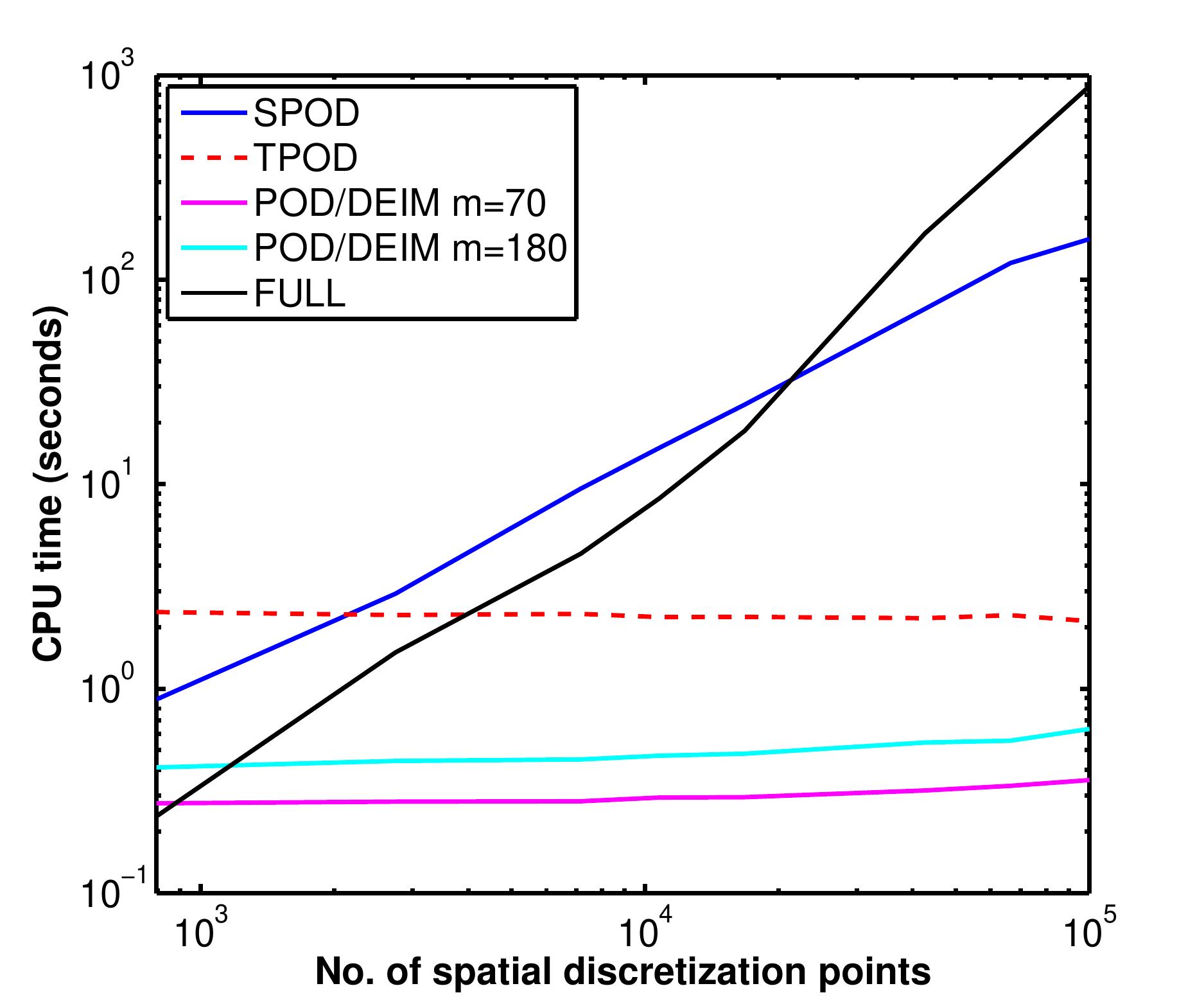}}
  \subfigure[Off-line stage ]{\includegraphics[scale=0.35]{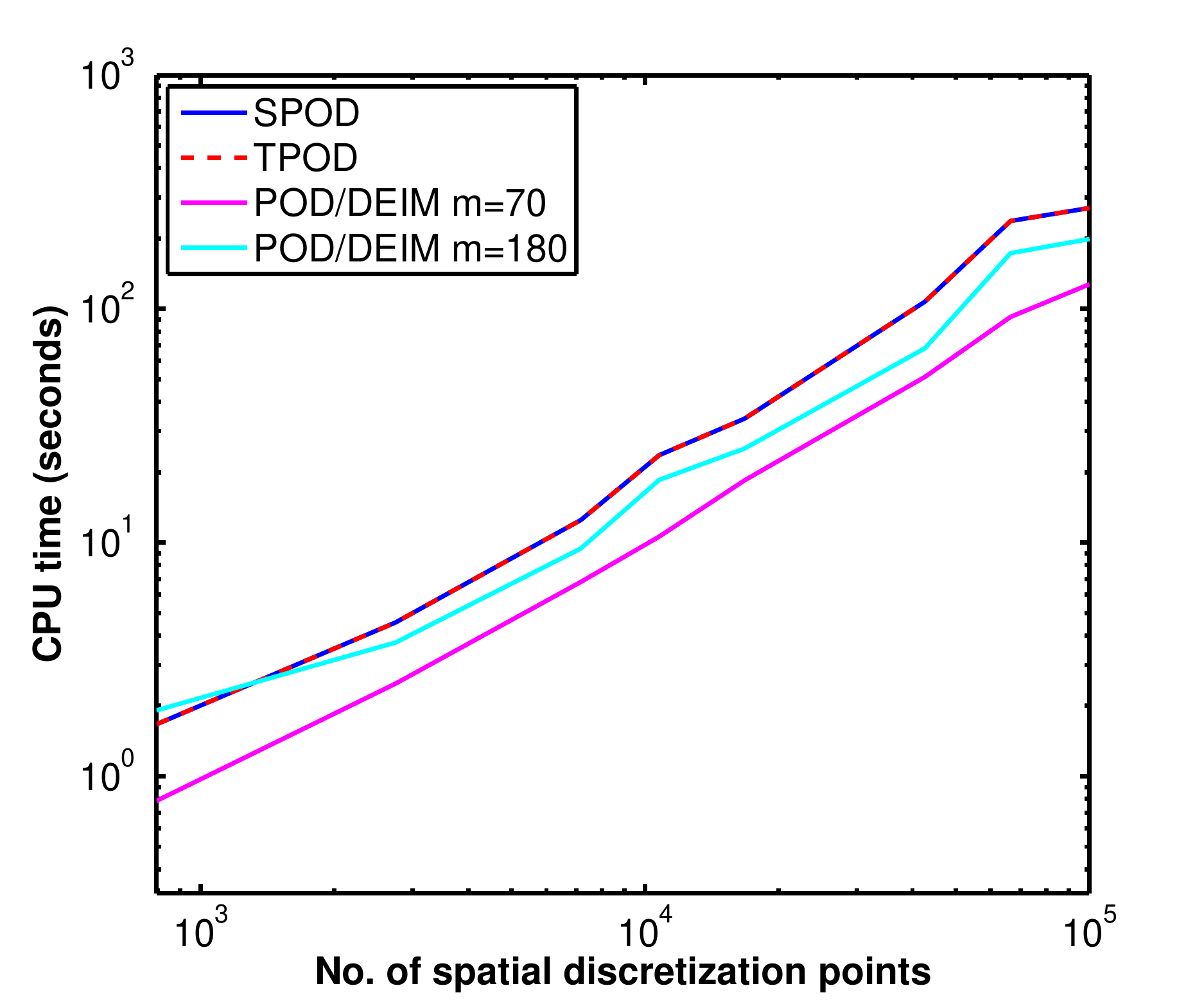}}
\caption{Cpu time vs. the number of spatial discretization points for $t_{\rm f}=3$h ; number of POD modes = 50; two different numbers of DEIM points $70$ and $180$ have been employed.\label{Fig::9}}
\end{figure}

Figure \ref{Fig::9} depicts the efficiency of tensorial POD and POD/DEIM SWE schemes as a function of spatial discretization points in the case of $t_{\rm f} = 3$h. We compare the results obtained for $8$ different mesh configurations $n=31\times 23,~61\times 45,~101\times 71,~121\times 89,~151\times 111,~241\times 177,~301\times 221,~376 \times 276$. CPU time performances of the off-line and on-line stages of the ROMs SWE schemes are compared since reduced order optimization algorithms include both phases.

For the on-line stage, once the number of spatial discretization points is larger than $151\times111$ tensorial POD scheme is $10 \times$ faster than the standard POD scheme. The performances of POD/DEIM depends on the number of DEIM points and the numerical results displays a $10 \times$ time reduction of the CPU costs in comparison with the standard POD outcome when $n \geq 61\times 45$ and $n \geq 101\times 71$ for $m=70$ and $m=180$ respectively.

The new algorithm introduced in Section $3$ relying on DEIM interpolation points delivers fast tensorial calculations required for computing the reduced Jacobian in the on-line stage and thus allowing POD/DEIM SWE scheme to have the fastest off-line stage (Figure 9b). This gives a good advantage of ROM optimization based on Discrete Empirical Interpolation Method supposing that quality approximations of nonlinear terms and reduced Jacobians are delivered since during optimization input data are different than the ones used to generate the DEIM interpolation points. DEIM was first employed by \citet{Baumann2013} to solve a reduced 4D-Var data assimilation problem and good results were obtained for a 1D Burgers  model. Extensions to 2D models are still not available in the literature.

%

\section{Conclusions}\label{sec:Conclusions}

It is well known that in standard POD the cost of evaluating nonlinear terms during the on-line stage depends on the full space dimension,
and this constitutes a major  efficiency bottleneck. The present manuscript applies tensorial calculus techniques which allows fast computations of standard POD reduced quadratic nonlinearities. We show that tensorial POD
can be applied to all type of polynomial nonlinearities and the resulting nonlinear terms have a complexity of $\mathcal{O}(k^{p+1})$ operations, where $k$ is the dimension of POD subspace and $p$ is the polynomial degree.
Consequently, this approach eliminates the dependency on the full space dimension, while yielding the same reduced solution accuracy as standard POD. Despite being independent of number of mesh points, tensorial POD is efficient only for quadratic nonlinear terms since for higher nonlinearities standard POD proves to be less time consuming once the POD basis dimension $k$ is increased.

The efficiency of tensorial POD is compared against that of standard POD and of POD/DEIM. We theoretically analyze the number of floating-point operations required as a function of polynomial degree $p$, the number of degrees of freedom of the high-fidelity model $n$, of POD modes $k$, and of DEIM interpolation points $m$. For quadratic nonlinearities and $k$ between $10$--$50$ modes, the tensorial POD needs $10$--$40$ times fewer operations than the standard POD approach and $10$--$20$ times more operations than the POD/DEIM with $m=100$.
But these performances are not translated into the same CPU time rates for solving the ROMs solutions since other more time consuming calculations are needed.

Numerical experiments are carried out using a two dimensional ADI SWE finite difference model.
Reduced order models were developed using each of the three ROM methods and Galerkin projection. The spectral analysis of snapshots matrices reveals that local versions of ROMs lead to more accurate results. Consequently, we focus on three hours time integration windows. The tensorial POD SWE model becomes considerably faster than standard POD when the dimension of the full model increases. For example, for $100,000$ spatial points the tensorial POD SWE model yields the same solutions accuracy as standard POD but is $76$ times faster. Numerical experiments of POD/DEIM SWE scheme revealed a considerable reduction of the computational complexity. For a number of $70$ DEIM points, POD/DEIM SWE model is $450$ times faster than standard POD, but only $6$ times faster than tensorial POD.

For models depending only on quadratic nonlinearities, the tensorial POD represents a solid alternative to POD/DEIM where the implementation effort is considerably larger. However, for cubic and higher order nonlinearities, tensorial POD loses its ability to deliver fast calculations and the POD/DEIM approach should be employed.

We also propose a new DEIM-based algorithm that allows fast computations of the tensors needed by reduced Jacobians calculations in the on-line stage. The resulting off-line POD/DEIM stage is the fastest among the ones considered here, even if additional SVD decompositions and low-rank terms are computed. This is an important advantage in optimization problems based on POD/DEIM surrogates where the reduced order bases need to be updated multiple times.

On-going work by the authors focuses on reduced order constrained optimization. The current research represents an important step toward developing tensorial POD and POD/DEIM four dimensional variational data assimilation systems, which are not available in the literature for complex models.

{\centering
\section*{Acknowledgments}
}
The work of Dr. R\u azvan Stefanescu and Prof. Adrian Sandu was supported by the NSF CCF--1218454, AFOSR FA9550--12--1--0293--DEF, AFOSR 12-2640-06, and by the Computational Science Laboratory at Virginia Tech. Prof. I.M. Navon acknowledges the support of NSF grant ATM-0931198. R\u azvan \c Stef\u anescu would like to thank Dr. Bernd R. Noack for his valuable suggestions on the current research topic that partially inspired the present manuscript.

\newpage
\bibliographystyle{plainnat}
\bibliography{Software,ROM_state_of_the_art,POD_bib,CDS_E_proposal,sandu,comprehensive_bibliography1,Razvan_bib,Razvan_bib_ROM_IP,NSF_KB,data_assim_weak-fdvar,reduced_models}

\begin{thebibliography}{88}
\providecommand{\natexlab}[1]{#1}
\providecommand{\url}[1]{\texttt{#1}}
\expandafter\ifx\csname urlstyle\endcsname\relax
  \providecommand{\doi}[1]{doi: #1}\else
  \providecommand{\doi}{doi: \begingroup \urlstyle{rm}\Url}\fi

\bibitem[Aanonsen(2009)]{Aanonsen2009}
T.~O. Aanonsen.
\newblock \emph{Empirical interpolation with application to reduced basis
  approximations}.
\newblock PhD thesis, Norwegian University of Science and Technology, 2009.

\bibitem[Amsallem et~al.(2011)Amsallem, Cortial, Carlberg, and
  Farhat]{Amsallem_2011}
D.~Amsallem, J.~Cortial, K.~Carlberg, and C.~Farhat.
\newblock A method for interpolating on manifolds structural dynamics
  reduced-order models.
\newblock \emph{International Journal for Numerical Methods in Engineering},
  80\penalty0 (9):\penalty0 1241--1257, 2011.

\bibitem[Antoulas(2005)]{antoulas2005als}
A.C. Antoulas.
\newblock \emph{{Approximation of Large-Scale Dynamical Systems}}.
\newblock Advances in Design and Control. SIAM, Philadelphia, 2005.

\bibitem[Antoulas(2009)]{Antoulas2009}
A.C. Antoulas.
\newblock {Approximation of large-scale dynamical systems}.
\newblock \emph{Society for Industrial and Applied Mathematics}, 6:\penalty0
  376--377, 2009.

\bibitem[Astrid et~al.(2008)Astrid, Weiland, Willcox, and Backx]{Astrid_2008}
P.~Astrid, S.~Weiland, K.~Willcox, and T.~Backx.
\newblock {Missing point estimation in models described by Proper Orthogonal
  Decomposition}.
\newblock \emph{IEEE Transactions on Automatic Control}, 53\penalty0
  (10):\penalty0 2237--2251, 2008.

\bibitem[Baker et~al.(1996)Baker, Mingori, and Goggins]{Bakeretal1996}
M.~Baker, D.~Mingori, and P.~Goggins.
\newblock {Approximate Subspace Iteration for Constructing Internally Balanced
  Reduced-Order Models of Unsteady Aerodynamic Systems}.
\newblock \emph{AIAA Meeting Papers on Disc}, pages 1070--1085, 1996.

\bibitem[Barrault et~al.(2004)Barrault, Maday, Nguyen, and Patera]{BMN2004}
M.~Barrault, Y.~Maday, N.C. Nguyen, and A.T. Patera.
\newblock An 'empirical interpolation' method: application to efficient
  reduced-basis discretization of partial differential equations.
\newblock \emph{Comptes Rendus Mathematique}, 339\penalty0 (9):\penalty0
  667--672, 2004.

\bibitem[Barrett et~al.(1994)Barrett, Berry, Chan, Demmel, Donato, Dongarra,
  Eijkhout, Pozo, Romine, and der Vorst]{Barrett94}
R.~Barrett, M.~Berry, T.~F. Chan, J.~Demmel, J.~Donato, J.~Dongarra,
  V.~Eijkhout, R.~Pozo, C.~Romine, and H.~Van der Vorst.
\newblock \emph{Templates for the Solution of Linear Systems: Building Blocks
  for Iterative Methods, 2nd Edition}.
\newblock SIAM, Philadelphia, PA, 1994.

\bibitem[Baumann(2013)]{Baumann2013}
M.M. Baumann.
\newblock {Nonlinear Model Order Reduction using POD/DEIM for Optimal Control
  of Burgersâquation}.
\newblock Master's thesis, Delft University of Technology, Netherlands, 2013.

\bibitem[Belzen and Weiland(2012)]{Belzen2012}
F.~Belzen and S.~Weiland.
\newblock {A tensor decomposition approach to data compression and
  approximation of ND systems}.
\newblock \emph{Multidimensional Systems and Signal Processing}, 23\penalty0
  (1-2):\penalty0 209--236, 2012.

\bibitem[Benner and Breiten(2012)]{MPIMD12-12}
P.~Benner and T.~Breiten.
\newblock Two-sided moment matching methods for nonlinear model reduction.
\newblock Technical Report MPIMD/12-12, Max Planck Institute Magdeburg
  Preprint, June 2012.

\bibitem[Benner and Sokolov(2006)]{Benner_Sokolov2006}
P.~Benner and V.I. Sokolov.
\newblock {Partial realization of descriptor systems}.
\newblock \emph{Systems $\&$ Control Letters}, 55\penalty0 (11):\penalty0 929
  --938, 2006.

\bibitem[Benner et~al.(2013)Benner, Gugercin, and Willcox]{MPIMD13-14}
P.~Benner, S.~Gugercin, and K.~Willcox.
\newblock A survey of model reduction methods for parametric systems.
\newblock Technical Report MPIMD/13-14, Max Planck Institute Magdeburg
  Preprint, August 2013.

\bibitem[Bo(2011)]{Nguyen_Van_Bo}
Nguyen~Van Bo.
\newblock \emph{Computational simulation of detonation waves and model
  reduction for reacting flows}.
\newblock PhD thesis, Singapore-MIT alliance, National University of Singapore,
  2011.

\bibitem[Bui-Thanh et~al.(2004)Bui-Thanh, Damodaran, and
  Willcox]{Bui-thanh04aerodynamicdata}
T.~Bui-Thanh, M.~Damodaran, and K.~Willcox.
\newblock Aerodynamic data reconstruction and inverse design using proper
  orthogonal decomposition.
\newblock \emph{AIAA Journal}, pages 1505--1516, 2004.

\bibitem[Bultheel and Moor(2000)]{Bultheel_Moor2000}
A.~Bultheel and B.~De Moor.
\newblock {Rational approximation in linear systems and control}.
\newblock \emph{Journal of Computational and Applied Mathematics},
  121:\penalty0 355--378, 2000.

\bibitem[Carlberg and Farhat(2011)]{Carlberg_2011}
K.~Carlberg and C.~Farhat.
\newblock A low-cost, goal-oriented 'compact proper orthogonal decomposition'
  basis for model reduction of static systems.
\newblock \emph{International Journal for Numerical Methods in Engineering},
  86\penalty0 (3):\penalty0 381--402, 2011.

\bibitem[Carlberg et~al.(2011)Carlberg, Bou-Mosleh, and Farhat]{Carlberg2_2011}
K.~Carlberg, C.~Bou-Mosleh, and C.~Farhat.
\newblock Efficient non-linear model reduction via a least-squares
  {P}etrov-–{G}alerkin projection and compressive tensor approximations.
\newblock \emph{International Journal for Numerical Methods in Engineering},
  86\penalty0 (2):\penalty0 155--181, 2011.

\bibitem[Carlberg et~al.(2012)Carlberg, Tuminaro, and Boggsz]{Carlberg_2012}
K.~Carlberg, R.~Tuminaro, and P.~Boggsz.
\newblock { Efficient structure-preserving model reduction for nonlinear
  mechanical systems with application to structural dynamics}.
\newblock preprint, Sandia National Laboratories, Livermore, CA 94551, USA,
  2012.

\bibitem[Chaturantabut(2008)]{Cha2008}
S.~Chaturantabut.
\newblock Dimension {R}eduction for {U}nsteady {N}onlinear {P}artial
  {D}ifferential {E}quations via {E}mpirical {I}nterpolation {M}ethods.
\newblock Technical Report TR09-38,CAAM, Rice University, 2008.

\bibitem[Chaturantabut and Sorensen(2011)]{ChaSor2011}
S.~Chaturantabut and D~.C. Sorensen.
\newblock Application of {POD} and {DEIM} on dimension reduction of non-linear
  miscible viscous fingering in porous media.
\newblock \emph{Mathematical and Computer Modelling of Dynamical Systems},
  17\penalty0 (4):\penalty0 337--353, 2011.

\bibitem[Chaturantabut and Sorensen(2010)]{ChaSor2010}
S.~Chaturantabut and D.C. Sorensen.
\newblock Nonlinear model reduction via discrete empirical interpolation.
\newblock \emph{SIAM Journal on Scientific Computing}, 32\penalty0
  (5):\penalty0 2737--2764, 2010.

\bibitem[Chaturantabut and Sorensen(2012)]{ChaSor2012}
S.~Chaturantabut and D.C. Sorensen.
\newblock A state space error estimate for {POD-DEIM} nonlinear model
  reduction.
\newblock \emph{SIAM Journal on Numerical Analysis}, 50\penalty0 (1):\penalty0
  46--63, 2012.

\bibitem[Dihlmann and Haasdonk(2013)]{Dihlmann_2013}
M.~Dihlmann and B.~Haasdonk.
\newblock Certified {PDE}-constrained parameter optimization using reduced
  basis surrogate models for evolution problems.
\newblock \emph{Submitted to the Journal of Computational Optimization and
  Applications}, 2013.
\newblock URL
  \url{http://www.agh.ians.uni-stuttgart.de/publications/2013/DH13}.

\bibitem[Everson and Sirovich(1995)]{Everson_1995}
R.~Everson and L.~Sirovich.
\newblock Karhunen–-{L}oeve procedure for gappy data.
\newblock \emph{Journal of the Optical Society of America A}, 12:\penalty0
  1657--64, 1995.

\bibitem[Fairweather and Navon(1980)]{FN1980}
G.~Fairweather and I.M. Navon.
\newblock A linear {ADI} method for the shallow water equations.
\newblock \emph{Journal of Computational Physics}, 37:\penalty0 1--18, 1980.

\bibitem[Feldmann and Freund(1995)]{Feldmann_Freund1995}
P.~Feldmann and R.W. Freund.
\newblock {Efficient linear circuit analysis by Pade approximation via the
  Lanczos process}.
\newblock \emph{IEEE Transactions on Computer-Aided Design of Integrated
  Circuits and Systems}, 14:\penalty0 639--649, 1995.

\bibitem[Freund(2003)]{Freund2003}
R.W. Freund.
\newblock {Model reduction methods based on Krylov subspaces}.
\newblock \emph{Acta Numerica}, 12:\penalty0 267--319, 2003.

\bibitem[Gallivan et~al.(1994)Gallivan, Grimme, and Dooren]{Gallivan94}
K.~Gallivan, E.~Grimme, and P.~Van Dooren.
\newblock Pad\'{e} approximation of large-scale dynamic systems with lanczos
  methods.
\newblock In \emph{Decision and Control, 1994., Proceedings of the 33rd IEEE
  Conference on}, volume~1, pages 443--448 vol.1, Dec 1994.

\bibitem[Gragg(1972)]{Gragg1972}
W.B. Gragg.
\newblock {The Pad\'{e} table and its relation to certain algorithms of
  numerical analysis}.
\newblock \emph{SIAM Review}, 14:\penalty0 1--62, 1972.

\bibitem[Gragg and Lindquist(1983)]{Gragg_Lindquist1983}
W.B. Gragg and A.~Lindquist.
\newblock {On the partial realization problem}.
\newblock \emph{Linear Algebra and Its Applications, Special Issue on Linear
  Systems and Control}, 50:\penalty0 277 --319, 1983.

\bibitem[Grammeltvedt(1969)]{Gram1969}
A.~Grammeltvedt.
\newblock {A survey of finite difference schemes for the primitive equations
  for a barotropic fluid}.
\newblock \emph{Monthly Weather Review}, 97\penalty0 (5):\penalty0 384--404,
  1969.

\bibitem[Grepl and Patera(2005)]{grepl2005posteriori}
M.A. Grepl and A.T. Patera.
\newblock A posteriori error bounds for reduced-basis approximations of
  parametrized parabolic partial differential equations.
\newblock \emph{ESAIM: Mathematical Modelling and Numerical Analysis},
  39\penalty0 (01):\penalty0 157--181, 2005.

\bibitem[Grepl et~al.(2007)Grepl, Maday, Nguyen, and Patera]{Grep_Mad_Nguy_Pat}
M.A. Grepl, Y.~Maday, N.C. Nguyen, and A.T. Patera.
\newblock Efficient reduced-basis treatment nofonaffine and nonlinear partial
  differential equations.
\newblock \emph{Mod\'{e}lisation Math\'{e}matique et Analyse Num\'{e}rique},
  41\penalty0 (3):\penalty0 575--605, 2007.

\bibitem[Grimme(1997)]{Grimme1997}
E.J. Grimme.
\newblock \emph{Krylov projection methods for model reduction}.
\newblock PhD thesis, Univ. Illinois, Urbana-Champaign, 1997.

\bibitem[Gudmundsson and Laub(1994)]{Gudmundsson_Laub1994}
T.~Gudmundsson and A.~Laub.
\newblock {Approximate Solution of Large Sparse Lyapunov Equations}.
\newblock \emph{IEEE Transactions on Automatic Control}, 39\penalty0
  (5):\penalty0 1110--1114, 1994.

\bibitem[Gustafsson(1971)]{Gus1971}
B.~Gustafsson.
\newblock {An alternating direction implicit method for solving the shallow
  water equations}.
\newblock \emph{Journal of Computational Physics}, 7:\penalty0 239--254, 1971.

\bibitem[Gutknecht(1994)]{Gutknecht1994}
M.H. Gutknecht.
\newblock {The Lanczos process and Pad\'{e} approximation}.
\newblock \emph{Proc. Cornelius Lanczos Intl. Centenary Conference, edited by
  J.D. Brown et al., SIAM, Philadelphia}, pages 61--75, 1994.

\bibitem[Hinze and Kunkel(2012)]{HK2012}
M.~Hinze and M.~Kunkel.
\newblock Discrete {E}mpirical {I}nterpolation in {POD} {M}odel {O}rder
  {R}eduction of {D}rift-{D}iffusion {E}quations in {E}lectrical {N}etworks.
\newblock \emph{Scientific Computing in Electrical Engineering SCEE 2010
  Mathematics in Industry}, 16\penalty0 (5):\penalty0 423--431, 2012.

\bibitem[Hochman et~al.(2011)Hochman, Bond, and White]{HBW2011}
A.~Hochman, B.N. Bond, and J.K. White.
\newblock A stabilized discrete empirical interpolation method for model
  reduction of electrical thermal and microelectromechanical systems.
\newblock \emph{Design Automation Conference (DAC), 48th ACM/EDAC/IEEE}, pages
  540--545., 2011.

\bibitem[Hodel(1991)]{Hodel1991}
A.S. Hodel.
\newblock {Least Squares Approximate Solution of the Lyapunov Equation}.
\newblock \emph{Proceedings of the 30th IEEE Conference on Decision and
  Control, IEEE Publications, Piscataway, NJ}, 1991.

\bibitem[Hotelling(1933)]{hotelling1939acs}
H.~Hotelling.
\newblock Analysis of a complex of statistical variables with principal
  components.
\newblock \emph{Journal of Educational Psychology}, 24:\penalty0 417--441,
  1933.

\bibitem[Jaimoukha and Kasenally(1994)]{Jaimoukha_Kasenally1994}
I.~Jaimoukha and E.~Kasenally.
\newblock {Krylov Subspace Methods for Solving Large Lyapunov Equations}.
\newblock \emph{SIAM Journal of Numerical Analysis}, 31\penalty0 (1):\penalty0
  227--251, 1994.

\bibitem[Kaiser et~al.(2013)Kaiser, Noack, Cordier, Spohn, Segond, Abel,
  Daviller, and Niven]{Kaiseretal2013}
E.~Kaiser, Bernd~R. Noack, L.~Cordier, A.~Spohn, M.~Segond, M.~Abel,
  G.~Daviller, and R.K. Niven.
\newblock {Cluster-based reduced-order modelling of a mixing layer}.
\newblock Technical Report arXiv:1309.0524 [physics.flu-dyn], Cornell
  University, September 2013.

\bibitem[Karhunen(1946)]{karhunen1946zss}
K.~Karhunen.
\newblock Zur spektraltheorie stochastischer prozesse.
\newblock \emph{Annales Academiae Scientarum Fennicae}, 37, 1946.

\bibitem[Kellems et~al.(2010)Kellems, Chaturantabut, Sorensen, and
  Cox]{ChaSor2_2010}
A.~R. Kellems, S.~Chaturantabut, D.~C. Sorensen, and S.~J. Cox.
\newblock Morphologically accurate reduced order modeling of spiking neurons.
\newblock \emph{Journal of Computational Neuroscience}, 28:\penalty0 477--494,
  2010.

\bibitem[Kelley(1995)]{Kelley95}
C.~T. Kelley.
\newblock \emph{Iterative Methods for Linear and Nonlinear Equations}.
\newblock Number~16 in Frontiers in Applied Mathematics. SIAM, 1995.

\bibitem[Kunisch and Volkwein(1999)]{Kunisch_Volkwein_1999}
K.~Kunisch and S.~Volkwein.
\newblock Control of the {B}urgers {E}quation by a {R}educed-{O}rder {A}pproach
  {U}sing {P}roper {O}rthogonal {D}ecomposition.
\newblock \emph{Journal of Optimization Theory and Applications}, 102\penalty0
  (2):\penalty0 345--371, 1999.

\bibitem[Kunisch and Volkwein(2002)]{Kunisch_Volkwein_POD2002}
K.~Kunisch and S.~Volkwein.
\newblock Galerkin {P}roper {O}rthogonal {D}ecomposition {M}ethods for a
  {G}eneral {E}quation in {F}luid {D}ynamics.
\newblock \emph{SIAM Journal on Numerical Analysis}, 40\penalty0 (2):\penalty0
  492--515, 2002.

\bibitem[Kunisch et~al.(2004)Kunisch, Volkwein, and Xie]{Kunisch_Volkwein_2004}
K.~Kunisch, S.~Volkwein, and L.~Xie.
\newblock {HJB}-{POD}-{B}ased {F}eedback {D}esign for the {O}ptimal {C}ontrol
  of {E}volution {P}roblems.
\newblock \emph{SIAM J. Appl. Dyn. Syst}, 3\penalty0 (4):\penalty0 701--722,
  2004.

\bibitem[Lass and Volkwein(2012)]{LW2012}
O.~Lass and S.~Volkwein.
\newblock {POD} {G}alerkin schemes for nonlinear elliptic-parabolic systems.
\newblock \emph{Konstanzer Schriften in Mathematik}, 301:\penalty0 1430--3558,
  2012.

\bibitem[Lloyd(1957)]{Lloyd_1957}
S.~Lloyd.
\newblock Least squares quantization in {PCM}.
\newblock \emph{IEEE Trans. Inform. Theory}, 28:\penalty0 129--137, 1957.

\bibitem[Lo\`eve(1955)]{loeve1955pt}
M.M. Lo\`eve.
\newblock \emph{Probability Theory}.
\newblock Van Nostrand, Princeton, NJ, 1955.

\bibitem[Lorenz(1956)]{lorenz1956eof}
E.N. Lorenz.
\newblock {Empirical Orthogonal Functions and Statistical Weather Prediction}.
\newblock Technical report, Massachusetts Institute of Technology, Dept. of
  Meteorology, 1956.

\bibitem[MacQueen(1967)]{MacQueen_1967}
J.~MacQueen.
\newblock Some methods for classification and analysis of multivariate
  observations.
\newblock \emph{Proceedings of the Fifth Berkeley Symposium on Mathematical
  Statistics and Probability}, 1:\penalty0 281--297, 1967.

\bibitem[Maday et~al.(2009)Maday, Nguyen, Patera, and Pau]{Mad_Nguy_Pat_Pau}
Y.~Maday, N.C. Nguyen, A.T. Patera, and G.S.H. Pau.
\newblock A {G}eneral {M}ultipurpose {I}nterpolation {P}rocedure: the {M}agic
  {P}oints.
\newblock \emph{Communications on Pure and Applied Analysis}, 8\penalty0
  (1):\penalty0 383--404, 2009.

\bibitem[Moore(1981)]{Moore1981}
B.C. Moore.
\newblock {Principal component analysis in linear systems: Controllability,
  observability, and model reduction}.
\newblock \emph{IEEE Transactions on Automatic Control}, 26\penalty0
  (1):\penalty0 17--32, 1981.

\bibitem[Mullis and Roberts(1976)]{Mullis_Roberts1976}
C.T. Mullis and R.A. Roberts.
\newblock {Synthesis of Minimum Roundoff Noise Fixed Point Digital Filters}.
\newblock \emph{IEEE Transactions on Circuits and Systems}, CAS-23:\penalty0
  551--562, 1976.

\bibitem[Navon and Villiers(1986)]{NVG1986}
I.~M. Navon and R.~De Villiers.
\newblock Gustaf: A {Q}uasi-{N}ewton nonlinear {ADI} fortran iv program for
  solving the shallow-water equations with augmented lagrangians.
\newblock \emph{Computers and Geosciences}, 12\penalty0 (2):\penalty0 151--173,
  1986.

\bibitem[Nguyen et~al.(2008)Nguyen, Patera, and Peraire]{NPP2008}
N.C. Nguyen, A.T. Patera, and J.~Peraire.
\newblock {A 'best points' interpolation method for efficient approximation of
  parametrized function}.
\newblock \emph{International Journal for Numerical Methods in Engineering},
  73:\penalty0 521--543, 2008.

\bibitem[Noack et~al.(2003)Noack, Afanasiev, Morzynski, Tadmor, and
  Thiele]{NAMTT03}
B.R. Noack, K.~Afanasiev, M.~Morzynski, G.~Tadmor, and F.~Thiele.
\newblock A hierarchy of low-dimensional models for the transient and
  post-transient cylinder wake.
\newblock \emph{Journal of Fluid Mechanics}, 497:\penalty0 335--363, 2003.
\newblock ISSN 0022-1120.

\bibitem[Noack et~al.(2010)Noack, Schlegel, Morzynski, and Tadmor]{Noack2010}
B.R. Noack, M.~Schlegel, M.~Morzynski, and G.~Tadmor.
\newblock System reduction strategy for galerkin models of fluid flows.
\newblock \emph{International Journal for Numerical Methods in Fluids},
  63\penalty0 (2):\penalty0 231--248, 2010.

\bibitem[Patera and Rozza(2007)]{patera2007reduced}
A.T. Patera and G.~Rozza.
\newblock Reduced basis approximation and a posteriori error estimation for
  parametrized partial differential equations, 2007.

\bibitem[Peherstorfer et~al.(2013)Peherstorfer, Butnaru, Willcox, and
  Bungartz]{Peherstorfer_2013}
B.~Peherstorfer, D.~Butnaru, K.~Willcox, and H.J. Bungartz.
\newblock { Localized Discrete Empirical Interpolation Method}.
\newblock MIT Aerospace Computational Design Laboratory Technical Report
  TR-13-1, 2013.

\bibitem[Rap\'{u}n and Vega(2010)]{Rapun_2010}
M.L. Rap\'{u}n and J.M. Vega.
\newblock Reduced order models based on local {POD} plus {G}alerkin projection.
\newblock \emph{Journal of Computational Physics}, 229\penalty0 (8):\penalty0
  3046--3063, 2010.

\bibitem[Rewie\'{n}ski and White(2001)]{Rewienski2001}
M.~Rewie\'{n}ski and J.~White.
\newblock {A Trajectory Piecewise-linear Approach to Model Order Reduction and
  Fast Simulation of Nonlinear Circuits and Micromachined Devices}.
\newblock In \emph{Proceedings of the 2001 IEEE/ACM International Conference on
  Computer-aided Design}, ICCAD '01, pages 252--257, Piscataway, NJ, USA, 2001.
  IEEE Press.

\bibitem[Rowley et~al.(2004)Rowley, Colonius, and Murray]{Rowley2004}
C.~W. Rowley, T.~Colonius, and R.~M. Murray.
\newblock Model reduction for compressible flows using {POD} and {G}alerkin
  projection.
\newblock \emph{Physica D. Nonlinear Phenomena}, 189\penalty0 (1--2):\penalty0
  115--129, 2004.

\bibitem[Rowley(2005)]{Rowley2005}
C.W. Rowley.
\newblock Model {R}eduction for {F}luids, using {B}alanced {P}roper
  {O}rthogonal {D}ecomposition.
\newblock \emph{International Journal of Bifurcation and Chaos (IJBC)},
  15\penalty0 (3):\penalty0 997--1013, 2005.

\bibitem[Rowley et~al.(2009)Rowley, Mezic, Bagheri, P.Schlatter, and
  Henningson]{Rowley2009}
C.W. Rowley, I.~Mezic, S.~Bagheri, P.Schlatter, and D.S. Henningson.
\newblock Spectral analysis of nonlinear flows.
\newblock \emph{Journal of Fluid Mechanics}, 641:\penalty0 115--127, 2009.

\bibitem[Rozza et~al.(2008)Rozza, Huynh, and Patera]{rozza2008reduced}
G.~Rozza, D.B.P. Huynh, and A.T. Patera.
\newblock Reduced basis approximation and a posteriori error estimation for
  affinely parametrized elliptic coercive partial differential equations.
\newblock \emph{Archives of Computational Methods in Engineering}, 15\penalty0
  (3):\penalty0 229--275, 2008.

\bibitem[Saad(1994)]{Saad1994}
Y.~Saad.
\newblock {Sparsekit: a basic tool kit for sparse matrix computations}.
\newblock {Technical Report, Computer Science Department, University of
  Minnesota}, 1994.

\bibitem[Saad(2003)]{Saad2003}
Y.~Saad.
\newblock \emph{Iterative Methods for Sparse Linear Systems}.
\newblock Society for Industrial and Applied Mathematics, Philadelphia, PA,
  USA, 2nd edition, 2003.

\bibitem[San and Iliescu(2013)]{San_Iliescu2013}
O.~San and T.~Iliescu.
\newblock {Proper orthogonal decomposition closure models for fluid flows:
  Burgers equation}.
\newblock Technical Report arXiv:1308.3276 [physics.flu-dyn], Cornell
  University, August 2013.

\bibitem[Schmid(2010)]{Schmid2010}
P.J. Schmid.
\newblock Dynamic mode decomposition of numerical and experimental data.
\newblock \emph{Journal of Fluid Mechanics}, 656:\penalty0 5--28, 2010.

\bibitem[Sirovich(1987{\natexlab{a}})]{Sir87a}
L.~Sirovich.
\newblock Turbulence and the dynamics of coherent structures. {I}. {C}oherent
  structures.
\newblock \emph{Quarterly of Applied Mathematics}, 45\penalty0 (3):\penalty0
  561--571, 1987{\natexlab{a}}.
\newblock ISSN 0033-569X.

\bibitem[Sirovich(1987{\natexlab{b}})]{Sir87b}
L.~Sirovich.
\newblock Turbulence and the dynamics of coherent structures. {II}.
  {S}ymmetries and transformations.
\newblock \emph{Quarterly of Applied Mathematics}, 45\penalty0 (3):\penalty0
  573--582, 1987{\natexlab{b}}.
\newblock ISSN 0033-569X.

\bibitem[Sirovich(1987{\natexlab{c}})]{Sir87c}
L.~Sirovich.
\newblock Turbulence and the dynamics of coherent structures. {III}. {D}ynamics
  and scaling.
\newblock \emph{Quarterly of Applied Mathematics}, 45\penalty0 (3):\penalty0
  583--590, 1987{\natexlab{c}}.
\newblock ISSN 0033-569X.

\bibitem[Sorensen and Antoulas(2002)]{Sorensen_Antoulas2002}
D.C. Sorensen and A.C. Antoulas.
\newblock {The Sylvester equation and approximate balanced reduction}.
\newblock \emph{Linear Algebra and its Applications}, 351-352\penalty0
  (0):\penalty0 671--700, 2002.

\bibitem[Stefanescu and Navon(2013)]{Stefanescu2013}
R.~Stefanescu and I.M. Navon.
\newblock {POD/DEIM} {N}onlinear model order reduction of an {ADI} implicit
  shallow water equations model.
\newblock \emph{Journal of Computational Physics}, 237:\penalty0 95--114, 2013.

\bibitem[Steinhaus(1956)]{Steinhaus_1956}
H.~Steinhaus.
\newblock Sur la division des corps mat\'{e}riels en parties.
\newblock \emph{Bulletin of the Polish Academy of Sciences}, 4\penalty0
  (12):\penalty0 801--804, 1956.

\bibitem[Suwartadi(2012)]{EkaSuwartadi2012}
E.~Suwartadi.
\newblock \emph{Gradient-based {M}ethods for {P}roduction {O}ptimization of
  {O}il {R}eservoirs}.
\newblock PhD thesis, Mathematics and Electrical Engineering, Department of
  Engineering Cybernetics,Norwegian University of Science and Technology, 2012.

\bibitem[Tissot et~al.(2013)Tissot, Cordier, Benard, and Noack]{Tissot2013}
G.~Tissot, L.~Cordier, N.~Benard, and B.R. Noack.
\newblock {Dynamic mode decomposition of PIV measurements for cylinder wake
  flow in turbulent regime}.
\newblock In \emph{Proceedings of the 8th International Symposium On Turbulent
  and Shear Flow Phenomena}, TSFP-8, 2013.

\bibitem[Tonn(2012)]{Tonn_2011}
T.~Tonn.
\newblock { Reduced-Basis Method (RBM) for Non-Affine Elliptic Parametrized
  {PDE}s}.
\newblock (PhD), Ulm University, 2012.

\bibitem[{Van Dooren}(1995)]{VanDooren_1995}
P.~{Van Dooren}.
\newblock {The Lanczos algorithm and Pad\'{e} approximations}.
\newblock In \emph{Short Course, Benelux Meeting on Systems and Control}, 1995.

\bibitem[Willcox and Peraire(2002)]{Willcox02balancedmodel}
K.~Willcox and J.~Peraire.
\newblock Balanced model reduction via the {Proper Orthogonal Decomposition}.
\newblock \emph{AIAA Journal}, pages 2323--2330, 2002.

\bibitem[Wirtz et~al.(2012)Wirtz, Sorensen, and Haasdonk]{Wirtz2012}
D.~Wirtz, D.C. Sorensen, and B.~Haasdonk.
\newblock {A-posteriori error estimation for DEIM reduced nonlinear dynamical
  systems}.
\newblock SRC SimTech Preprint Series, 2012.

\bibitem[Xiao et~al.(2014)Xiao, Fang, Buchan, Pain, Navon, Du, and
  Hu]{Xiao2014}
D.~Xiao, F.~Fang, A.G. Buchan, C.C. Pain, I.M. Navon, J.~Du, and G.~Hu.
\newblock {Non-linear model reduction for the Navier-Stokes equations using
  residual DEIM method}.
\newblock \emph{Journal of Computational Physics}, 263:\penalty0 1--18, 2014.
\newblock ISSN 0021-9991.

\bibitem[Zhou(2012)]{Zhou_2012}
Y.B. Zhou.
\newblock {Model reduction for nonlinear dynamical systems with parametric
  uncertainties}.
\newblock (M.S), Massachusetts Institute of Technology, Dept. of Aeronautics
  and Astronautics, 2012.

\end{thebibliography}

\end{document}